\documentclass[10pt]{article}
\usepackage{mathrsfs}

\usepackage{epsfig,epsf}

\usepackage{graphics}
\usepackage{amscd}
\usepackage{graphics,graphicx}
\usepackage{latexsym}
\usepackage{amsthm}
\usepackage{amsmath}
\usepackage{amssymb}
\usepackage{mathrsfs}
\usepackage{bbold}

\newcommand{\grad}{\nabla }

\graphicspath{{Figures/}}
\newcommand{\N}[1]{\mbox{I\hspace{-.15em}N}^{#1}}

\newcommand{\R}[1]{\mbox{I\hspace{-.15em}R}^{#1}}

\newcommand{\PH}[1]{{\phi^{*}}}
\newcommand{\PS}[1]{{\psi^{*}}}
\newcommand{\PTH}[1]{\tilde{\phi}}
\newcommand{\PTS}[1]{\tilde{\psi}}
\newcommand{\uT}[1]{\tilde{u}}
\newcommand{\vT}[1]{\tilde{v}}
\newcommand{\wT}[1]{\tilde{w}}
\newcommand{\ucT}[1]{\tilde{U}}
\newcommand{\vcT}[1]{\tilde{V}}
\newcommand{\h}[1]{{\bf h}}
\newcommand{\U}[1]{{\bf u}}
\newcommand{\UL}[1]{\textcolor{green}{\textcolor{green}{{\bf u}_{L}}}}
\newcommand{\UT}[1]{\tilde{\bf u}}
\newcommand{\hL}[1]{\textcolor{green}{\textcolor {green}{h_{L}}}}
\newcommand{\V}[1]{{\bf v}}
\newcommand{\W}[1]{{\bf w}}
\newcommand{\f}[1]{{\bf f}}
\newcommand{\g}[1]{{\bf g}}
\newcommand{\Ze}[1]{\bf{0}}
\newcommand{\Uk}[1]{{\bf u}_{k}}
\newcommand{\Vk}[1]{{\bf v}_{k}}
\newcommand{\fk}[1]{{\bf f}_{k}}
\newcommand{\gk}[1]{{\bf g}_{k}}
\newcommand{\aUk}[1]{\mid {\bf u}_{k}\mid ^2}
\newcommand{\ahk}[1]{\mid h_{k}\mid^2}
\newcommand{\UUk}[1]{X_{k}}
\newcommand{\UU}[1]{{\bf U}}
\newcommand{\RR}[1]{{\bf R}}
\newcommand{\aUUk}[1]{\mid X_{k}\mid ^2}
\newcommand{\VVk}[1]{X_{k}}
\newcommand{\ffk}[1]{{\bf F}_{k}}
\newcommand{\ffe}[1]{{\bf F}_{\epsilon}}
\newcommand{\hk}[1]{h_{k}}
\newcommand{\UUe}[1]{{\bf U}_{\epsilon}}
\newcommand{\aUUe}[1]{\mid {\bf U}_{\epsilon}\mid_{2}^2}
\newcommand{\We}[1]{{\bf w}_{\epsilon}}
\newcommand{\fe}[1]{{\bf f}_{\epsilon}}
\newcommand{\He}[1]{H_{\epsilon}}
\newcommand{\aHe}[1]{\mid H_{\epsilon}\mid ^2}
\newcommand{\RRe}[1]{{\bf R}_{\epsilon}}
\newcommand{\SSe}[1]{S_{\epsilon}}
\newcommand{\aRRe}[1]{\mid {\bf R}_{\epsilon}\mid^2}
\newcommand{\aSSe}[1]{\mid S_{\epsilon}\mid ^2}
\newcommand{\pr}[1]{\parallel}
\newcommand{\fld}[1]{\longrightarrow}
\newcommand{\p}[1]{\partial}
\newcommand{\fb}[1]{\rightharpoonup}
\newcommand{\ft}[1]{\longrightarrow}
\newcommand{\ITT}[1]{\int_0^T}
\newcommand{\ev}[1]{L^2(0,T,V)}
\newcommand{\eh}[1]{L^{\infty}(0,T,H)}
\newcommand{\evh}[1]{L^2(0,T,V) \cap L^{\infty}(0,T,H)}
\newcommand{\demi}[1]{\frac{1}{2}}
\newcommand{\ueps}[1]{\frac{1}{\epsilon}}
\newcommand{\Je}[1]{J_{\epsilon}}
\newcommand{\Xe}[1]{X_{\epsilon}}
\newcommand{\aXe}[1]{\mid X_{\epsilon}\mid ^2}
\newcommand{\Ye}[1]{Y_{\epsilon}}
\newcommand{\aYe}[1]{\mid Y_{\epsilon}\mid ^2}
\newcommand{\TRRe}[1]{\widetilde{{\bf R}_{\epsilon}}}
\newcommand{\TSSe}[1]{\widetilde{S_{\epsilon}}}
\newcommand{\aTRRe}[1]{\mid \widetilde{{\bf R}_{\epsilon}}\mid^2}
\newcommand{\aTSSe}[1]{\mid \widetilde{S_{\epsilon}}\mid ^2}

\newcommand{\Ug}[1]{{\bf U}_{g}}
\newcommand{\n}[1]{{\bf n}}
\newcommand{\Hg}[1]{H_{g}}
\newcommand{\Uw}[1]{{\bf U}_{w}}
\newcommand{\Hw}[1]{H_{w}}
\newcommand{\VV}[1]{\vec{\bf v}}

\newtheorem{theo}{Theorem}
\newtheorem{pps}{Proposition}

\newtheorem{req}{Remark}

\begin{document}
\date{}
%\markboth{Aziz Belmiloudi}{Instructions for Typing Manuscripts (Paper's Title)}

%%%%%%%%%%\catchline{}{}{}{}{}
%
%%%%%%%%%%%%%%%%%%%%%%%%%%%%%%%%%%%%%%%%%%%%%%%%%%%%%%%%%%%%%%%%%%%%

\title{Optimal Control for Chemotaxis Systems and Adjoint-Based Optimization with Multiple-Relaxation-Time Lattice Boltzmann Models}

\author{Aziz Belmiloudi\footnote{Mathematics research institute of Rennes (IRMAR), European University of Brittany (UEB),  20 av des Buttes de Co\"esmes, 
CS 14315, 35043 Rennes C\'{e}dex, France (Aziz.Belmiloudi@math.cnrs.fr)}}

\maketitle

\begin{abstract}
This paper is devoted to continuous and discrete adjoint-based optimization approaches for optimal control problems governed by an important class of Nonlinear Coupled Anisotropic Convection-Diffusion Chemotaxis-type System (NCACDCS). This study is motivated by the fact that the considered complex systems (with complex geometries) appear in diverse biochemical, biological and biosocial criminology problems. To solve numerically the corresponding nonlinear optimization problems, the primal problem  NCACDCS is discretised by a coupled Lattice Boltzmann Method with a general Multiple-Relaxation-Time  collision operators (MRT) while for the adjoint problem, an Adjoint Multiple-Relaxation-Time lattice Boltzmann model (AMRT) is proposed and investigated.

First, the optimal control problems are formulated and first-order necessary optimality conditions are established by using sensitivity and adjoint calculus. 
The resulting problems are discretised by the coupled MRT and AMRT models and solved via gradient descent methods. First of all, an efficient and stable modified MRT model for NCACDCS is developed, and through the Chapman-Enskog analysis we show that NCACDCS can be correctly recovered from the proposed MRT model. For the adjoint problem, the discretisation strategy is based on AMRT model, which is found to be as simple as MRT model with also highly-efficient parallel nature. The derivation of AMRT model and the discrete cost functional gradient are derived mathematically in detail using the developed MRT model. The obtained method is reliable, efficient, practical to implement and can be easily incorporated into any existing MRT code.
\end{abstract}
{\bf keywords} {Optimal control, nonlinear coupled anisotropic convection-diffusion chemotaxis-type system, self and cross diffusion tensors, adjoint system,  adjoint-based multiple-relaxation-time lattice Boltzmann method, multiscale Chapman-Enskog expansion, optimization.}
\section{Introduction and mathematical setting of the problem}
The mathematical problem considered in this paper derives from the modeling of directed motion of biological organisms in response to changes of a number of chemical and physical cues in their environment. 
During the last few years, the effect of chemical gradients on organism migration has been widely studied, and complex systems have been developed to analyze organism behavior under simultaneous chemical, electrical, and mechanical stimuli.

In contrast to random diffusion without orientation, chemotaxis describes directed movements and migrations of bacteria, biological cells and organisms under the effect of a chemical stimulus (produced by a substance therein inhomogeneously distributed), in which organisms move either toward or away from the stimulus (i.e., attractive or repulsive movements); it is represented by the gradients of concentration of chemical substances (which have been characterized at molecular level). It is also by far the best understood mode of directed organism movements. 
The chemical signals are either secreted by the organisms themselves (autogenous) or come from foreign sources (xenogenous), and lead to aggregation of organisms and to formation of complex patterns. 

Chemotaxis plays a crucial fundamental role, in physiology and pathological conditions, in many biological processes (see e.g., \cite{Eal}) such as tissue formation, regeneration or immune functions, and pathological processes like cancer metastasis, as well as in many social activities of micro-organisms, e.g. virulence, social motility, fruiting body formation and biofilm formation.

Recently, numerous and various problems concerning chemotaxis have been intensively studied either from a theoretical or  from a numerical point of view (the literature is extensive in this field). Many mathematical models have been developed to describe the chemotactic organism. The most generalized and widely-used model of chemotaxis is in connection with Keller-Segel models \cite{KS1,KS2} which are a coupled nonlinear/degenerate parabolic equations. Well-known examples are the dynamics and aggregation of Escherichia coli colonies and of the slime mould amoebae Dictyostelium discoideum. Also, endothelial cells react to vascular endothelial growth factor to form the primary blood vessels through aggregation (phenomenon of vasculogenesis) \cite{ABP, Sal}. Keller-Segel model has been successfully used by many scientists for the analysis of various biological phenomena. Others, after evaluations of this classical (minimal) chemotaxis model in specific situations, have concluded that many of the hypotheses (foundational to the model) are not sufficient to describe fundamental phenomena (in chemotaxis) in a satisfactory way. Then these latter ones have modified and generalized the model to adequate systems by incorporating additional informations (such as describing aggregation and blow up of attractive chemotaxis and collapse of repulsive chemotaxis, tumor angiogenesis and invasion, vasculogenesis, organism pattern formation, morphogenesis, attraction and repulsion between mutually interacting species in ecology and epidemiology, criminal behavior, etc) and obtain more complicated and more realistic models, see e.g., \cite{ABP, AC,  AM, ASV,  Ba, BL1, BL2, BBI, EJDW, Eal, GC, HP1, Ho2, Ho3, Ho4, I1, LCEM, LCA, MM, Mu, Mu1, OS,  PS1, PH2, Pe, RST, Sal, SOPT, VM, Wr1, Wr2, XO}, and the references therein. Furthermore, well-posedness results for these types of models, which are highly nonlinear systems, include, for example, existence, uniqueness and regularity of solutions (local, global, weak or radial), qualitative properties of positive solutions, uniform boundedness/blow-up for finite/infinite times and convergence of solutions to an equilibrium, and asymptotic behavior, see e.g., \cite{BDI, BKZ, BFD, Ef, EZ, FS, KK, MWZZ, NO, OS, PX, RB, SY,TW, Wi, ZC}, and the references therein. Recently, several numerical algorithms based on finite-volume and finite-element methods, 
upwind difference methods, operator-splitting methods, discontinuous Galerkin methods or simple-relaxation time lattice Boltzmann methods are used to solve these types of models (see e.g., 
\cite{ABS, CB, CBA, Ep, EK, Fi, Sa, SSK, SST, YSC, ZZL}, and the references therein). 
%ICI

Construction of quantitative models and predict optimal behaviors are among the main objectives of quantitative biology, which aim to understand complex cellular and organismal phenomena, at the system level lie, in a systematic, quantitative, and optimal manner. Optimal behavior is defined as an action that optimizes and evaluates the costs and benefits that influence the outcome of a decision, and contributes to an understanding of the system. The biologists control, in general, their experimental devices by using a certain number of functions or parameters of control which enable them to optimize and to stabilize the system, in accordance with the desired performance and targets. 

To predict the dynamic response of interconnected biological systems from given parameters, data and source terms, requires a mathematical model of the behavior of process under investigation and
physical and biological theory linking the state variables of the model to data and parameters. This prediction of the observation constitutes the so-called direct problem (primal
problem, or also forward problem).

If any of the conditions necessary to define a direct problem are unknown or rather badly known, a control problem (or inverse problem) results, typically when
modeling  physico-chemical biology situations where the model parameters (intervening either in the boundary conditions, in initial conditions or in equations model itself) or material properties are unknown, difficult to measure biologically or partially known. Certain parameters or data can influence considerably the material behavior or modify phenomena in biological or medical matter; then their knowledge is an invaluable help to, for example, medical, biological and biosocial criminology applications. The resolution of these inverse problems thus provides them essential informations 
which are necessary to comprehension of various processes which can intervene in these models. In addition to mathematical model (primal problem), the formulation of control problem in form of a known cost (evaluation or objective)  function ${\cal J}$ is a pre-requirement for achieving desired optimal conditions. This resolution, in order to ensure sufficiently smoothness of cost function with respect to the variation of design parameters and data, needs some regularity and additional conditions, and partial informations of some unknown parameters and fields (observations) given, for example, by experiment measurements. For optimal control in chemotaxis-type systems, we can mention e.g., \cite{B1, B0, LMa, Mg, UYa}, and the references therein.

The problem studied in this paper derives from nonlinear optimal control of a class of dynamic systems which are characterized by chemotaxis. In order to target a wide range of applications, we have written the continuous chemotaxis-type models in a generalized form by considering a class of coupled nonlinear/degenerate anisotropic convection-diffusion chemotaxis-type systems with Neumann-type boundary conditions (for simplicity here we use null flux conditions). We have also opted to place the control variables in all equations of the coupled systems. To facilitate the presentation, we  consider a system with only three equations (the developed analysis will remains valid if we consider a system with $N$ equations ($N\geq 4$) or a coupled systems of reaction-diffusion
and ordinary differential equations, as for example the system modeling atherogenesis or habitat-destruction tumor model). More precisely, we consider the following class of system
\begin{equation}\label{eq:ReacDiff} 
\begin{array}{lcr}
\displaystyle \frac{\partial u}{\partial t}-div(D_{u}(\mathbf{x},t;{\bf u})\nabla u) +div(\mathfrak{T}_{u}(\mathbf{x},t;u))= \Phi_{u}(\mathbf{x},t,{\bf f}_{1};{\bf u})\\
\displaystyle \hspace{2cm}-div(D_{u,1}(\mathbf{x},t;{\bf u})\nabla v)-div(D_{u,2}(\mathbf{x},t;{\bf u})\nabla w) ~~~on~{\cal Q}, \\%[0.2cm]\\

\displaystyle \frac{\partial v}{\partial t}-div(D_{v}(\mathbf{x},t;{\bf u})\nabla v) +div(\mathfrak{T}_{v}(\mathbf{x},t;v))= \Phi_{v}(\mathbf{x},t,{\bf f}_{2};{\bf u})\\
\displaystyle \hspace{2cm}-div(D_{v,1}(\mathbf{x},t;{\bf u})\nabla u)-div(D_{v,2}(\mathbf{x},t;{\bf u})\nabla w)~~~on~{\cal Q}, \\%[0.2cm] \\

\displaystyle \frac{\partial w}{\partial t}-div(D_{w}(\mathbf{x},t;{\bf u})\nabla w) +div(\mathfrak{T}_{w}(\mathbf{x},t;w))= \Phi_{w}(\mathbf{x},t,{\bf f}_{3};{\bf u})\\
\displaystyle \hspace{2cm}-div(D_{w,1}(\mathbf{x},t;{\bf u})\nabla u)-div(D_{w,2}(\mathbf{x},t;{\bf u})\nabla v)~~~on~{\cal Q}, \\%[0.2cm] \\
\text{with the initial condition}\\
(u,v,w)(.,0)=(u_{0},v_{0},w_{0})~~~on~ \Omega,
\end{array}
\end{equation}
subject to the following homogeneous Neumann boundary conditions
\begin{equation}\label{E2}
\begin{array}{lcr}
\displaystyle (D_{u}\nabla u).\n\ =(D_{u,1}\nabla v).\n\ +(D_{u,2}\nabla w).\n\ ~~on~ \Sigma,\\
\displaystyle (D_{v}\nabla v).\n\ =(D_{v,1}\nabla u).\n\ +(D_{v,2}\nabla w).\n\ ~~on~ \Sigma,\\
\displaystyle (D_{w}\nabla w).\n\ =(D_{w,1}\nabla u).\n\ +(D_{w,2}\nabla v).\n\ ~~on~ \Sigma,
\end{array}
\end{equation}
where $(\mathbf{x},t)$ are the space-time variables, $\n\ $ being the outward normal to $\Gamma$, ${\bf f}=({\bf f}_{1},{\bf f}_{2}, {\bf f}_{3})$, where  ${\bf f}_{i} : (\mathbf{x},t)  \ft\ {\bf f}_{i} (\mathbf{x},t) \in \R{d_{i}}$ with $d_{i}\geq 1$, represents the input of the system (source terms, parameters and others) and
${\bf u}=(u, v,w)$ (where $u$,$v$ and $w$ are a scalar function) represents the state or output of the system (e.g., attractant concentrations, attractiveness values, density of the bacterial/cellular/criminal species, nutrient concentrations, etc), $(\Phi_{u},\Phi_{v},\Phi_{w})$ represents the kinetics/source term. The operators $D_{\Theta}$, for $\Theta=u,v$ or $w$ are the self diffusion tensors and $D_{\Theta,k}$, for $\Theta=u,v$ or $w$ and $k=1,3$ are cross-diffusion tensors. These last tensors are assumed to be symmetric and sufficiently regular.
The domain $\Omega$ is an open bounded subset of $\R{m}$, $m\leq 3$, with the Lipschitz-continuous boundary $\Gamma=\p\ \Omega$ and the time $T>0$ is a fixed constant (a given final time). 
Furthermore we define ${\cal Q}=\Omega \times(0,T)$ and $\Sigma=\Gamma\times(0,T)$. 

The nonlinear operators, $D$ (tensor-valued), $\mathfrak{T}$ (vector-valued) and $\Phi$ (scalar-valued) in \eqref{eq:ReacDiff} are assumed to be Carath\'eodory functions.
\begin{req}
If the tensor operators $D_{\Theta}$ and $D_{\Theta,k}$ (for $\Theta=u,v$ or $w$ and $k=1,2$) are of the form $D_{\Theta}=\sigma_{\Theta}{\bf I}_{d}$ and $D_{\Theta,k}=\sigma_{\Theta,k}{\bf I}_{d}$, with ${\bf I}_{d}$ identity matrix and $\sigma_{\Theta}$, $\sigma_{\Theta,k}$ scalar operators, to close the system \eqref{eq:ReacDiff} we can impose the following boundary conditions:
$\displaystyle \nabla u.\n\ =0,~~ \nabla v.\n\ =0 \text{~~~and~~~}\nabla w.\n\ =0~~on~ \Sigma.$ \text{\hfill $\Box$}
\end{req}
\begin{req}
To act in accordance with a desired realistic situation, one can change the boundary conditions \eqref{E2} by considering, for example, mixed boundary conditions (Dirichlet, Neumann, linear/nonlinear Robin-type, etc) on different portions $\Gamma_{i}$ (i=1,p) of $\displaystyle \Gamma=\p\ \Omega=\cup_{i=1,p}\Gamma_{i}$ (such that $\displaystyle \cap_{i=1,p}\Gamma_{i}=\emptyset$). \hfill $\Box$
\end{req}
In order to solve numerically the considered optimal control problems, this work proposes and develops an adjoint-based coupling multiple-relaxation-time lattice Boltzmann model, reliable, efficient, stable and easy to implement in the context of general anisotropic reaction-diffusion systems, with Neumann (Robin) type boundary conditions in complex geometry boundaries, coupled to gradient-based algorithms. While a huge number of papers has been devoted to adjoint optimisation approach for PDE-constrained optimal control problems, there are few papers dedicated to combination with adjoint problem and lattice Boltzmann model, and most of them relates to the fluid dynamics with the classical and simple-relaxation-time lattice Boltzmann method (which is limited to the description of isotropic diffusion problems), see e.g., \cite{HM,KTH,LR,YY} and the references therein. Recently, in \cite{BIJB} we have developed a new mathematical framework which models the competition between tumor and normal cells under chemotherapy constraints in treatment of brain tumors. The goal is to predict the distribution and necessary quantity of drugs delivered in drug-therapy by using a control problem. The control estimates simultaneously blood perfusion rate, reabsorption rate of drug and drug dosage administered, which affect the effects of brain-tumor chemotherapy. In order to solve  numerically 
this problem, we have proposed and investigated an adjoint multiple-relaxation-time lattice Boltzmann method for a class of nonlinear coupled anisotropic convection-diffusion system (which includes the developed model for brain tumor targeted drug delivery system). One of the purposes of this paper is to extend the numerical method developed in \cite{BIJB} to systems derive from the modeling
of directed motion of biological organisms in response to changes of a number of chemical and physical cues in their environment. The main focus of this method is to express the gradient of cost functional from the adjoint of coupled multiple-relaxation-time lattice Boltzmann system. 

Lattice Boltzmann model (LBM) was originated from Boltzmann's kinetic theory of gases (70s), and attracts more and more attentions for simulating complex fluid flows since 90s. More recently, LBM has been extended successfully to simulate different types of parabolic reaction-diffusion equation in biological model as monodomain and bidomain model in cardiac electrophysiology, see \cite{CB, CO} and  
Keller-Segel chemotaxis  model see \cite{YSC} and the references therein. In order to take into account anisotropic diffusion problems, multiple-relaxation-time LBM has been developed for nonlinear anisotropic convection-diffusion equations, see e.g.  \cite{BIJB, HW,ZBZ} and the references therein. Some LBM models have been developed to be kinetic models (LBKM) so that the non-equilibrium flow systems can be better accessed (the kinetic moment relations used in developing this method are physical in comparison with thus used for LBM), see e.g., \cite{LLFS, XLZ} and the references therein. 

The concept of LBM, which is efficient in terms of parallelization, is based on Boltzmann equations which describe the evolution of particles in kinetic theory. Indeed, traditional numerical
methods as finite difference method or finite element method directly solve governing equations for deriving macroscopic variable, whereas LBM is based on the particle (the discrete) distribution
function and numerical solving of continuous Boltzmann transport equation. Then the macroscopic variables of the considered system can be recovered from the discrete equations through multi-scaling
Chapman-Enskog expansion procedure. LBM has two main phases: the first one, local, models the collision between particles and the second one, along each direction of interpolation, models the transport phase.

The rest of paper is organized as follows. In the next section, we start by formulating the control problem in terms of a general control configuration. Afterwards, we establish the (primal-dual) first order optimality conditions for the optimal solution by using sensitivity and adjoint calculus. The optimality system  requires calculation of gradients of the given cost functional which are also necessary to develop numerical optimization algorithms. For ease of understanding, some examples are presented. Some generalities about a possible strategy for numerical realization of optimal control problems based on gradient optimization algorithms (by using adjoint variables) are also given. In Section 3, first, we develop our coupled multiple-relaxation-time LBM (MRT) for numerical simulations of the primal problem and investigate its asymptotic behavior. Next, we describe the adjoint-based strategy, derive a multiple-relaxation-time adjoint system associated to MRT and then  gradients of the cost functional. Finally, comments are presented in Section 4.
%
%%%%%%%%%
%%%%%%%%%%
\section{Study of control framework}
\subsection{Formulation of the problem}
We suppose that there exist ${\cal V}_{0}(\Omega)$, ${\cal U}={\cal U}_{1}\times {\cal U}_{2}\times {\cal U}_{3}$ (where ${\cal U}_{i}=\prod_{j=1}^{d_{i}}{\cal U}^{j}_{i}$, for $i=1,3$) and ${\cal W}({\cal Q})$ three Banach spaces such that for initial variables and input variables, respectively, in ${\cal V}_{0}(\Omega)$ and ${\cal U}$, the problem \eqref{eq:ReacDiff} with boundary conditions \eqref{E2} (under some hypotheses for the data and some regularity of the nonlinear operators) is well-posed in the sense of Hadamard   and the unique solution ${\bf u}$ is in ${\cal W}({\cal Q})$. Moreover, we assume that the nonlinear operator are Fr\'echet differentiable in order to ensure the Fr\'echet differentiability of the control-to-state operator in suitable Banach spaces.

Assume now that the initial condition ${\bf u}_{0}$ is a given function in ${\cal V}_{0}(\Omega)$ (fixed) and introduce the mapping ${\cal F}$ which maps the input ${\bf f}\in {\cal U}$ into the corresponding solution ${\bf u}={\cal F}({\bf f})\in {\cal W}({\cal Q})$ of \eqref{eq:ReacDiff}-\eqref{E2}. We denote by ${\cal K}_{ad}={\cal K}^{ad}_{1}\times {\cal K}^{ad}_{2}\times {\cal K}^{ad}_{3}\subset {\cal U}$ (where ${\cal K}^{ad}_{i}=\prod_{j=1}^{d_{i}}{\cal K}^{ad,j}_{i}$, for $i=1,3$) the set of admissible controls which is a given non-empty, closed and convex  subset of ${\cal U}$.

The nonlinear optimal problem that we are going to study  may be present in a very large class of problems, which can be formally described as follows : 
 find ${\bf f}^{*}\in {\cal K}_{ad}$ such that the optimal criteria which is expressed by the following cost functional ${\cal J}$ (in the reduced form) 

\begin{equation}\label{E15}%\label{COSTLBM}
{\cal J}({\bf f})=\int_{\Omega }J_{0}({\bf u}(T))dx+\int_{0}^{T}\int_{\Omega}{\cal Z}({\bf u},{\bf f})dxdt,
\end{equation}
is minimized with respect to ${\bf f}$ subject to system (\ref{eq:ReacDiff})-(\ref{E2}).
In an other way, we will study 
\begin{equation}\label{EPS} 
\begin{array}{lcr}
\text{Find ${\bf f}^{*}\in {\cal K}_{ad}$ such that}\\
\displaystyle {\cal J}({\bf f}^{*})=\inf_{{\bf f}\in {\cal K}_{ad}}{\cal J}({\bf f}).
\end{array}
\end{equation}
The Fr\'echet differential operator $J_{0}$ and ${\cal Z}$ are assumed to be sufficiently regular to ensure sufficiently smoothness of cost function ${\cal J}$ with respect to the variation of design
${\bf f}$.
\begin{req}
\begin{enumerate}
\item The cost functional ${\cal J}$ describing the control problem depends on function ${\bf f}$ and state function ${\bf u}$ in domain $\Omega$ over the time interval under consideration $[0,T]$. It also depends on a given measurement data $\mathfrak{m}_{data}({\bf x},t)\in \R{s}$ and reference data $\mathfrak{r}_{data}({\bf x},t)\in \R{r}$. In order to simplify the presentation, we have used in expression \eqref{E15}, the reduced form of the functional i.e., ${\cal J}({\bf f})$ in the place of classical form ${\cal J}({\bf f}, {\cal F}({\bf f}))$. 
More classical cost functional $\cal J$ is given by 
\begin{equation}\label{classicJ}
\begin{array}{lcr}
\displaystyle{\cal J}({\bf f})=\frac{a_{1}}{2}\int_{\Omega}\!\!\mid u(T)-u_{f}\mid^{2}dx+\frac{a_{2}}{2}\int_{\Omega}\!\!\mid v(T)-v_{f}\mid^{2}dx+\frac{a_{3}}{2}\int_{\Omega}\!\!\mid w(T)-w_{f}\mid^{2}dx\\
\displaystyle+\frac{b_{1}}{2}\int_{0}^{T}\!\!\!\int_{\Omega}\!\!\mid u-u_{d}\mid^{2}dxdt+\frac{b_{2}}{2}\int_{0}^{T}\!\!\!\int_{\Omega}\!\!\mid v-v_{d}\mid^{2}dxdt+\frac{b_{3}}{2}\int_{0}^{T}\!\!\!\int_{\Omega}\!\!\mid w-w_{d}\mid^{2}dxdt\\
\displaystyle +\sum_{i=1}^{3}\sum_{j=1}^{d_{i}}\alpha_{i,j}\int_{0}^{T}\!\!\!\int_{\Omega}\!\!\mid f^{(j)}_{i}-f^{(j)}_{r,i}\mid^{p_{i,j}}dxdt
\end{array}
\end{equation}
where $p_{i,j}\geq 2$, the functions $f^{(j)}_{r,i}$, for $i=1,3$ and $j=1,d_{i}$, are given reference data, the data $u_{f}, v_{f}, w_{f},u_{d}, v_{d}$ and $w_{d}$ are the target profile for state functions (observations) and the coefficients $a_{k}\geq 0$, $b_{k}\geq 0$ ($k =1,3$) and $\alpha_{i,j}\geq 0$  ($i=1,3$ and $j=1,d_{i}$) such that $\sum_{k=1}^{3}(a_{k}+b_{k})>0$ and $\sum_{i=1}^{3}\sum_{j=1}^{d_{i}}\alpha_{i,j}>0$ are chosen as constants to establish the relative weight of the terms that appear in the definition of the functional.
\item  In accordance with practical applications and available clinical or experimental observations, it is clear that, we can consider other observations and/or controls, which can appear in boundary conditions, in initial conditions or in cross-diffusion tensors (see e.g., \cite{B2,BTT,B8}). The study carried out in this paper, by changing the boundary conditions accordingly,  will remain valid without major modifications.\hfill $\Box$
\end{enumerate}
\end{req}
\begin{req}
If ${\cal K}^{ad,j}_{i}$ is for example of the form $\{g\in L^{2}({\cal Q})~:~m_{1}\leq g\leq m_{2} ~~a.e. ~in~{\cal Q}\}$ (where $m_{1}$ and $m_{2}$ 	are two constants), although it is a subset of $L^{\infty}({\cal Q})$, we prefer to use the standard norms of the space $L^{2}({\cal Q})$. The reason is that we would like to take advantages of the differentiability of the latter norm away from the origin to perform our variational analysis.\hfill $\Box$
\end{req}
%
%%%%%%%%%%%%%%%%%%%%%%%%%%%%%%%%%%%%%%%%%%%%%%%%%%%%%%%%%%%%%%
\subsection{Some examples}
To illustrate the considered general problem, in this paragraph, we give some application examples. 
\subsubsection{Crime problem}
In this application we consider a reaction-diffusion system involving mobile criminal offenders as
\begin{equation}\label{C} 
\begin{array}{lcr}
\displaystyle \frac{\partial u}{\partial t}-div(\sigma \nabla u)= -u+f^{(1)}_{1}uv+f^{(2)}_{1} ~~~on~{\cal Q}, \\%[0.2cm]\\
\displaystyle \frac{\partial v}{\partial t}-div(\nabla v) = -div(2\frac{v}{u}\nabla u)+f_{2}-uv~~~on~{\cal Q}. \\%[0.2cm] \\
\end{array}
\end{equation}
The state $u$ is the attractiveness field, the state $v$ is the expectation value of the criminal density, and the number of burglaries being committed at time $t$ and location ${\bf x}$ is given by 
$u({\bf x},t)v({\bf x},t)$. The parameter  $f^{(1)}_{1}$ is an attractive force pulling offenders back to locations where they have successfully committed crimes. The criminal agents are being created at a rate $f_{2}$ and are removed from the model when a burglary is committed. The function $f^{(2)}_{1}$ represents for example a static, spatially varying, component of the attractiveness. The function ${\bf f}=(f^{(1)}_{1},f^{(2)}_{1},f_{2})$ influences considerably the dynamics of criminals and of attractiveness field, and then their control will be very beneficial and great help  to predict accurately several regimes of aggregation, including hotspots of high criminal activity.
\subsubsection{Attraction-repulsion chemotaxis type model with logistic source}
Our second example concerns an attraction-repulsion chemotaxis-type model with fast diffusion and logistic source
\begin{equation}\label{ARC} 
\begin{array}{lcr}
\displaystyle \frac{\partial u}{\partial t}-div(\sigma_{1}(u) \nabla u)+div(\vec{\omega} u)=-div(\chi_{1}(.;u,v)\nabla  v)-div(\chi_{2}(.;u,w)\nabla  w)+F(u) ~~~on~{\cal Q},\\
\displaystyle \frac{\partial v}{\partial t}-div(\sigma_{2} \nabla v)+div(\vec{\omega} v)=f^{(1)}_{2}u-f^{(2)}_{2}v-f^{(3)}_{2}uv~~~on~{\cal Q},\\
\displaystyle \frac{\partial w}{\partial t}-div(\sigma_{3} \nabla w)+div(\vec{\omega} w)=f^{(1)}_{3}u-f^{(2)}_{3}w+f^{(3)}_{3}~~~on~{\cal Q},
\end{array}
\end{equation}
where $F$ can be given by $F(u)=a_{1}u-a_{2}u^{\kappa}$ (with $\kappa\geq 1$ and the positive parameters $a_{1}$, $a_{2}$ describing respectively the organism growth rate and the carrying capacity), the states $u$, $v$, $w$ are respectively the cell density, the concentration of an attractive chemical signal, and the concentration of a repulsive chemical signal, and the vector field $\vec{\omega}$ is the flow velocity which is generated by a force due to aggregation of cells. The positive parameters representing the chemotactic behavior are $f^{(1)}_{2}$ and $f^{(1)}_{3}$ which represent the 
chemical production rate,  $f^{(2)}_{2}$ and $f^{(2)}_{3}$ which represent the chemical degradation rate, $f^{(3)}_{2}$ which represents the chemical consumed rate, 
$f^{(3)}_{3}$ which corresponds source function, and the matrix-valued operators $\chi_{k}$, $k=1,2$ which represent chemotactic signal of attraction and repulsion, respectively.  The operator $\sigma_{1}$ can be  described for example by $\sigma_{1}(u)=\frac{1}{(1-u)^{\alpha}}$ with $\alpha\geq 2$, by $\sigma_{1}(u)=\frac{1+u}{(1-u+ln(u))}$ or by $\sigma_{1}(u)=u^{p-1}$ (with $p>1$). The operators 
$\sigma_{2}$ and $\sigma_{3}$ are assumed to be constant. Finally, the operator $\chi_{k}$ can be described for example by $\chi_{k}({\bf x},t;u,\zeta)=\mathbb{X}({\bf x},t)u{\bf I}_{d}$, or by 
$\chi_{k}({\bf x},t;u,\zeta)=\mathbb{X}({\bf x},t)\frac{u}{1+\zeta^{2}}{\bf I}_{d}$, where $\mathbb{X}$ is a sufficiently regular function. 

We can control some of parameters/source  $f^{(j)}_{i}$ from par example observations given by measurements.
\subsubsection{Two species chemotaxis blow up}
The third example deals with two-species chemotaxis system with competitive kinetics
\begin{equation}\label{TS} 
\begin{array}{lcr}
\displaystyle \frac{\partial u}{\partial t}-div(\sigma_{1}(u) \nabla u)+div(\vec{\omega} u)=-div(\chi_{1}(.;u,w)\nabla  w)+F_{1}(u,v) ~~~on~{\cal Q},\\
\displaystyle \frac{\partial v}{\partial t}-div(\sigma_{2}(v) \nabla v)+div(\vec{\omega} v)=-div(\chi_{2}(.;v,w)\nabla  w) +F_{2}(u,v)~~~on~{\cal Q},\\
\displaystyle \frac{\partial w}{\partial t}-div(\sigma_{3} \nabla w)+div(\vec{\omega} w)=F_{3}(u,v,w;{\bf f}_{3})~~~on~{\cal Q},
\end{array}
\end{equation}
 where $F_{i}$ (for $i=1,3$) can be given by $F_{1}(u,v)=\mu_{1}u(1-u-a_{1}v)$, $F_{2}(u,v)=\mu_{2}v(1-v-a_{2}u)$ and $F_{3}(u,v,w,;{\bf f}_{3})=-w+f^{(1)}_{3}u+f^{(2)}_{3}v +f^{(3)}_{3} $or $-w(f^{(1)}_{3}u+f^{(2)}_{3}v)+f^{(3)}_{3}$ with ${\bf f}_{3}=(f^{(1)}_{3},f^{(2)}_{3},f^{(3)}_{3})$, the state $u$ and $v$ denote the cell densities of the first and second species, respectively, and $w$ presents the concentration of the chemical signal. The parameters present in the model are positive and the operators are of the same type as those of the previous application.
\subsubsection{Tumor invasion into surrounding environment and habitat-destruction type model}
In our last example, we propose a modified version of the haptotaxis model introduced in \cite{CL}. The model, which takes into account the competition for space between
cancer cells and extracellular matrix  macromolecules (ECM), consists of an ODE reflecting the degradation of ECM, coupled with a reaction-diffusion chemotaxis-haptotaxis system describing the evolution of cancer cell density (which uses ECM for movement), the matrix degrading enzyme concentration and the oxygen concentration (which is produced by the matrix molecules). The system is of the form 
\begin{equation}\label{TIS} 
\begin{array}{lcr}
\displaystyle \frac{\partial u}{\partial t}-div(\sigma_{1}(u) \nabla u)+div(\vec{\omega} u)=\mu_{1}\frac{uw}{w+\theta}-\mu_{2}u+\mu_{3}u(1-m-u)\\
\displaystyle \hspace{3cm}-div(u\chi_{1}(.;m)\nabla  m)-div(u\chi_{2}(.;v)\nabla  v) ~~~on~{\cal Q},\\
\displaystyle \frac{\partial v}{\partial t}-div(\sigma_{2}\nabla v)+div(\vec{\omega} v)=r_{1}u(1-u) -r_{2}v+f~~~on~{\cal Q},\\
\displaystyle \frac{\partial w}{\partial t}-div(\sigma_{3} \nabla w)+div(\vec{\omega} w)=\gamma_{1}m-\gamma_{2}w-\gamma_{3}\frac{uw}{w+\theta}+g~~~on~{\cal Q},\\
\displaystyle \displaystyle \frac{\partial m}{\partial t}=-\nu_{1} w m+\nu_{2}u(1-m-u)~~~on~{\cal Q},
\end{array}
\end{equation}
where the states $u,v,w$ and $m$ are cancer cell density, matrix degradation enzyme concentration, oxygen concentration and extracellular matrix density, respectively. The parameters present in the model are positive.

We can identify (estimate) the best optimal prognostic values of some parameters appearing in the model or control the functions $f$ and $g$ from observations (desired target).
%
%
%%%%%%%%%%%%%%%%%%%%%%%%%%%%%%%%%%%%%%%%%%%%%%%%%%%%%%%%%%%%%%%%%%%%%%%%%%%%%ù
%
%%%%%%%%%%%%%%%%%%%%%%%%%%%%%%%%%%%%%%%%%%%%%%%%%%%%%%%%%%%%%%%%%%%%%%%%%
\subsection{Continuous adjoint-based optimality conditions}
Assume that the nonlinear control problem (\ref{EPS}) admits an optimal solution; the necessary conditions
for this optimum is given by the following theorem (see \cite{B0}).

\begin{theo} If ${\cal J}$ attains a (local) minimum at a point ${\bf f}^{*}\in {\cal K}_{ad}$, then the following first optimality conditions hold:
\begin{equation}\label{CO}
{\cal J}'({\bf f}^{*}).({\bf f}-{\bf f}^{*})\geq 0,~~\forall {\bf f}\in {\cal K}_{ad},
\end{equation}
where ${\cal J}'$ is the directional derivative of ${\cal J}$.
\end{theo}
In order to solve numerically the optimal control problem, it is necessary to derive the gradient of the cost functional ${\cal J}$𝐽 with respect to the control ${\bf f}=({\bf f}_{1},{\bf f}_{2},{\bf f}_{3})$. For this, we suppose that the operator solution ${\cal F}$ is continuously differentiable on ${\cal K}_{ad}$ and its derivative,𝑤 at ${\bf f}$ in direction ${\bf h}=({\bf h}_{1},{\bf h}_{2},{\bf h}_{3})$, $\displaystyle {\wp}=(p,q,\rho)= {\cal F}'({\bf f}).{\bf h}=\lim_{\theta\ft\ 0^{+}}\frac{({\cal F}({\bf f}+\theta {\bf h})-{\cal F}({\bf f}))}{𝜖\theta}$
is the unique solution of the following system (which is called the tangent linear model (TLM) or sensitivity problem)
%\newpage
\begin{equation}\label{E1d} 
\begin{array}{lcr}
\displaystyle \frac{\partial p}{\partial t}-div(D_{u}(.;{\bf u})\nabla p)-div((\tilde{H}_{u;u}(.;{\bf u}, \nabla {\bf u})-\frac{\p\ \mathfrak{T}_{u}}{\p\ u}(.;u))p)-\frac{\p\ \Phi_{u}}{\p\ u}(.; {\bf f}_{1},{\bf u})p\\
\displaystyle \hspace{1.5cm}=div(\tilde{H}_{u;v}(.;{\bf u}, \nabla {\bf u})q)+div(\tilde{H}_{u;w}(.;{\bf u}, \nabla {\bf u})\rho)
+\frac{\p\ \Phi_{u}}{\p\ v}(.; {\bf f}_{1},{\bf u})q+\frac{\p\ \Phi_{u}}{\p\ w}(.; {\bf f}_{1},{\bf u})\rho \\
\displaystyle \hspace{2cm}+\frac{\p\ \Phi_{u}}{\p\ {\bf f}_{1}}(.; {\bf f}_{1},{\bf u}){\bf h}_{1}-div(D_{u,1}(.;{\bf u})\nabla q)-div(D_{u,2}(.;{\bf u})\nabla \rho),\\ %[0.2cm] \\
\displaystyle \frac{\partial q}{\partial t}-div(D_{v}(.;{\bf u})\nabla q)-div((\tilde{H}_{v;v}(.;{\bf u}, \nabla {\bf u})-\frac{\p\ \mathfrak{T}_{v}}{\p\ v}(.;v))q)-\frac{\p\ \Phi_{v}}{\p\ v}(.; {\bf f}_{2},{\bf u})q\\
\displaystyle \hspace{1.5cm}=div(\tilde{H}_{v;u}(.;{\bf u}, \nabla {\bf u})p)+div(\tilde{H}_{v;w}(.;{\bf u}, \nabla {\bf u})\rho)
+\frac{\p\ \Phi_{v}}{\p\ u}(.; {\bf f}_{2},{\bf u})p+\frac{\p\ \Phi_{v}}{\p\ w}(.;{\bf f}_{2},{\bf u})\rho \\
\displaystyle \hspace{2cm}+\frac{\p\ \Phi_{v}}{\p\ {\bf f}_{2}}(.; {\bf f}_{2},{\bf u}){\bf h}_{2}-div(D_{v,1}(.;{\bf u})\nabla p)-div(D_{v,2}(.;{\bf u})\nabla \rho),\\ %[0.2cm] \\
\displaystyle \frac{\partial \rho}{\partial t}-div(D_{w}(.;{\bf u})\nabla \rho))-div((\tilde{H}_{w;w}(.;{\bf u}, \nabla {\bf u})-\frac{\p\ \mathfrak{T}_{w}}{\p\ w}(.;w))\rho)-\frac{\p\ \Phi_{w}}{\p\ w}(.; {\bf f}_{3},{\bf u})\rho\\
\displaystyle \hspace{1.5cm}=div(\tilde{H}_{w;u}(.;{\bf u}, \nabla {\bf u})p)+div(\tilde{H}_{w;v}(.;{\bf u}, \nabla {\bf u})q)+\frac{\p\ \Phi_{w}}{\p\ u}(.; {\bf f}_{3},{\bf u})p+\frac{\p\ \Phi_{w}}{\p\ v}(.;{\bf f}_{3},{\bf u})q \\
\displaystyle \hspace{2cm}+\frac{\p\ \Phi_{w}}{\p\ {\bf f}_{3}}(.; {\bf f}_{3},{\bf u}){\bf h}_{3}-div(D_{w,1}(.;{\bf u})\nabla p)-div(D_{w,2}(.;{\bf u})\nabla q), %[0.2cm] \\
\end{array}
\end{equation}
under the initial and Robin-type boundary conditions: 
\begin{equation}\label{E1dB} 
\begin{array}{lcr}
(u,v,w)(.,0)=0~~~on~ \Omega,\\
\displaystyle (D_{u}\nabla p).\n\ +\left(p\tilde{H}_{u;u}(.;{\bf u}, \nabla {\bf u})
+q\tilde{H}_{u;v}(.;{\bf u}, \nabla {\bf u})+\rho\tilde{H}_{u;w}(.;{\bf u}, \nabla {\bf u})\right).\n\ \\
\displaystyle \hspace{1.5cm}=(D_{u,1}\nabla q).\n\ +(D_{u,2}\nabla \rho).\n\ , ~~~on~  \Sigma\\
\displaystyle (D_{v}\nabla q).\n\ +\left(p\tilde{H}_{v;u}(.;{\bf u}, \nabla {\bf u})+q\tilde{H}_{v;v}(.;{\bf u}, \nabla {\bf u})+\rho\tilde{H}_{v;w}(.;{\bf u}, \nabla {\bf u})\right).\n\ \\ 
\displaystyle\hspace{1.5cm}=(D_{v,1}\nabla p).\n\ +(D_{v,2}\nabla \rho).\n\ ,  ~~~on~  \Sigma\\
\displaystyle (D_{w}\nabla \rho).\n\ +\left(p\tilde{H}_{w;u}(.;{\bf u}, \nabla {\bf u})+q\tilde{H}_{w;v}(.;{\bf u}, \nabla {\bf u})+\rho\tilde{H}_{w;w}(.;{\bf u}, \nabla {\bf u})\right).\n\ \\
\displaystyle \hspace{1.5cm}=(D_{w,1}\nabla p).\n\ +(D_{w,2}\nabla q).\n\ , ~~~on~  \Sigma
\end{array}
\end{equation}
where  
\begin{equation*}
\begin{array}{lcr}
\displaystyle \tilde{D}_{u;z}(.;{\bf u}, \nabla u)= \frac{\p\ D_{u}}{\p\ z}(.;{\bf u})\nabla u,~~ 
\displaystyle \tilde{D}_{v;z}(.;{\bf u}, \nabla v)= \frac{\p\ D_{v}}{\p\ z}(.;{\bf u})\nabla v,\\
\displaystyle \tilde{D}_{w;z}(.;{\bf u}, \nabla w)= \frac{\p\ D_{w}}{\p\ z}(.;{\bf u})\nabla w,\\
\displaystyle \tilde{T}_{u;z}(.;{\bf u}, \nabla v, \nabla w)=\frac{\p\ D_{u,1}}{\p\ z}(.;{\bf u})\nabla v+\frac{\p\ D_{u,2}}{\p\ z}(.;{\bf u})\nabla w,\\
\displaystyle \tilde{T}_{v;z}(.;{\bf u}, \nabla u, \nabla w)=\frac{\p\ D_{v,1}}{\p\ z}(.;{\bf u})\nabla u+\frac{\p\ D_{v,2}}{\p\ z}(.;{\bf u})\nabla w,\\
\displaystyle \tilde{T}_{w;z}(.;{\bf u}, \nabla u, \nabla v)=\frac{\p\ D_{w,1}}{\p\ z}(.;{\bf u})\nabla u+\frac{\p\ D_{w,2}}{\p\ z}(.;{\bf u})\nabla v,\\

\displaystyle \tilde{H}_{u;z}(.;{\bf u}, \nabla {\bf u})=\tilde{D}_{u;z}-\tilde{T}_{u;z},~~\tilde{H}_{v;z}(.;{\bf u}, \nabla {\bf u})=\tilde{D}_{v;z}-\tilde{T}_{v;z},\\
\displaystyle\tilde{H}_{w;z}(.;{\bf u}, \nabla {\bf u})=\tilde{D}_{w;z}-\tilde{T}_{w;z}.
\end{array}
\end{equation*}
We can now show the first-order necessary conditions (optimality conditions) and calculate the gradient of ${\cal J}$ by using TLM and by introducing an intermediate adjoint (dual ou costate)  model. 
\begin{theo}
Suppose that $({\bf u}^{*}, {\bf f}^{*})\in {\cal W}({\cal Q})\times {\cal K}_{ad}$, is an optimal solution such that ${\bf f}^{*}=({\bf f}^{*}_{1}, {\bf f}^{*}_{2}, {\bf f}^{*}_{3})$ is defined by (\ref{EPS}) and ${\bf u}^{*}={\cal F}({\bf f}^{*})$ is the solution of \eqref{eq:ReacDiff}-\eqref{E2}. Then there exists an adjoint variable $\tilde{\bf u}^{*}=( \tilde{u}^{*},\tilde{v}^{*},\tilde{w}^{*})$ 
satisfying the following adjoint problem
\begin{equation}\label{FVM-PF-Adj}
\begin{array}{lcr}
\displaystyle -\frac{\partial \tilde{u}^{*}}{\partial t}-div(D_{u}(.;{\bf u}^{*})\nabla \tilde{u}^{*})-\frac{\p\ \Phi_{u}}{\p\ u}(.; {\bf f}^{*}_{1},{\bf u}^{*})\tilde{u}^{*}
+(\tilde{H}_{u;u}(.;{\bf u}^{*}, \nabla {\bf u}^{*})-\frac{\p\ \mathfrak{T}_{u}}{\p\ u}(.;u^{*}))\nabla \tilde{u}^{*}\\
\displaystyle\hspace{1.cm} +\tilde{H}_{v;u}(.;{\bf u}^{*}, \nabla {\bf u}^{*})\nabla \tilde{v}^{*} -\frac{\p\ \Phi_{v}}{\p\ u}(.;{\bf f}^{*}_{2},{\bf u}^{*}) \tilde{v}^{*} 
\displaystyle+\tilde{H}_{w;u}(.;{\bf u}^{*}, \nabla {\bf u}^{*})\nabla \tilde{w}^{*} -\frac{\p\ \Phi_{w}}{\p\ u}(.; {\bf f}^{*}_{3},{\bf u}^{*}) \tilde{w}^{*}\\
\displaystyle\hspace{1.7cm}=-div(D_{v,1}(.;{\bf u}^{*})\nabla \tilde{v}^{*}) -div(D_{w,1}(.;{\bf u}^{*})\nabla \tilde{w}^{*})+\frac{\p\ {\cal Z}}{\p\ u}({\bf u}^{*},{\bf f}^{*}),\\ %[0.2cm] \\
\displaystyle -\frac{\partial \tilde{v}^{*}}{\partial t}-div(D_{v}(.;{\bf u}^{*})\nabla \tilde{v}^{*})-\frac{\p\ \Phi_{v}}{\p\ v}(.; {\bf f}^{*}_{2},{\bf u}^{*})\tilde{v}^{*}
+(\tilde{H}_{v;v}(.;{\bf u}^{*}, \nabla {\bf u}^{*})-\frac{\p\ \mathfrak{T}_{v}}{\p\ v}(.;v^{*}))\nabla \tilde{v}^{*}\\
\displaystyle\hspace{1.cm}+\tilde{H}_{u;v}(.;{\bf u}^{*}, \nabla {\bf u}^{*})\nabla \tilde{u}^{*} -\frac{\p\ \Phi_{u}}{\p\ v}(.; {\bf f}^{*}_{1},{\bf u}^{*}) \tilde{u}^{*}
+\tilde{H}_{w;v}(.;{\bf u}^{*}, \nabla {\bf u}^{*})\nabla \tilde{w}^{*} -\frac{\p\ \Phi_{w}}{\p\ v}(.; {\bf f}^{*}_{3},{\bf u}^{*}) \tilde{w}^{*}\\
\displaystyle\hspace{1.7cm}=-div(D_{u,1}(.;{\bf u}^{*})\nabla \tilde{u}^{*})-div(D_{w,2}(.;{\bf u}^{*})\nabla \tilde{w}^{*})+\frac{\p\ {\cal Z}}{\p\ v}({\bf u}^{*},{\bf f}^{*}),\\ %[0.2cm] \\
\displaystyle -\frac{\partial \tilde{w}^{*}}{\partial t}-div(D_{w}(.;{\bf u}^{*})\nabla \tilde{w}^{*})-\frac{\p\ \Phi_{w}}{\p\ w}(.; {\bf f}^{*}_{3},{\bf u}^{*})\tilde{w}^{*}
\displaystyle+(\tilde{H}_{w;w}(.;{\bf u}^{*}, \nabla {\bf u}^{*})-\frac{\p\ \mathfrak{T}_{w}}{\p\ w}(.;w^{*}))\nabla \tilde{w}^{*}\\
\displaystyle\hspace{1.cm}+\tilde{H}_{u;w}(.;{\bf u}^{*}, \nabla {\bf u}^{*})\nabla \tilde{u}^{*} -\frac{\p\ \Phi_{u}}{\p\ w}(.; {\bf f}^{*}_{1},{\bf u}^{*}) \tilde{u}^{*} 
\displaystyle+\tilde{H}_{v;w}(.;{\bf u}^{*}, \nabla {\bf u}^{*})\nabla \tilde{v}^{*} -\frac{\p\ \Phi_{v}}{\p\ w}(.; {\bf f}^{*}_{2},{\bf u}^{*}) \tilde{v}^{*} \\
\displaystyle\hspace{1.7cm}=-div(D_{u,2}(.;{\bf u})\nabla \tilde{u}^{*})-div(D_{v,2}(.;{\bf u}^{*})\nabla \tilde{v}^{*})+\frac{\p\ {\cal Z}}{\p\ w}({\bf u}^{*},{\bf f}^{*}),\\ %[0.2cm] \\
\text{under the boundary conditions }\\
\displaystyle (\tilde{u}^{*}\frac{\p\ \mathfrak{T}_{u}}{\p\ u}(.;u^{*})-\tilde{u}^{*}\tilde{H}_{u;u}(.;{\bf u}^{*}, \nabla {\bf u}^{*})+D_{u}(.;{\bf u}^{*})\nabla \tilde{u}^{*}).\n\ \\
\displaystyle\hspace{1.5cm}=(D_{v,1}(.;{\bf u}^{*})\nabla \tilde{v}^{*}+D_{w,1}(.;{\bf u}^{*})\nabla \tilde{w}^{*}).\n\ =0,\\
\displaystyle (\tilde{v}^{*}\frac{\p\ \mathfrak{T}_{v}}{\p\ v}(.;v^{*})-\tilde{v}^{*}\tilde{H}_{v;v}(.;{\bf u}^{*}, \nabla {\bf u}^{*})+D_{v}(.;{\bf u}^{*})\nabla \tilde{v}^{*}).\n\ \\
\displaystyle\hspace{1.5cm}=(D_{u,1}(.;{\bf u}^{*})\nabla \tilde{u}^{*}+D_{w,2}(.;{\bf u}^{*})\nabla \tilde{w}^{*}).\n\ =0,\\
\displaystyle (\tilde{w}^{*}\frac{\p\ \mathfrak{T}_{w}}{\p\ v}(.;w^{*})-\tilde{w}^{*}\tilde{H}_{w;w}(.;{\bf u}^{*}, \nabla {\bf u}^{*})+D_{w}(.;{\bf u}^{*})\nabla \tilde{w}^{*}).\n\ \\
\displaystyle\hspace{1.5cm}=(D_{u,2}(.;{\bf u}^{*})\nabla \tilde{u}^{*}-D_{v,2}(.;{\bf u}^{*})\nabla \tilde{v}^{*}).\n\ =0\\
\text{and the final condition}\\ 
\displaystyle \tilde{\bf u}^{*}(T) =\frac{\p\ J_{0}}{\p\ {\bf u}}({\bf u}^{*}(T))
\end{array}
\end{equation}
and the following variational inequality (for all ${\bf f}=({\bf f}_{1}, {\bf f}_{2}, {\bf f}_{3})\in {\cal K}_{ad}$)
\begin{equation}\label{COS}
\begin{array}{lcr}
\displaystyle \int_{0}^{T}\!\!\!\int_{\Omega}(\frac{\p\ {\cal Z}}{\p\ {\bf f}_{1}}({\bf u}^{*},{\bf f}^{*})+\tilde{u}^{*}\frac{\p\ \Phi_{u}}{\p\ {\bf f}_{1}}(.; {\bf f}^{*}_{1},{\bf u}^{*})).({\bf f}_{1}-{\bf f}^{*}_{1})dxdt\geq 0,\\
\displaystyle \int_{0}^{T}\!\!\!\int_{\Omega}(\frac{\p\ {\cal Z}}{\p\ {\bf f}_{2}}({\bf u}^{*},{\bf f}^{*})+\tilde{v}^{*}\frac{\p\ \Phi_{v}}{\p\ {\bf f}_{2}}(.; {\bf f}^{*}_{2},{\bf u}^{*})).({\bf f}_{2}-{\bf f}^{*}_{2})dxdt\geq 0,\\
\displaystyle \int_{0}^{T}\!\!\!\int_{\Omega}(\frac{\p\ {\cal Z}}{\p\ {\bf f}_{3}}({\bf u}^{*},{\bf f}^{*})+\tilde{w}^{*}\frac{\p\ \Phi_{w}}{\p\ {\bf f}_{3}}(.; {\bf f}^{*}_{3},{\bf u}^{*})).({\bf f}_{3}-{\bf f}^{*}_{3})dxdt\geq 0.
\end{array}
\end{equation}
\end{theo}
{\bf Proof}. By using the same technique as in \cite{B0}, we start by calculating the variation of ${\cal J}$.

According to the regularity of operator solution ${\cal F}$ and the nature of cost function ${\cal J}$ (which is the composition of Fr\'echet
differentiable mappings), we have that ${\cal J}$ is differentiable and the directional derivative of ${\cal J}$ at point ${\bf f}$ along the direction ${\bf h}=({\bf h}_{1},{\bf h}_{2},{\bf h}_{3})$ can be given by
\begin{equation}\label{DJ}
\begin{array}{lcr}
\displaystyle {\cal J}'({\bf f}).{\bf h}=\lim_{\theta\ft\ 0^{+}}\frac{({\cal J}({\bf f}+\theta {\bf h})-{\cal J}({\bf f}))}{𝜖\theta}\\
\hspace{2cm}\displaystyle=\int_{\Omega}\frac{\p\ J_{0}}{\p\ {\bf u}}({\bf u}(T)).\wp(T)dx+\int_{0}^{T}\!\!\int_{\Omega}\frac{\p\ {\cal Z}}{\p\ {\bf u}}({\bf u},{\bf f}).\wp~dxdt\\
\displaystyle \hspace{2cm}+\int_{0}^{T}\!\!\int_{\Omega}\frac{\p\ {\cal Z}}{\p\ {\bf f}}({\bf u},{\bf f}).{\bf h}~ dxdt.
\end{array}
\end{equation}

Now, we simplify the directional derivative of ${\cal J}$. For this we multiply the first part of \eqref{E1d} by some regular function $\tilde{{\bf u}}=(\tilde{u},\tilde{v},\tilde{w})$, integrating
over ${\cal Q}$, using Green's formula and integrating by parts in times, we obtain (according to the boundary and initial conditions for ${\wp}$)
\begin{equation*}
\begin{array}{lcr}
\displaystyle \int_{\Omega}p(T)\tilde{u} (T)dx+\int_{0}^{T}\!\!\!\int_{\Gamma}p(\frac{\p\ \mathfrak{T}_{u}}{\p\ u}(.;u).\n\  )\tilde{u} d\Gamma dt
+\int_{0}^{T}\!\!\!\int_{\Gamma}p(D_{u}(.;{\bf u})\nabla \tilde{u}).\n\ d\Gamma dt\\
\displaystyle -\int_{0}^{T}\!\!\!\int_{\Gamma}p(\tilde{H}_{u;u}.\n\  )\tilde{u} d\Gamma dt
\displaystyle-\int_{0}^{T}\!\!\!\int_{\Gamma}q(D_{u,1}(.;{\bf u})\nabla \tilde{u}).\n\  d\Gamma dt-\int_{0}^{T}\!\!\!\int_{\Gamma}\rho(D_{u,2}(.;{\bf u})\nabla \tilde{u}).\n\ d\Gamma dt\\
\displaystyle +\int_{0}^{T}\!\!\!\int_{\Omega}p\left[-\frac{\partial \tilde{u}}{\partial t}-div(D_{u}(.;{\bf u})\nabla \tilde{u})-\frac{\p\ \Phi_{u}}{\p\ u}(.; {\bf f}_{1},{\bf u})\tilde{u}
+(\tilde{H}_{u;u}-\frac{\p\ \mathfrak{T}_{u}}{\p\ u}(.;u))\nabla \tilde{u} \right]dxdt\\
\displaystyle=\int_{0}^{T}\!\!\!\int_{\Omega}q\left[-\tilde{H}_{u;v}\nabla \tilde{u} +\frac{\p\ \Phi_{u}}{\p\ v}(.; {\bf f}_{1},{\bf u}) \tilde{u} \right]dxdt
+\int_{0}^{T}\!\!\!\int_{\Omega}\rho\left[-\tilde{H}_{u;w}\nabla \tilde{u} +\frac{\p\ \Phi_{u}}{\p\ w}(.; {\bf f}_{1},{\bf u}) \tilde{u} \right]dxdt\\
\displaystyle+\int_{0}^{T}\!\!\!\int_{\Omega}\left[\frac{\p\ \Phi_{u}}{\p\ {\bf f}_{1}}(.; {\bf f}_{1},{\bf u}){\bf h}_{1}\right]\tilde{u}dxdt\\
\displaystyle-\int_{0}^{T}\!\!\!\int_{\Omega}div(D_{u,1}(.;{\bf u})\nabla \tilde{u})q dxdt-\int_{0}^{T}\!\!\!\int_{\Omega}div(D_{u,2}(.;{\bf u})\nabla \tilde{u})\rho dxdt,  %[0.2cm] \\
\end{array}
\end{equation*}
\begin{equation*}
\begin{array}{lcr}
 %%%%%%%%%Eq2
%
\displaystyle \int_{\Omega}q(T)\tilde{v} (T)dx+\int_{0}^{T}\!\!\!\int_{\Gamma}q(\frac{\p\ \mathfrak{T}_{v}}{\p\ v}(.;v).\n\  )\tilde{v} d\Gamma dt
+\int_{0}^{T}\!\!\!\int_{\Gamma}q(D_{v}(.;{\bf u})\nabla \tilde{v}).\n\ d\Gamma dt\\
\displaystyle-\int_{0}^{T}\!\!\!\int_{\Gamma}q(\tilde{H}_{v;v}.\n\  )\tilde{v} d\Gamma dt
\displaystyle-\int_{0}^{T}\!\!\!\int_{\Gamma}p(D_{v,1}(.;{\bf u})\nabla \tilde{v}).\n\  d\Gamma dt-\int_{0}^{T}\!\!\!\int_{\Gamma}\rho(D_{v,2}(.;{\bf u})\nabla \tilde{v}).\n\ d\Gamma dt\\
\displaystyle +\int_{0}^{T}\!\!\!\int_{\Omega}q\left[-\frac{\partial \tilde{v}}{\partial t}-div(D_{v}(.;{\bf u})\nabla \tilde{v})-\frac{\p\ \Phi_{v}}{\p\ v}(.; {\bf f}_{2},{\bf u})\tilde{v}
+(\tilde{H}_{v;v}-\frac{\p\ \mathfrak{T}_{v}}{\p\ v}(.;v))\nabla \tilde{v} \right]dxdt\\
\displaystyle=\int_{0}^{T}\!\!\!\int_{\Omega}p\left[-\tilde{H}_{v;u}\nabla \tilde{v} +\frac{\p\ \Phi_{v}}{\p\ u}(.; {\bf f}_{2},{\bf u}) \tilde{v} \right]dxdt
+\int_{0}^{T}\!\!\!\int_{\Omega}\rho\left[-\tilde{H}_{v;w}\nabla \tilde{v} +\frac{\p\ \Phi_{v}}{\p\ w}(.; {\bf f}_{2}, {\bf u}) \tilde{v} \right]dxdt\\
\displaystyle+\int_{0}^{T}\!\!\!\int_{\Omega}\left[\frac{\p\ \Phi_{v}}{\p\ {\bf f}_{2}}(.; {\bf f}_{2},{\bf u}){\bf h}_{2}\right]\tilde{v}dxdt\\
\displaystyle-\int_{0}^{T}\!\!\!\int_{\Omega}div(D_{v,1}(.;{\bf u})\nabla \tilde{v})p dxdt-\int_{0}^{T}\!\!\!\int_{\Omega}div(D_{v,2}(.;{\bf u})\nabla \tilde{v})\rho dxdt
 %
%%%%%%%%%%%Eq3
\end{array}
\end{equation*}
and
\begin{equation*}
\begin{array}{lcr}
\displaystyle \int_{\Omega}\rho(T)\tilde{w} (T)dx+\int_{0}^{T}\!\!\!\int_{\Gamma}\rho(\frac{\p\ \mathfrak{T}_{w}}{\p\ w}(.;w).\n\  )\tilde{w} d\Gamma dt
+\int_{0}^{T}\!\!\!\int_{\Gamma}\rho(D_{w}(.;{\bf u})\nabla \tilde{w}).\n\ d\Gamma dt\\
\displaystyle-\int_{0}^{T}\!\!\!\int_{\Gamma}\rho(\tilde{H}_{w;w}.\n\  )\tilde{w} d\Gamma dt
\displaystyle-\int_{0}^{T}\!\!\!\int_{\Gamma}p(D_{w,1}(.;{\bf u})\nabla \tilde{w}).\n\  d\Gamma dt-\int_{0}^{T}\!\!\!\int_{\Gamma}q(D_{w,2}(.;{\bf u})\nabla \tilde{w}).\n\  d\Gamma dt\\
\displaystyle +\int_{0}^{T}\!\!\!\int_{\Omega}\rho\left[-\frac{\partial \tilde{v}}{\partial t}-div(D_{w}(.;{\bf u})\nabla \tilde{w})-\frac{\p\ \Phi_{w}}{\p\ w}(.; {\bf f}_{3},{\bf u})\tilde{w}
+(\tilde{H}_{w;w}-\frac{\p\ \mathfrak{T}_{w}}{\p\ w}(.;w))\nabla \tilde{w} \right]dxdt\\
\displaystyle=\int_{0}^{T}\!\!\!\int_{\Omega}p\left[-\tilde{H}_{w;u}\nabla \tilde{w} +\frac{\p\ \Phi_{w}}{\p\ u}(.; {\bf f}_{3},{\bf u})\tilde{w} \right]dxdt
+\int_{0}^{T}\!\!\!\int_{\Omega}q\left[-\tilde{H}_{w;v}\nabla \tilde{w} +\frac{\p\ \Phi_{w}}{\p\ v}(.; {\bf f}_{3},{\bf u}) \tilde{w} \right]dxdt\\
\displaystyle+\int_{0}^{T}\!\!\!\int_{\Omega}\left[\frac{\p\ \Phi_{w}}{\p\ {\bf f}_{3}}({\bf f}_{3},{\bf u}){\bf h}_{3}\right]\tilde{w}dxdt\\
\displaystyle-\int_{0}^{T}\!\!\!\int_{\Omega}div(D_{w,1}(.;{\bf u})\nabla \tilde{w})p dxdt-\int_{0}^{T}\!\!\!\int_{\Omega}div(D_{w,2}(.;{\bf u})\nabla \tilde{w})q dxdt.
\end{array}
\end{equation*}
Assume that the function $\tilde{\bf u}$ satisfies the following boundary conditions on $\Sigma$ :
\begin{equation} \label{ACB}
\begin{array}{lcr}
\displaystyle (\tilde{u}\frac{\p\ \mathfrak{T}_{u}}{\p\ u}(.;u)-\tilde{u}\tilde{H}_{u;u}+D_{u}(.;{\bf u})\nabla \tilde{u}-D_{v,1}(.;{\bf u})\nabla \tilde{v}-D_{w,1}(.;{\bf u})\nabla \tilde{w}).\n\ =0,\\
\displaystyle (\tilde{v}\frac{\p\ \mathfrak{T}_{v}}{\p\ v}(.;v)-\tilde{v}\tilde{H}_{v;v}+D_{v}(.;{\bf u})\nabla \tilde{v}-D_{u,1}(.;{\bf u})\nabla \tilde{u}-D_{w,2}(.;{\bf u})\nabla \tilde{w}).\n\ =0,\\
\displaystyle (\tilde{w}\frac{\p\ \mathfrak{T}_{w}}{\p\ w}(.;w)-\tilde{w}\tilde{H}_{w;w}+D_{w}(.;{\bf u})\nabla \tilde{w}-D_{u,2}(.;{\bf u})\nabla \tilde{u}-D_{v,2}(.;{\bf u})\nabla \tilde{v}).\n\ =0
\end{array}
\end{equation}
and the final condition 
\begin{equation} \label{ACF}
\begin{array}{lcr}
\displaystyle \tilde{\bf u}(T) =\frac{\p\ J_{0}}{\p\ {\bf u}}({\bf u}(T)).
\end{array}
\end{equation}
Then, by summing the three relations of previous system, we obtain (according to \eqref{ACB} and \eqref{ACF})
\begin{equation}\label{Int}
\begin{array}{lcr}
\displaystyle \int_{\Omega}\frac{\p\ J_{0}}{\p\ {\bf u}}({\bf u}(T)).\wp (T)dx
\displaystyle-\int_{0}^{T}\!\!\!\int_{\Omega}\tilde{u}\frac{\p\ \Phi_{u}}{\p\ {\bf f}_{1}}(.; {\bf f}_{1},{\bf u}){\bf h}_{1}dxdt\\
\displaystyle-\int_{0}^{T}\!\!\!\int_{\Omega}\tilde{v}\frac{\p\ \Phi_{v}}{\p\ {\bf f}_{2}}(.; {\bf f}_{2},{\bf u}){\bf h}_{2}dxdt
\displaystyle-\int_{0}^{T}\!\!\!\int_{\Omega}\tilde{w}\frac{\p\ \Phi_{w}}{\p\ {\bf f}_{3}}(.; {\bf f}_{3},{\bf u}){\bf h}_{3}dxdt
\\
\displaystyle +\int_{0}^{T}\!\!\!\int_{\Omega}p\left[-\frac{\partial \tilde{u}}{\partial t}-div(D_{u}(.;{\bf u})\nabla \tilde{u})-\frac{\p\ \Phi_{u}}{\p\ u}(.; {\bf f}_{1},{\bf u})\tilde{u}
+(\tilde{H}_{u;u}-\frac{\p\ \mathfrak{T}_{u}}{\p\ u}(.;u))\nabla \tilde{u} \right]dxdt\\
\displaystyle-\int_{0}^{T}\!\!\!\int_{\Omega}p\left[-\tilde{H}_{v;u}\nabla \tilde{v} +\frac{\p\ \Phi_{v}}{\p\ u}(.; {\bf f}_{2},{\bf u}) \tilde{v} \right]dxdt
\displaystyle-\int_{0}^{T}\!\!\!\int_{\Omega}p\left[-\tilde{H}_{w;u}\nabla \tilde{w} +\frac{\p\ \Phi_{w}}{\p\ u}(.; {\bf f}_{3},{\bf u}) \tilde{w} \right]dxdt\\
\displaystyle+\int_{0}^{T}\!\!\!\int_{\Omega}div(D_{v,1}(.;{\bf u})\nabla \tilde{v})p dxdt+\int_{0}^{T}\!\!\!\int_{\Omega}div(D_{w,1}(.;{\bf u})\nabla \tilde{w})p dxdt\\ 

\displaystyle +\int_{0}^{T}\!\!\!\int_{\Omega}q\left[-\frac{\partial \tilde{v}}{\partial t}-div(D_{v}(.;{\bf u})\nabla \tilde{v})-\frac{\p\ \Phi_{v}}{\p\ v}(.; {\bf f}_{2},{\bf u})\tilde{v}
+(\tilde{H}_{v;v}-\frac{\p\ \mathfrak{T}_{v}}{\p\ v}(.;v))\nabla \tilde{v} \right]dxdt\\
\displaystyle-\int_{0}^{T}\!\!\!\int_{\Omega}q\left[-\tilde{H}_{u;v}\nabla \tilde{u} +\frac{\p\ \Phi_{u}}{\p\ v}(.; {\bf f}_{1},{\bf u}) \tilde{u} \right]dxdt
-\int_{0}^{T}\!\!\!\int_{\Omega}q\left[-\tilde{H}_{w;v}\nabla \tilde{w} +\frac{\p\ \Phi_{w}}{\p\ v}(.; {\bf f}_{3},{\bf u}) \tilde{w} \right]dxdt\\
\displaystyle+\int_{0}^{T}\!\!\!\int_{\Omega}div(D_{u,1}(.;{\bf u})\nabla \tilde{u})q dxdt+\int_{0}^{T}\!\!\!\int_{\Omega}div(D_{w,2}(.;{\bf u})\nabla \tilde{w})q dxdt\\

\displaystyle +\int_{0}^{T}\!\!\!\int_{\Omega}\rho\left[-\frac{\partial \tilde{w}}{\partial t}-div(D_{w}(.;{\bf u})\nabla \tilde{w})-\frac{\p\ \Phi_{w}}{\p\ w}(.; {\bf f}_{3},{\bf u})\tilde{w}
\displaystyle+(\tilde{H}_{w;w}-\frac{\p\ \mathfrak{T}_{w}}{\p\ w}(.;w))\nabla \tilde{w} \right]dxdt\\
\displaystyle-\int_{0}^{T}\!\!\!\int_{\Omega}\rho\left[-\tilde{H}_{u;w}\nabla \tilde{u} +\frac{\p\ \Phi_{u}}{\p\ w}(.; {\bf f}_{1},{\bf u}) \tilde{u} \right]dxdt
\displaystyle-\int_{0}^{T}\!\!\!\int_{\Omega}\rho\left[-\tilde{H}_{v;w}\nabla \tilde{v} +\frac{\p\ \Phi_{v}}{\p\ w}(.; {\bf f}_{2},{\bf u}) \tilde{v} \right]dxdt\\
\displaystyle+\int_{0}^{T}\!\!\!\int_{\Omega}div(D_{u,2}(.;{\bf u})\nabla \tilde{u})\rho dxdt+\int_{0}^{T}\!\!\!\int_{\Omega}div(D_{v,2}(.;{\bf u})\nabla \tilde{v})\rho dxdt 
=0. 
\end{array}
\end{equation}
In order to simplify \eqref{DJ}, according to \eqref{Int}, we suppose that $\tilde{\bf u}=(\tilde{u}, \tilde{v},\tilde{w})$ is the solution of the following adjoint system
\begin{equation}\label{FVM-PF-Adj-K}
\begin{array}{lcr}
\displaystyle -\frac{\partial \tilde{u}}{\partial t}-div(D_{u}(.;{\bf u})\nabla \tilde{u})-\frac{\p\ \Phi_{u}}{\p\ u}(.; {\bf f}_{1},{\bf u})\tilde{u}
+(\tilde{H}_{u;u}-\frac{\p\ \mathfrak{T}_{u}}{\p\ u}(.;u))\nabla \tilde{u}\\
\displaystyle+\tilde{H}_{v;u}\nabla \tilde{v} -\frac{\p\ \Phi_{v}}{\p\ u}(.;{\bf f}_{2},{\bf u}) \tilde{v} 
\displaystyle+\tilde{H}_{w;u}\nabla \tilde{w} -\frac{\p\ \Phi_{w}}{\p\ u}(.; {\bf f}_{3},{\bf u}) \tilde{w}\\
\displaystyle=-div(D_{v,1}(.;{\bf u})\nabla \tilde{v}) -div(D_{w,1}(.;{\bf u})\nabla \tilde{w})+\frac{\p\ {\cal Z}}{\p\ u}({\bf u},{\bf f}),\\ [0.2cm] 
\displaystyle -\frac{\partial \tilde{v}}{\partial t}-div(D_{v}(.;{\bf u})\nabla \tilde{v})-\frac{\p\ \Phi_{v}}{\p\ v}(.; {\bf f}_{2},{\bf u})\tilde{v}
+(\tilde{H}_{v;v}-\frac{\p\ \mathfrak{T}_{v}}{\p\ v}(.;v))\nabla \tilde{v}\\
\displaystyle+\tilde{H}_{u;v}\nabla \tilde{u} -\frac{\p\ \Phi_{u}}{\p\ v}(.; {\bf f}_{1},{\bf u}) \tilde{u}
+\tilde{H}_{w;v}\nabla \tilde{w} -\frac{\p\ \Phi_{w}}{\p\ v}(.; {\bf f}_{3},{\bf u}) \tilde{w}\\
\displaystyle=-div(D_{u,1}(.;{\bf u})\nabla \tilde{u})-div(D_{w,2}(.;{\bf u})\nabla \tilde{w})+\frac{\p\ {\cal Z}}{\p\ v}({\bf u},{\bf f}),\\ [0.2cm] 
\displaystyle -\frac{\partial \tilde{w}}{\partial t}-div(D_{w}(.;{\bf u})\nabla \tilde{w})-\frac{\p\ \Phi_{w}}{\p\ w}(.; {\bf f}_{3},{\bf u})\tilde{w}
\displaystyle+(\tilde{H}_{w;w}-\frac{\p\ \mathfrak{T}_{w}}{\p\ w}(.;w))\nabla \tilde{w}\\
\displaystyle+\tilde{H}_{u;w}\nabla \tilde{u} -\frac{\p\ \Phi_{u}}{\p\ w}(.; {\bf f}_{1},{\bf u}) \tilde{u} 
\displaystyle+\tilde{H}_{v;w}\nabla \tilde{v} -\frac{\p\ \Phi_{v}}{\p\ w}(.; {\bf f}_{2},{\bf u}) \tilde{v} \\
\displaystyle=-div(D_{u,2}(.;{\bf u})\nabla \tilde{u})-div(D_{v,2}(.;{\bf u})\nabla \tilde{v})+\frac{\p\ {\cal Z}}{\p\ w}({\bf u},{\bf f}),\\ [0.2cm] 
\text{under the final and boundary conditions \eqref{ACF}-\eqref{ACB}}.
\end{array}
\end{equation}
Then the previous relation becomes
\begin{equation}
\begin{array}{r}
\displaystyle \int_{\Omega}\frac{\p\ J_{0}}{\p\ {\bf u}}({\bf u}(T)).\wp (T)dx+\int_{0}^{T}\!\!\int_{\Omega}\frac{\p\ {\cal Z}}{\p\ {\bf u}}({\bf u},{\bf f}).\wp dxdt
=\int_{0}^{T}\!\!\!\int_{\Omega}\tilde{u}\frac{\p\ \Phi_{u}}{\p\ {\bf f}_{1}}(.; {\bf f}_{1},{\bf u}){\bf h}_{1}dxdt\\
\displaystyle\hspace{1.2cm}+\int_{0}^{T}\!\!\!\int_{\Omega}\tilde{u}\frac{\p\ \Phi_{v}}{\p\ {\bf f}_{2}}(.; {\bf f}_{2},{\bf u}){\bf h}_{2}dxdt
+\int_{0}^{T}\!\!\!\int_{\Omega}\tilde{u}\frac{\p\ \Phi_{w}}{\p\ {\bf f}_{3}}(.; {\bf f}_{3},{\bf u}){\bf h}_{3}dxdt
\end{array}
\end{equation}
and then 
\begin{equation}\label{DJR}
\begin{array}{lcr}
\displaystyle \frac{\p\ {\cal J}}{\p\ {\bf f}_{1}}({\bf f}).{\bf h}_{1}=
\int_{0}^{T}\!\!\!\int_{\Omega}(\frac{\p\ {\cal Z}}{\p\ {\bf f}_{1}}({\bf u},{\bf f})+\tilde{u}\frac{\p\ \Phi_{u}}{\p\ {\bf f}_{1}}(.; {\bf f}_{1},{\bf u})){\bf h}_{1}dxdt,\\
\displaystyle \frac{\p\ {\cal J}}{\p\ {\bf f}_{2}}({\bf f}).{\bf h}_{2}=
\int_{0}^{T}\!\!\!\int_{\Omega}(\frac{\p\ {\cal Z}}{\p\ {\bf f}_{2}}({\bf u},{\bf f})+\tilde{v}\frac{\p\ \Phi_{v}}{\p\ {\bf f}_{2}}(.; {\bf f}_{2},{\bf u})){\bf h}_{2}dxdt,\\
\displaystyle \frac{\p\ {\cal J}}{\p\ {\bf f}_{3}}({\bf f}).{\bf h}_{3}=
\int_{0}^{T}\!\!\!\int_{\Omega}(\frac{\p\ {\cal Z}}{\p\ {\bf f}_{3}}({\bf u},{\bf f})+\tilde{w}\frac{\p\ \Phi_{w}}{\p\ {\bf f}_{3}}(.; {\bf f}_{3},{\bf u})){\bf h}_{3}dxdt.
\end{array}
\end{equation}
Otherwise (in weak sense)
\begin{equation}\label{DJRW}
\begin{array}{lcr}
\displaystyle \frac{\p\ {\cal J}}{\p\ {\bf f}_{1}}({\bf f})=
\frac{\p\ {\cal Z}}{\p\ {\bf f}_{1}}({\bf u},{\bf f})+\tilde{u}\frac{\p\ \Phi_{u}}{\p\ {\bf f}_{1}}(.; {\bf f}_{1},{\bf u}),\\
\displaystyle \frac{\p\ {\cal J}}{\p\ {\bf f}_{2}}({\bf f})=
\frac{\p\ {\cal Z}}{\p\ {\bf f}_{2}}({\bf u},{\bf f})+\tilde{v}\frac{\p\ \Phi_{v}}{\p\ {\bf f}_{2}}(.; {\bf f}_{2},{\bf u}),\\
\displaystyle \frac{\p\ {\cal J}}{\p\ {\bf f}_{3}}({\bf f})=
\frac{\p\ {\cal Z}}{\p\ {\bf f}_{3}}({\bf u},{\bf f})+\tilde{w}\frac{\p\ \Phi_{w}}{\p\ {\bf f}_{3}}(.; {\bf f}_{3},{\bf u}).
\end{array}
\end{equation}
Since ${\bf f}^{*}$ is an optimal solution we have ($\forall {\bf f}=({\bf f}_{1}, {\bf f}_{2}, {\bf f}_{3})\in {\cal K}_{ad}$)
\begin{equation*}
\begin{array}{lcr}
\displaystyle \int_{0}^{T}\!\!\!\int_{\Omega}(\frac{\p\ {\cal Z}}{\p\ {\bf f}_{1}}({\bf u}^{*},{\bf f}^{*})+\tilde{u}^{*}\frac{\p\ \Phi_{u}}{\p\ {\bf f}_{1}}(.; {\bf f}^{*}_{1},{\bf u}^{*})).({\bf f}_{1}-{\bf f}^{*}_{1})dxdt\geq 0,\\
\displaystyle \int_{0}^{T}\!\!\!\int_{\Omega}(\frac{\p\ {\cal Z}}{\p\ {\bf f}_{2}}({\bf u}^{*},{\bf f}^{*})+\tilde{v}^{*}\frac{\p\ \Phi_{v}}{\p\ {\bf f}_{2}}(.; {\bf f}^{*}_{2},{\bf u}^{*})).({\bf f}_{2}-{\bf f}^{*}_{2})dxdt\geq 0,\\
\displaystyle \int_{0}^{T}\!\!\!\int_{\Omega}(\frac{\p\ {\cal Z}}{\p\ {\bf f}_{3}}({\bf u}^{*},{\bf f}^{*})+\tilde{w}^{*}\frac{\p\ \Phi_{w}}{\p\ {\bf f}_{3}}(.; {\bf f}^{*}_{3},{\bf u}^{*})).({\bf f}_{3}-{\bf f}^{*}_{3})dxdt\geq 0,
\end{array}
\end{equation*}
with ${\bf u}^{*}={\cal F}({\bf f}^{*})$ and $\tilde{\bf u}^{*}$ the solution of \eqref{FVM-PF-Adj-K} (corresponding to $({\bf f}^{*}, {\bf u}^{*})$).
This completes the proof. \hfill $\Box$
%
%%%%%%%%%%%%%%
%%%%%%%%%%%%%%%%%%%%%%%%
%
\begin{req}
1. The solution $\tilde{\bf u}$ of \eqref{FVM-PF-Adj-K} is said to be the costate, adjoint or dual solution corresponding to the primal solution ${\bf u}={\cal F}({\bf f})$ and will be denoted in the sequel by $\tilde{\bf u}=\tilde{\cal F}({\bf f},{\bf u})$.\\
2. If ${\cal K}_{ad}={\cal U}$ (i.e., optimal control without constraints) the optimality condition (\ref{COS}) become
\begin{equation}\label{COSSC}
\begin{array}{lcr}
\displaystyle 0=\frac{\p\ {\cal Z}}{\p\ {\bf f}_{1}}({\bf u}^{*},{\bf f}^{*})+\tilde{u}^{*}\frac{\p\ \Phi_{u}}{\p\ {\bf f}_{1}}(.; {\bf f}^{*}_{1},{\bf u}^{*}),\\
\displaystyle 0=\frac{\p\ {\cal Z}}{\p\ {\bf f}_{2}}({\bf u}^{*},{\bf f}^{*})+\tilde{v}^{*}\frac{\p\ \Phi_{v}}{\p\ {\bf f}_{2}}(.; {\bf f}^{*}_{2},{\bf u}^{*}),\\
\displaystyle 0=\frac{\p\ {\cal Z}}{\p\ {\bf f}_{3}}({\bf u}^{*},{\bf f}^{*})+\tilde{w}^{*}\frac{\p\ \Phi_{w}}{\p\ {\bf f}_{3}}(.; {\bf f}^{*}_{3},{\bf u}^{*}). 
\end{array}
\end{equation}
\end{req}
\begin{req}
In order to solve the nonlinear control problem numerically, by using the adjoint variables, we can combine, for example, the obtained optimal necessary conditions and the gradient-iterative algorithm or Newton-iterative algorithm (see e.g., \cite{B0}). This resolution requires, at each iteration of optimization algorithm, the numerical resolution of direct problem and its corresponding adjoint problem. The direct, sensitivity and adjoint problems, which are coupled systems of anisotropic convection-diffusion chemotaxis-type equations, can be solved by using a multiple-relaxation-time lattice Boltzmann method.\hfill $\Box$
\end{req}
We end this section by a description of a gradient algorithm to solve Problem \eqref{EPS}. The method is formulated in terms of continuous variables which are independent of a specific numerical discretization. 
\subsection{Gradient-iterative algorithm and continuous adjoint models}
For a given  function $J_{0}$ and ${\cal Z}$ (and then measurement data $\mathfrak{m}_{data}$ and reference data $\mathfrak{r}_{data}$) and initial states $(u_{0},v_{0},w_{0})$, we present a gradient-iterative algorithm where the descent direction is calculated by using adjoint variables, particularly by choosing an admissible step size. 
The gradient algorithm  is given bellows (for $k=0,\ldots$, (iteration index) we denote by ${\bf f}^{(k)}=({\bf f}^{(k)}_{1}, {\bf f}^{(k)}_{2}, {\bf f}^{(k)}_{3})$ the approximation of control variables at the $kth$ iteration of the algorithm):\\
{\bf Step 1:} Initialization: ${\bf f}^{(0)}=({\bf f}^{(0)}_{1}, {\bf f}^{(0)}_{2}, {\bf f}^{(0)}_{3})$ (given initial guess).\\
{\bf Step 2:} Solve problem (\ref{eq:ReacDiff})-(\ref{E2}) with source term ${\bf f}^{(k)}$, gives ${\bf u}^{(k)}=(u^{(k)},v^{(k)},w^{(k)})={\cal F}({\bf f}^{(k)})$.\\% solution of
{\bf Step 3:} Solve problem (\ref{FVM-PF-Adj}) (based on ${\bf u}^{(k)}$),  gives adjoint state $\tilde{\cal F}({\bf f}^{(k)}, {\bf u}^{(k)})=(\uT\ ^{(k)},\vT\ ^{(k)},\wT\ ^{(k)})$.\\
{\bf Step 4:} Gradient of ${\cal J}$, $E_{k}={\cal J}'({\bf f}^{(k)})=(\frac{\p\ {\cal J}}{\p\ {\bf f}_{1}}({\bf f}^{(k)}),\frac{\p\ {\cal J}}{\p\ {\bf f}_{2}}({\bf f}^{(k)}),\frac{\p\ {\cal J}}{\p\ {\bf f}_{3}}({\bf f}^{(k)}))$ at point ${\bf f}^{(k)}$ is given by
\begin{equation}\label{GJA}
\begin{array}{lcr}
\displaystyle \frac{\p\ {\cal J}}{\p\ {\bf f}_{1}}({\bf f}^{(k)})= \frac{\p\ {\cal Z}}{\p\ {\bf f}_{1}}({\bf u}^{(k)},{\bf f}^{(k)})+\tilde{u}^{(k)}\frac{\p\ \Phi_{u}}{\p\ {\bf f}_{1}}({\bf f}_{1}^{(k)},{\bf u}^{(k)}),\\
\displaystyle \frac{\p\ {\cal J}}{\p\ {\bf f}_{2}}({\bf f}^{(k)})= \frac{\p\ {\cal Z}}{\p\ {\bf f}_{2}}({\bf u}^{(k)},{\bf f}^{(k)})+\tilde{v}^{(k)}\frac{\p\ \Phi_{v}}{\p\ {\bf f}_{2}}({\bf f}_{2}^{(k)},{\bf u}^{(k)}),\\
\displaystyle \frac{\p\ {\cal J}}{\p\ {\bf f}_{3}}({\bf f}^{(k)})= \frac{\p\ {\cal Z}}{\p\ {\bf f}_{3}}({\bf u}^{(k)},{\bf f}^{(k)})+\tilde{w}^{(k)}\frac{\p\ \Phi_{w}}{\p\ {\bf f}_{3}}({\bf f}_{3}^{(k)},{\bf u}^{(k)}),\\
\end{array}
\end{equation}
{\bf Step 5:} Determine ${\bf f}^{(k+1)}$ : $\displaystyle {\bf f}^{(k+1)}={\bf f}^{(k)}-\zeta_{k} E_{k}$,\\ 
where $\zeta_{k}$ with $0<m_{0}\leq \zeta_{k}\leq M_{0}$ are the sequences of step lengths.\\
{\bf Step 6:} If gradient of ${\cal J}$ is sufficiently small: end; else set  $k:=k+1$ and goto {\bf Step 2}.

Optimal Solution: ${\bf f}^{*}:={\bf f}^{(k)}$ and $(u^{*},v^{*},w^{*}):=(u^{(k)},v^{(k)},w^{(k)})$.
\begin{req}
1. In order to obtain an algorithm which is numerically efficient, the best choice of $\zeta_{k}$ will be the result of a line minimization algorithm. Otherwise, at each iteration step $k$ of the previous algorithm, we solve the one-dimensional optimization problem of parameter $\zeta_{k}$:
\begin{equation}\label{LM}
\begin{array}{lcr}
\displaystyle \zeta_{k}=\min_{\lambda>0}{\cal J}({\bf f}^{(k)}-\lambda E_{k}).
\end{array}
\end{equation}
To derive an approximation for $\zeta_{k}$ we can use a purely heuristic approach, e.g., by taking $\zeta_{k}=min(1,\pr\ E_{k}\pr\ ^{-1}_{\infty})$ or by using the linearization of ${\cal F}({\bf f}^{(k)}-\lambda E_{k})$ at ${\bf f}^{(k)}$ by : $\displaystyle{\cal F}({\bf f}^{(k)}-\lambda E_{k})\approx {\cal F}({\bf f}^{(k)})-\lambda {\wp}^{(k)}={\bf u}^{(k)}-\lambda {\wp}^{(k)}$, where ${\wp}^{(k)}={\cal F}'({\bf f}^{(k)}).E_{k}=(p^{(k)},q^{(k)},\rho^{(k)})$ is the solution of the sensitivity problem (\ref{E1d}). For example in the case of quadratic cost function ${\cal J}$, according to the previous approximation, we can approximate (\ref{LM}) by %\label{PE-CBI-PF-K}
\begin{equation}\label{LMA}
\begin{array}{lcr}
\displaystyle \zeta_{k}=\min_{\lambda>0}P(\lambda),\\
\end{array}
\end{equation}
where $P$ is a polynomial function of degree $2$ (since ${\cal J}$ is quadratic), then problem (\ref{LMA}) can be solved exactly. 
Consequently, we obtain explicitly the value of $\zeta_{k}$.\\
2. We can write similarly the conjugate gradient algorithm combined with the Wolfe-Powell line search procedure for computing admissible step-sizes along
the descent direction. The advantage of this method, compared to the gradient method, is that it performs
a soft reset whenever the GC search direction yields no significant progress. In general, the method has the following form:
\begin{equation*}
\begin{array}{lcr}
D_{k}=\left\{
\begin{array}{lcr}
-G_{k}~~for~~k=0,\\
-G_{k}+r_{k-1}D_{k-1}~~for~~k\geq1,
\end{array}
\right.
\\
z_{k+1}=z_{k}-\zeta_{k}D_{k}
\end{array}
\end{equation*}
where $G_{k}$ denotes gradient of the functional to be optimized at point $z_{k}$, $\zeta_{k}$ is a step length obtained by a line search, $D_{k}$ is the search direction and $r_{k}$ is a constant. Several varieties of this method differ in the way of selecting $r_{k}$. Some well-known formula for $r_{k}$ are given by Fletcher-Reeves (FR), Polak-Ribi\`ere (PR), Hestenes-Stiefel (HS) 
and Dai-Yuan (DY). \hfill $\Box$
\end{req}
\section{Adjoint-based numerical optimization approaches}\label{Section Opti}
For given functions $J_{0}$ and ${\cal Z}$, and initial states ${\bf u}_{0}=(u_{0},v_{0},w_{0})$, in numerical treatments of nonlinear optimal control problems, the direct system, sensitivity system, adjoint system and objective functional must be discretized (reduction of the infinite-dimensional dynamics to finite-dimensional problems). The discretized formulation for direct, sensitivity and adjoint problems can be performed by combining  the classical Galerkin and finite element methods, to the variational formulations, for the space discretization and Euler method for the time discretization, or by using lattice Boltzmann methods. It is clear that the amount of computation is crucial in solution of control problems, because the computation of discrete gradient of objective functional by discrete adjoint methods requires one forward solve of discrete state system and one backward solve of adjoint systems in which the state trajectory is an input. In order to reduce the implementation cost to evaluate discrete gradients, we propose and develop an adjoint multiple-relaxation-time coupled lattice Boltzmann method which is efficient, stable and easy to implement and to parallelize. In objective functional, the integrals with respect to time can be approximated by the composition trapezoidal rules.

Based on work presented in \cite{BIJB,YSC}, in next section we develop the multiple-relaxation-time lattice Boltzmann method for \eqref{eq:ReacDiff}-\eqref{E2} and show that the method can recover \eqref{eq:ReacDiff} through multi-scaling Chapman-Enskog procedure.
\subsection{Multiple-relaxation-time LBM for general coupled anisotropic convection-diffusion chemotaxis-type equations}\label{part:LBEreacdiff}
In order to develop multiple-relaxation-time LBM for \eqref{eq:ReacDiff} and show that this method can recover \eqref{eq:ReacDiff}, let us introduce the evolution equation of the LBM (see \cite{Boltz1872}) for solving a reaction-diffusion equation, with the macroscopic variable $\Theta$,
\begin{equation}\label{EqBoltz}
\displaystyle \frac{\partial \varphi(\mathbf{x},t;\mathbf{e})}{\partial t}+\mathbf{e}\cdot\nabla \varphi(\mathbf{x},t;\mathbf{e})={\cal J}(\mathbf{x},t;\mathbf{e}), 
\end{equation}
with
\begin{equation} \label{eq:IntphitoTheta}
\begin{array}{lcr}
\displaystyle \Theta(\mathbf{x},t)=\int \varphi(\mathbf{x},t;\mathbf{e}) d\mathbf{e},
\end{array}
\end{equation}
where $\displaystyle {\cal J}(\mathbf{x},t;\mathbf{e})=Q^{col}(\varphi(\mathbf{x},t;\mathbf{e})-\varphi^{eq}(\mathbf{x},t;\mathbf{e}))+P(\mathbf{x},t;\mathbf{e})$, $\varphi(\mathbf{x},t;\mathbf{e})$ is  the distribution function of the single particle moving with velocity $\mathbf{e}$ at position $\mathbf{x}$ and time $t$,  $P$ is the distribution type function of particle of macroscopic external force $\Phi$ moving with velocity $\mathbf{e}$, $Q^{col}$ is a Bhatnagar-Gross-Krook (BGK) type collision operator (see \cite{BGK}) and $\varphi^{eq}$ is the  Maxwell-Boltzmann equilibrium distribution function. 

LBM leads us to approximate \eqref{EqBoltz} to recover reaction-diffusion equation \eqref{eq:ReacDiff} with multi-scaling Chapman-Enskog expansion (see e.g., \cite{ChapEnsk}). In method of Chapman-Enskog, the distribution function $\varphi$ is expanded as : $\varphi=\varphi^{(0)}+\varepsilon \varphi^{(1)}+\varepsilon^{2}\varphi^{(2)}+....=\sum_{k=0}^{\infty}\varepsilon^{k}\varphi^{(k)}$, which can be regarded as a power series in a small variable $\varepsilon$ or an expression that keeps track of the relative orders of magnitude of the  different terms through scaling parameter $\varepsilon$ (the so-called Knudsen number). The distribution functions $\varphi^{(0)}$, $\varphi^{(1)}$ and $\varphi^{(2)}$ represent the zero, first and second  approximation to the distribution function $\varphi$, and so on.

The numerical solution of Boltzmann equation  (\ref{EqBoltz}) requires to discretize the configuration spaces, velocities and time. For that, we discretize $\overline{\Omega}\times(0,T)$ in time and space. For a given space discretization parameters $h$ and lattice time step size $\tau={T}/{N}$, $N\in \N{*}$, the domain $\overline{\Omega}=\Omega\cup \p\ \Omega$ is approximated by a uniform, regular grid with spacing $h$ (which is commonly called the lattice), and the interval $[0,T]$ is split by the equidistant partitioning using the following points $0=t_{0}\leq t_{1}\leq\cdots\leq t_{N}=T$
with $t_{n}=n\tau$, for $n=0,\ldots,N$. The steps $h$ and $\tau$ are chosen sufficiently small to guarantee both the time accuracy and convergence of solution.

Finally, we define the discrete time interval $\daleth_{\tau}=\{t\in [0,T]: t=n\tau,~~n\in \N{}\}$, streaming lattice speed $c=h/\tau$ and lattice sound speed $C_s$ by $3C_s^2=c^2$.  

We note $(\mathbf{e}_{i})_{i=0,q-1}$ the discrete set of $q$ admissible particle velocities (the magnitude of each $\mathbf{e}_{i}=(e_{i}^{(k)})_{k=1,m}\in \R{m}$ depends on speed $c$) and we assume that for each node $\textbf{x}$ of lattice ${\cal L}_{h}$, and each velocity $\mathbf{e}_{i}$, the point $\mathbf{x}+\mathbf{e}_{i}\tau$ is also a node of the lattice ${\cal L}_{h}$. For the latter, we use the so-called $DmQq$ LBM scheme (i.e., m-dimensional and $q$ velocity vectors).
First, we consider the finite discrete-velocity system of Boltzmann equation with finite discrete velocity $\mathbf{e}_{i}$, for $i=0,q-1$ (by discretizing the velocity space by  discrete set of microscopic velocities), for system \eqref{eq:ReacDiff} (for $\Theta=u, v$ or $w$)
\begin{equation}\label{EqBoltz-D}
\begin{array}{lcr}
\displaystyle \frac{\partial \varphi_{\Theta,i}(\mathbf{x},t)}{\partial t}+\mathbf{e}_{i}\cdot\nabla \varphi_{\Theta,i}(\mathbf{x},t)={\cal J}_{\Theta,i}(\mathbf{x},t),\\
\end{array}
\end{equation}
where, for each particle on lattice, we associate discrete functions, for $i=0,q-1$ and $\Theta=u, v$ or $w$, $\varphi_{\Theta,i}(\mathbf{x},t)= \varphi_{\Theta}(\mathbf{x},t;\mathbf{e}_{i})$, $\varphi^{eq}_{\Theta,i}(\mathbf{x},t)= 
\varphi_{\Theta}^{eq}(\mathbf{x},t;\mathbf{e}_{i})$ and  $\displaystyle {\cal J}_{\Theta, i}=\sum_{j}\hat{\Lambda}^{(\Theta)}_{ij}(\varphi_{\Theta, j}-\varphi^{eq}_{\Theta, j})+\Phi_{\Theta,i}+{\cal S}_{\Theta,i}$, where $\Phi_{\Theta,i}$ (resp. ${\cal S}_{\Theta,i}$) is the discrete operator of $\Phi_{\Theta}$ (resp. the discrete operator corresponding to cross-diffusion operators), which describes the probability of streaming in one particular direction and $\hat{\bf \Lambda}_{\Theta}$, defined by $\hat{\bf \Lambda}_{\Theta}(\mathbf{x},t)=(\hat{\Lambda}^{(\Theta)}_{ij}(\mathbf{x},t))_{i,j=0,q-1}$, is the collision matrix (see e.g., \cite{DHG}). 

The aim of the so-called multiple-relaxation-time LBM scheme is to precise the distribution function of  particle $(\varphi_{\Theta, i}(\mathbf{x},t))_{i=0,q-1}$ for $\mathbf{x}\in {\cal L}_{h}$ and $t\in \daleth_{\tau}$ by solving the discretization of Boltzmann equation \eqref{EqBoltz-D} in two steps (for each time): collision and streaming processes. Then, by using a second order time integration scheme to approximate \eqref{EqBoltz-D} (by limiting physical space to a lattice and velocity space to discrete set of microscopic velocities), we can obtain 

\begin{equation}\label{eq:LBE}
\begin{array}{lcr}
\displaystyle  \varphi_{\Theta,i}(\mathbf{x}+\mathbf{e}_i\tau,t+\tau)=\displaystyle\varphi_{\Theta,i}(\mathbf{x},t)-\sum_{j}({\cal M}^{-1}\Pi_{\Theta} {\cal M})_{ij}\left( \varphi_{\Theta,j}(\mathbf{x},t)-\varphi_{\Theta,j}^ {eq}(\mathbf{x},t)\right)  \\
\hspace{2cm}\displaystyle +\tau {\cal S}_{\Theta,i}(\mathbf{x},t)+  \tau\Phi_{\Theta,i}(\mathbf{x},t)+ \frac{\tau^2}{2}(\kappa_{t}\frac{\partial}{\partial t}+\kappa_{x}\mathbf{e}_{i}\cdot\nabla)\Phi_{\Theta,i}(\mathbf{x},t),
\end{array}
\end{equation}
where the discrete cross-diffusion operators are defined by
\begin{equation}\label{S}
\begin{array}{lcr}
\displaystyle {\cal S}_{u,i}=\frac{1}{C^{2}_{s}}\omega_{i}\left(\mathbf{e}_i.(A_{u,1}\nabla \mathfrak{h}_{u,1}(v))+\mathbf{e}_i.(A_{u,2}\nabla \mathfrak{h}_{u,2}(w))\right),\\
\displaystyle {\cal S}_{v,i}=\frac{1}{C^{2}_{s}}\omega_{i}\left(\mathbf{e}_i.(A_{v,1}\nabla \mathfrak{h}_{v,1}(u))+\mathbf{e}_i.(A_{v,2}\nabla \mathfrak{h}_{v,2}(w))\right),\\
\displaystyle {\cal S}_{w,i}=\frac{1}{C^{2}_{s}}\omega_{i}\left(\mathbf{e}_i.(A_{w,1}\nabla \mathfrak{h}_{w,1}(u))+\mathbf{e}_i.(A_{w,2}\nabla \mathfrak{h}_{w,2}(v))\right),
\end{array}
\end{equation}
and the discrete source terms are defined by ($\Theta=u, v$ or $w$)
\begin{equation}\label{PhiS}
\Phi_{\Theta,i}=\omega_{i}\Phi_{\Theta}\left(1 + \frac{{\vec{V}_{\Theta}}.\mathbf{e}_{i}}{C^{2}_{s}}\right),
\end{equation}
where the differential tensors $A_{\Theta,k}$ and ${\vec{V}_{\Theta}}$ (for $\Theta=u, v$ or $w$, $k=1,2$)  will be determined below and $\mathfrak{h}_{\Theta,k}$ is the identity operator, $\kappa_{x}$ and $\kappa_{t}$ are two parameters to be determined later on and $\Pi_{\Theta}=\tau {\cal M}\hat{\bf \Lambda}_{\Theta}{\cal M}^{-1}$ is the relaxation matrix in velocity space.  ${\cal M}$ is a transformation matrix (constant) which is constructed from discrete velocities via the Gram-Schmidt orthogonalization procedure. Then the collision process is carried out in moment space and ${\cal M}$ can be used to project $\varphi_i$ and $\varphi^{eq}_{\Theta,i}$ in discrete velocity space onto macroscopic variables in moment space as 
\begin{equation}\label{MF}
{\bf m}_{\Theta}={\cal M}{\bf F}_{\Theta}=(m_{\Theta,i})_{i=0,q-1},~~{\bf m}_{\Theta}^{eq}={\cal M}{\bf F}^{eq}_{\Theta}=(m^{eq}_{\Theta,i})_{i=0,q-1}  
\end{equation}
where ${\bf F}_{\Theta}=(\varphi_{\Theta, i})_{i=0,q-1}$ and ${\bf F}^{eq}_{\Theta}=(\varphi^{eq}_{\Theta, i})_{i=0,q-1}$.
Therefore, the collision process in moment space can be expressed as
\begin{equation}\label{MC}
{\bf m}^{col}_{\Theta}={\bf m}_{\Theta}-\Pi_{\Theta}({\bf m}_{\Theta}-{\bf m}^{eq}_{\Theta}).
\end{equation}
Introduce also the following family of matrix as (for $k=1,m$)
\begin{equation}\label{RFM} 
{\cal E}_{k}={\cal M}~\mathrm{diag}(e^{(k)}_{0},e^{(k)}_{1},...,e^{(k)}_{q-1}){\cal M}^{-1} \text{~and~} {\cal E}=({\cal E}_{k})_{k=1,m}.
\end{equation}
Then the equilibrium distribution functions $\varphi^{eq}_{\Theta,i}(\mathbf{x},t)$, for $\Theta=u,v$ or $w$, of the LBM for system \eqref{eq:ReacDiff}, associated with discrete velocity $\mathbf{e}_{i}$ at position $\mathbf{x}$ and time $t$, are defined as (for $i=0,q-1$)
\begin{equation}\label{phieq}
\varphi^{eq}_{\Theta,i}=\omega_{i}\left(\Theta+ \frac{\mathfrak{T}_{\Theta}.\mathbf{e}_{i}}{C^{2}_{s}}+\frac{(\mathfrak{C}_{\Theta}+d_{\Theta}C^{2}_{s}z{\bf I}_{d}-C^{2}_{s}\Theta {\bf I}_{d}):\mho_{i}}{2C^{4}_{s}}\right)
\end{equation}
where ${\bf I}_{d}$ is the unit matrix, $d_{\Theta}>0$ is a parameter related to the diffusion tensor $D_{\Theta}$, the matrix operator $(\vec{\nu}\otimes \vec{\mu})$, for $\vec{\nu}=(\nu_{i})_{i=1,m}$ and $\vec{\mu}=(\mu_{i})_{i=1,m}$, is defined by $(\vec{\nu}\otimes \vec{\mu})_{ij}=\nu_{i}\mu_{j}$, for $i,j=1,m$, the contraction product between two matrices $A=(a_{ij})_{i,j=1,m}$ and $B=(b_{ij})_{i,j=1,m}$ is defined by  $A:B=\sum_{i,j}a_{ij}b_{ij}$ and $\omega_{i}$ is a weight coefficient in direction $\mathbf{e}_{i}$. Finally the second-order differential tensor 
$\mathfrak{C}_{\Theta}$ is such that (for all $({\bf x},t)$ and $\Theta$)
\begin{equation}\label{tensorC}
\frac{\p\ \mathfrak{C}_{\Theta}}{\p\ \Theta}({\bf x},t;\Theta)=\frac{\p\ \mathfrak{T}_{\Theta}}{\p\ \Theta}\otimes \frac{\p\ \mathfrak{T}_{\Theta}}{\p\ \Theta}({\bf x},t;\Theta)
\end{equation}
and the matrix $\mho_{i}$ is given by 
\begin{equation}\label{amalg}
\mho_{i}=\mathbf{e}_i\otimes\mathbf{e}_i -C_s^2{\bf I}_{d}.
\end{equation}
\begin{req}
For different lattice models $DmQq$, the velocities $\mathbf{e}_{i}$ and weight coefficients $\omega_{i}$, for $i=0,\dots,q$  can be defined as (for example)
\begin{description}
\item[\underline{For $D2Q9$:}]
%\begin{eqnarray*}
$\displaystyle \displaystyle\mathbf{e}_0=c\begin{pmatrix}0\\0\end{pmatrix},~\mathbf{e}_1=c\begin{pmatrix}1\\0\end{pmatrix},~\mathbf{e}_2=c\begin{pmatrix}0\\1\end{pmatrix},~\mathbf{e}_3=-\mathbf{e}_{1}
~\mathbf{e}_4=-\mathbf{e}_{2},\displaystyle \displaystyle\mathbf{e}_5=c\begin{pmatrix}1\\1\end{pmatrix},~\mathbf{e}_6=c\begin{pmatrix}-1\\1\end{pmatrix},$\\$~\mathbf{e}_7=-\mathbf{e}_{5},~\mathbf{e}_8=-\mathbf{e}_6$
%\end{eqnarray*}
and 
\begin{equation}\label{sumWi}
\begin{array}{lcr}
\sum\limits_{i=0}^8w_i=1,~~ w_0=4/9,~~
w_i=1/9 \text{~(for~} i=1,\dots,4) \text{~and~} w_i=1/36 \text{~(for~} i=5,\dots,8).
\end{array}
\end{equation}
\item[\underline{For $D3Q19$:}]
%\begin{eqnarray*}
$\displaystyle \displaystyle\mathbf{e}_0=c\begin{pmatrix}0\\0\\0\end{pmatrix},~\mathbf{e}_1=c\begin{pmatrix}1\\0\\0\end{pmatrix},~\mathbf{e}_2=-\mathbf{e}_1,
~\mathbf{e}_3=c\begin{pmatrix}0\\1\\0\end{pmatrix},~\mathbf{e}_4=-\mathbf{e}_3,~\displaystyle \mathbf{e}_5=c\begin{pmatrix}0\\0\\1\end{pmatrix},~\mathbf{e}_6=-\mathbf{e}_5,$\\
$\mathbf{e}_7=c\begin{pmatrix}1\\1\\0\end{pmatrix},~\mathbf{e}_8=-\mathbf{e}_7,~\mathbf{e}_9=c\begin{pmatrix}1\\-1\\0\end{pmatrix},~\mathbf{e}_{10}=-\mathbf{e}_9,~\mathbf{e}_{11}=c\begin{pmatrix}1\\0\\1\end{pmatrix},~\mathbf{e}_{12}=-\mathbf{e}_{11},$\\ 
$~\mathbf{e}_{13}=c\begin{pmatrix}-1\\0\\1\end{pmatrix},~\mathbf{e}_{14}=-\mathbf{e}_{13},~\mathbf{e}_{15}=c\begin{pmatrix}0\\1\\1\end{pmatrix},~\mathbf{e}_{16}=-\mathbf{e}_{15},~\mathbf{e}_{17}=c\begin{pmatrix}0\\1\\-1\end{pmatrix},~\mathbf{e}_{18}=-\mathbf{e}_{17}$
%\end{eqnarray*}
and 
\begin{equation}\label{sumWi3}
\begin{array}{lcr}
\sum\limits_{i=0}^{18} \omega_i=1,~~ \omega_0=1/3,~~
\omega_i=1/18 \text{~(for~} i=1,\dots,6) \text{~and~} \omega_i=1/36 \text{~(for~} i=7,\dots,18).
\end{array}
\end{equation}
\hfill $\Box$
\end{description}
\end{req}
Finally, in order to derive precisely the macroscopic variable $\Theta$ in terms of microscopic particle distribution function  governed by lattice Boltzmann equation \eqref{eq:LBE},  the following conditions also need to be satisfied 
\begin{equation}\label{TF0}
\begin{array}{lcr}
\displaystyle \Theta=\sum_{i}\varphi_{\Theta,i}=\sum_{i}\varphi^{eq}_{\Theta,i},~~   \mathfrak{T}_{\Theta}=\sum_{i}{\bf e}_{i}\varphi^{eq}_{\Theta,i} =\sum_{i}{\bf e}_{i}\varphi_{\Theta,i}\text{~~(the momentum)~},\\
\displaystyle \sum\limits_{i} \mathbf{e}_i\otimes\mathbf{e}_i\varphi_{\Theta,i}^{eq} = (\mathfrak{C}_{\Theta}+ d_{\Theta}C_s^2 \Theta) \text{~~(the stress tensor)~},\\
\displaystyle \Phi_{\Theta}=\sum_{i}\Phi_{\Theta,i},~~\displaystyle\sum_{i}e_{i}\Phi_{\Theta,i}=\vec{V}_{\Theta} \Phi_{\Theta}, ~~\text{and}\\
\sum_{i}{\cal S}_{\Theta,i}=0,~~\sum_{i}\mathbf{e}_i{\cal S}_{\Theta,i}=A_{\Theta,1}\nabla \mathfrak{h}_{\Theta,1}+A_{\Theta,2}\nabla \mathfrak{h}_{\Theta,2}.\\
\end{array}
\end{equation}
\begin{req}
 During streaming and collision processes, in order to satisfy boundary conditions, boundary nodes need special treatments on distribution functions, which are essential to stability and accuracy of the method. We can use extrapolation methods (see e.g., \cite{GuoBoundary2002}) and modified bounce-back methods (see e.g.,\cite{ZhangBounceBack12,ImprovedBounceBack}).\hfill $\Box$
\end{req}
In the following proposition, we prove formally that governing macroscopic equation \eqref{eq:ReacDiff}  is recovered correctly from the developed multiple-relaxation-time LBM model. 

\begin{pps}
The LB equations \eqref{eq:LBE} with \eqref{S} and \eqref{PhiS} approximate solutions to continuous macroscopic system \eqref{eq:ReacDiff} via Chapman Enskog asymptotic expansion with   errors of order $O(\varepsilon^{3})$, where $\varepsilon$ is the so-called non-uniformity parameter $\varepsilon$ (which can be interpreted as the Knudsen number that controls the order of the approximation).
\end{pps}
{\bf Proof:} 
To derive the macroscopic equation from developed multiple-relaxation-time lattice BGK model, the Chapman-Enskog expansion is applied under the assumption of small Knudsen
number $\varepsilon$ (which can be a ratio between a characteristic length $L$ and a particular mean free path $l$) to determine equilibrium functions, source terms and link functions. Each
function is decomposed around its equilibrium state with different scale of perturbations. Then, the difference between associated LBE and Taylor expansion of the distribution function
lead us to recover the governing macroscopic equation.
For every  distribution functions $\varphi_{\Theta, i}$, the Chapman-Enskog expansion is applied with the small parameter $\varepsilon$ as:
 \begin{eqnarray}
\displaystyle \varphi_{\Theta,i} &=&\varphi^{(0)}_{\Theta, i}+\varepsilon \varphi^{(1)}_{\Theta, i}+ \varepsilon^2 \varphi^{(2)}_{\Theta, i}+... ,\label{Ansatz:phi}\\
\displaystyle \frac{\partial}{\partial t} &=&\displaystyle \varepsilon \frac{\partial}{\partial t_1} +\varepsilon^2\frac{\partial}{\partial t_2}, \label{Ansatz:t_inf}\\
\displaystyle  \grad &=&\displaystyle \varepsilon \grad_1~,~\text{div}_{1}=\grad_{1}\cdot, \label{Ansatz:grad_inf}
\end{eqnarray}
where $\varphi^{(k)}_{\Theta, i}$ are assumed to be sufficiently smooth functions in $({\bf x},t)$ and independent of $\varepsilon$, $t_{i}$ is the time scale (for $i=1,2$), ${\bf x}_{1}$ is the space scale, $\nabla_{1}=\nabla_{{\bf x}_{1}}$ is the gradient with respect to ${\bf x}_{1}$ and $\varphi_{\Theta, i}^{(0)}$ is associated to unperturbed state. Correspondingly, macroscopic quantity $\Theta$ is expanded  as
$$  \displaystyle \Theta =\Theta^{(0)}+\varepsilon \Theta^{(1)} + \varepsilon^2 \Theta^{(2)}+..., \text{~with~} \Theta^{(n)}=\sum_{i}\varphi^{(n)}_{\Theta, i}.$$
We assume also that every term $\Phi_{\Theta, i}$ and ${\cal S}_{\Theta,i}$ take the form
\begin{equation}\label{Ansatz:Phi}
\begin{array}{lcr}
\displaystyle \Phi_{\Theta, i} =\varepsilon \Phi_{\Theta, i}^{(1)}+\varepsilon^{2} \Phi_{\Theta, i}^{(2)},~~{\cal S}_{\Theta,i}=\varepsilon{\cal S}_{\Theta,i}^{(1)}+\varepsilon^{2}{\cal S}_{\Theta,i}^{(2)}\\
\text{and we denote by $\Phi_{\Theta}^{(n)}$ and ${\cal S}_{\Theta}^{(n)}$ the terms}\\
\displaystyle \Phi_{\Theta}^{(n)}=\sum\limits_{i} \Phi_{\Theta, i}^{(n)},~~ {\cal S}_{\Theta}^{(n)}=\sum\limits_{i} {\cal S}_{\Theta,i}^{(n)},~~\text{~for $n=1,2$}.
\end{array}
\end{equation}
The previous expansions are inserted  into LBE \eqref{eq:LBE} and Taylor-series expansion\footnote{Taylor-series expansion
with truncation to the $p$-order derivative terms and remainder term, $p\geq 2$ (depending on the regularity of $\varphi_{\Theta, i}$) can be applied.} is applied to $\varphi_{\Theta, i}(\mathbf{x}+\mathbf{e}_i\tau,t+\tau)$ ($\varphi_i$ moving with velocity $\mathbf{e}_i$) at point $({\bf x},t)$ as
\begin{equation}\label{eq:TaylorExp}
\begin{array}{rl}
\displaystyle \varphi_{\Theta, i}(\mathbf{x}+\mathbf{e}_i\tau,t+\tau)=\varphi_{\Theta, i}(\mathbf{x},t)+ \sum\limits_{k\geq 1}\frac{\tau^k}{k!}\left(\frac{\partial}{\partial t} + \mathbf{e}_i\cdot\nabla \right)^k\varphi_{\Theta, i}(\mathbf{x},t),
\end{array}
\end{equation}
where
\begin{equation*}
\displaystyle \left(\frac{\partial}{\partial t} + \mathbf{e}_i\cdot\nabla \right)^k= \sum\limits_{n=0}^k \begin{pmatrix}  k \\ n\end{pmatrix} \frac{\partial^{k-n}}{\partial t^{k-n}}\left(\mathbf{e}_i\cdot\nabla \right)^n,
\end{equation*}
with 
$\displaystyle (\mathbf{e}_i\cdot\nabla)^n= (\mathbf{e}_i\cdot\nabla) \left(\mathbf{e}_i\cdot\nabla\right)^{n-1}.$

So, by using the expressions \eqref{Ansatz:phi}-\eqref{Ansatz:Phi}, the expansion \eqref{eq:TaylorExp} and LBE \eqref{eq:LBE} become ($\forall i=0,\dots,q-1$)
\begin{equation}\label{eq:ChapmEnskTaylor}
\begin{array}{lcr}
\displaystyle \varphi_{\Theta, i}(\mathbf{x}+\mathbf{e}_i\tau,t+\tau)=\left(\varphi_{\Theta, i}^{(0)}(\mathbf{x},t)+\varepsilon\varphi_{\Theta, i}^{(1)}(\mathbf{x},t)+...\right)\\
\hspace{0.5cm}\displaystyle +\sum\limits_{k\geq 1}\frac{\tau^k}{k!}\left(\varepsilon \frac{\partial}{\partial t_1} +\varepsilon^2\frac{\partial}{\partial t_2}+ \mathbf{e}_i\cdot(\varepsilon \grad_1) \right)^k
\left(\varphi_{\Theta, i}^{(0)}(\mathbf{x},t)+\varepsilon\varphi_{\Theta, i}^{(1)}(\mathbf{x},t)+...\right) 
\end{array}
\end{equation}
and 
\begin{eqnarray}
\displaystyle \varphi_{\Theta, i}(\mathbf{x}+\mathbf{e}_i\tau,t+\tau)&=&\left(\varphi_{\Theta, i}^{(0)}(\mathbf{x},t)+\varepsilon\varphi_{\Theta, i}^{(1)}(\mathbf{x},t)+...\right) \nonumber \\
& &\displaystyle\hspace{-3.cm}-\sum_{j}({\cal M}^{-1}\Pi_{\Theta} {\cal M})_{ij}\left(\varphi_{\Theta, j}^{(0)}(\mathbf{x},t)+\varepsilon\varphi_{\Theta, j}^{(1)}(\mathbf{x},t)+...\right)+\sum_{j}({\cal M}^{-1}\Pi_{\Theta} {\cal M})_{ij}\varphi_{\Theta, j}^ {eq}(\mathbf{x},t) \nonumber \\
& &\displaystyle\hspace{-3.cm}+ \tau\left(\varepsilon \Phi_{\Theta, i}^{(1)}(\mathbf{x},t)+\varepsilon^{2} \Phi_{\Theta, i}^{(2)}(\mathbf{x},t)\right) 
+ \tau\left(\varepsilon{\cal S}_{\Theta,i}^{(1)}(\mathbf{x},t)+\varepsilon^{2}{\cal S}_{\Theta,i}^{(2)}(\mathbf{x},t)\right)\nonumber\\
 & &\displaystyle \hspace{-3.cm}+\frac{\tau^2}{2}\left( \varepsilon \kappa_{t}\frac{\partial}{\partial t_1} +\varepsilon^2\kappa_{t}\frac{\partial}{\partial t_2}+\varepsilon\kappa_{x}\mathbf{e}_i\cdot \grad_1\right) \displaystyle \left(\varepsilon \Phi_{\Theta, i}^{(1)}(\mathbf{x},t)+\varepsilon^{2} \Phi_{\Theta, i}^{(2)}(\mathbf{x},t)\right). \label{eq:ChapmEnskLBE}
\end{eqnarray}
By calculating the difference between relations \eqref{eq:ChapmEnskTaylor} and \eqref{eq:ChapmEnskLBE}, we can obtain
\begin{eqnarray}
0&=&\sum_{j}({\cal M}^{-1}\Pi_{\Theta} {\cal M})_{ij}\left(\varphi_{\Theta, j}^{(0)}(\mathbf{x},t)-\varphi_{\Theta, j}^ {eq}(\mathbf{x},t)\right)\nonumber\\
& &+\sum_{j}({\cal M}^{-1}\Pi_{\Theta} {\cal M})_{ij}\left(\varepsilon\varphi_{\Theta, j}^{(1)}(\mathbf{x},t)+\varepsilon^{2}\varphi_{\Theta, j}^{(2)}(\mathbf{x},t)+...\right)\nonumber\\
& &+\sum\limits_{k\geq 1}\frac{\tau^k}{k!}\left(\varepsilon \frac{\partial}{\partial t_1} +\varepsilon^2\frac{\partial}{\partial t_2}+ \mathbf{e}_i\cdot(\varepsilon \grad_1) \right)^k
\left(\varphi_{\Theta, i}^{(0)}(\mathbf{x},t)+\varepsilon\varphi_{\Theta, i}^{(1)}(\mathbf{x},t)+...\right)\nonumber\\
& &- \tau\left(\varepsilon \Phi_{\Theta, i}^{(1)}(\mathbf{x},t)+\varepsilon^{2} \Phi_{\Theta, i}^{(2)}(\mathbf{x},t)\right)-\tau\left(\varepsilon{\cal S}_{\Theta,i}^{(1)}(\mathbf{x},t)+\varepsilon^{2}{\cal S}_{\Theta,i}^{(2)}(\mathbf{x},t)\right) \nonumber\\
& &\displaystyle  - \displaystyle \frac{\tau^2}{2}\left( \varepsilon \kappa_{t}\frac{\partial}{\partial t_1} +\varepsilon^2\kappa_{t}\frac{\partial}{\partial t_2}+\varepsilon\kappa_{x}\mathbf{e}_i\cdot \grad_1\right)\left(\varepsilon \Phi_{\Theta, i}^{(1)}(\mathbf{x},t)+\varepsilon^{2} \Phi_{\Theta, i}^{(2)}(\mathbf{x},t)\right), \label{eq:DiffLBETaylor}
\end{eqnarray}
which can be rewritten as follows (according to different scales of $\varepsilon$)
\begin{equation}
\displaystyle 0=\Psi_{\Theta, i}^{(0)}(\mathbf{x},t) + \varepsilon \Psi_{\Theta, i}^{(1)}(\mathbf{x},t) + \varepsilon^2 \Psi_{\Theta, i}^{(2)}(\mathbf{x},t) +\dots,
\end{equation}
where probability density functions $\Psi_{\Theta, i}^{(n)}$ regroup all terms of \eqref{eq:DiffLBETaylor} in order of $\varepsilon^n$, $n=0,...,$ respectively (which are partial differential equations for $n\neq 0$)
\begin{eqnarray}
\hspace{-1.5cm} \displaystyle \Psi_{\Theta, i}^{(0)}&=&\displaystyle \sum_{j}({\cal M}^{-1}\Pi_{\Theta} {\cal M})_{ij}\left(\varphi_{\Theta, j}^{(0)}- \varphi_{\Theta, j}^{eq} \right)\!=\!0~(\text{then $\varphi_{\Theta, i}^{(0)}=\varphi_{\Theta, i}^{eq}, i=0,q-1$)}, \label{eq:Psi0} \\
\hspace{-1.5cm} \displaystyle  \Psi_{\Theta, i}^{(1)}&=&\tau\left(\frac{\partial}{\partial t_1}+\mathbf{e}_i\cdot\nabla_1 \right)\varphi_{\Theta, i}^{(0)}+\sum_{j}({\cal M}^{-1}\Pi_{\Theta} {\cal M})_{ij}\varphi_{\Theta, j}^{(1)}
-\tau \Phi_{\Theta, i}^{(1)}-\tau {\cal S}_{\Theta,i}^{(1)}=0,\label{eq:Psi1}\\
\hspace{-1.5cm} \displaystyle  \Psi_{\Theta, i}^{(2)}&=&\left(\frac{\tau^2}{2}\left(\frac{\partial}{\partial t_1}+\mathbf{e}_i\cdot\nabla_1 \right)^2 + \tau\frac{\partial}{\partial t_2}\right)\varphi_{\Theta, i}^{(0)}+\tau\left(\frac{\partial}{\partial t_1}+\mathbf{e}_i\cdot\nabla_1 \right)\varphi_{\Theta, i}^{(1)}\nonumber\\
\displaystyle &&\displaystyle +\sum_{j}({\cal M}^{-1}\Pi_{\Theta} {\cal M})_{ij}\varphi_{\Theta, j}^{(2)}-\tau \Phi_i^{(2)}-\tau {\cal S}_{\Theta,i}^{(2)}%\nonumber\\
%\displaystyle &&\displaystyle 
-\frac{\tau^2}{2}(\kappa_{t}\frac{\partial}{\partial t_1}+\kappa_{x}\mathbf{e}_i\cdot \grad_1)\Phi_{\Theta, i}^{(1)}=0,\label{eq:Psi2}\\
\dots \nonumber
\end{eqnarray}
To treat second order terms, we rewrite $\displaystyle \frac{\tau^2}{2}\left(\frac{\partial}{\partial t_1}+\mathbf{e}_i\cdot\nabla_1 \right)^2 \varphi_{\Theta, i}^{(0)}$ with using \eqref{eq:Psi1} as:
\begin{eqnarray*}
\displaystyle \frac{\tau^2}{2}\left(\frac{\partial}{\partial t_1}+\mathbf{e}_i\cdot\nabla_1 \right)^2 \varphi_{\Theta, i}^{(0)} &=& 
\frac{\tau}{2} \left(\frac{\partial}{\partial t_1}+\mathbf{e}_i\cdot\nabla_1 \right)\left(-\sum_{j}({\cal M}^{-1}\Pi_{\Theta} {\cal M})_{ij}\varphi_{\Theta,j}^{(1)} + \tau\Phi_{\Theta, i}^{(1)}
+\tau{\cal S}_{\Theta,i}^{(1)} \right),\nonumber\\
&=&-\tau\left(\frac{\partial}{\partial t_1}+\mathbf{e}_i\cdot\nabla_1 \right)\left(\frac{1}{2}\sum_{j}({\cal M}^{-1}\Pi_{\Theta} {\cal M})_{ij}\varphi_{\Theta, j}^{(1)}\right) \nonumber\\
& &+\frac{\tau^{2}}{2} \left(\frac{\partial}{\partial t_1}+\mathbf{e}_i\cdot\nabla_1 \right)\Phi_{\Theta, i}^{(1)}+\frac{\tau^{2}}{2} \left(\frac{\partial}{\partial t_1}+\mathbf{e}_i\cdot\nabla_1 \right){\cal S}_{\Theta,i}^{(1)} ,%\nonumber\\
\end{eqnarray*}
and then $\Psi_i^{(2)}$ becomes
\begin{eqnarray}
\displaystyle  \Psi_{\Theta, i}^{(2)}&=&\tau\left(\frac{\partial}{\partial t_1}+\mathbf{e}_i\cdot\nabla_1 \right)\left(\varphi_{\Theta, i}^{(1)}-\frac{1}{2}\sum_{j}({\cal M}^{-1}\Pi_{\Theta} {\cal M})_{ij}\varphi_{\Theta, j}^{(1)}\right) 
+\sum_{j}({\cal M}^{-1}\Pi_{\Theta} {\cal M})_{ij}\varphi_{\Theta, j}^{(2)} \nonumber\\
\displaystyle &&\displaystyle + \tau\frac{\partial}{\partial t_2}\varphi_{\Theta, i}^{(0)} -  \tau \Phi_{\Theta, i}^{(2)}-\tau {\cal S}_{\Theta,i}^{(2)}+\frac{\tau^{2}}{2} \left(\frac{\partial}{\partial t_1}+\mathbf{e}_i\cdot\nabla_1 \right){\cal S}_{\Theta,i}^{(1)}\nonumber\\
\displaystyle &&-\frac{\tau^2}{2}((\kappa_{t}-1)\frac{\partial}{\partial t_1}+(\kappa_{x}-1)\mathbf{e}_i\cdot \grad_1)\Phi_{\Theta, i}^{(1)}=0.\label{eq:Psi2B}
\end{eqnarray}
%%%%%
If we rewrite the previous equations \eqref{eq:Psi0}, \eqref{eq:Psi1} and \eqref{eq:Psi2B} in a vector form and we multiply on both sides of them, we can deduce by {\it retaining terms up to $O(\varepsilon^{3})$}
\begin{equation}\label{MPV}
\begin{array}{lcr}
%\begin{eqnarray}
\displaystyle \varepsilon^{0}:\hspace{1cm}{\bf m}^{(0)}_{\Theta}={\bf m}^{eq}_{\Theta},\\
\displaystyle  \varepsilon^{1}:\hspace{1cm} 0 =\tau\left({\bf I}_{d}\frac{\partial}{\partial t_1}+{\cal E}\nabla_1 \right){\bf m}^{(0)}_{\Theta}+\Pi_{\Theta} {\bf m}^{(1)}_{\Theta}-\tau {\cal H}^{(1)}_{\Theta}-\tau {\cal L}^{(1)}_{\Theta},\\
\displaystyle  \varepsilon^{2}:\hspace{1cm}0=\tau\left({\bf I}_{d}\frac{\partial}{\partial t_1}+{\cal E}\nabla_1 \right)\left({\cal I}_{d}-\frac{1}{2}\Pi_{\Theta} \right){\bf m}_{\Theta}^{(1)} 
+\Pi_{\Theta} {\bf m}^{(2)}_{\Theta} + \tau\frac{\partial}{\partial t_2}{\bf m}^{(0)}_{\Theta}\\
\displaystyle \hspace{2cm} -  \tau {\cal H}^{(2)}_{\Theta} \displaystyle -\frac{\tau^2}{2}\left((\kappa_{t}-1){\bf I}_{d}\frac{\partial}{\partial t_1}+(\kappa_{x}-1){\cal E}\grad_1\right){\cal H}^{(1)}_{\Theta}\\
\displaystyle \hspace{2cm}-\tau {\cal L}^{(2)}_{\Theta} +\frac{\tau^2}{2}\left({\bf I}_{d}\frac{\partial}{\partial t_1}+{\cal E}\grad_1\right){\cal L}^{(1)}_{\Theta},
\end{array}
\end{equation}
where $({\cal M}^{-1}{\cal H}^{(n)}_{\Theta})_{i}=\Phi_{\Theta, i}^{(n)}$ and  $({\cal M}^{-1}{\cal L}^{(n)}_{\Theta})_{i}={\cal S}_{\Theta,i}^{(n)}$, for $n=1,...$, and ${\cal E}\nabla_1=\sum_{k=1,m}{\cal E}_{k}\nabla_{1,k}$.

\begin{figure}[htbp]
\begin{center}
\hspace{0.5cm}\includegraphics[scale=0.32]{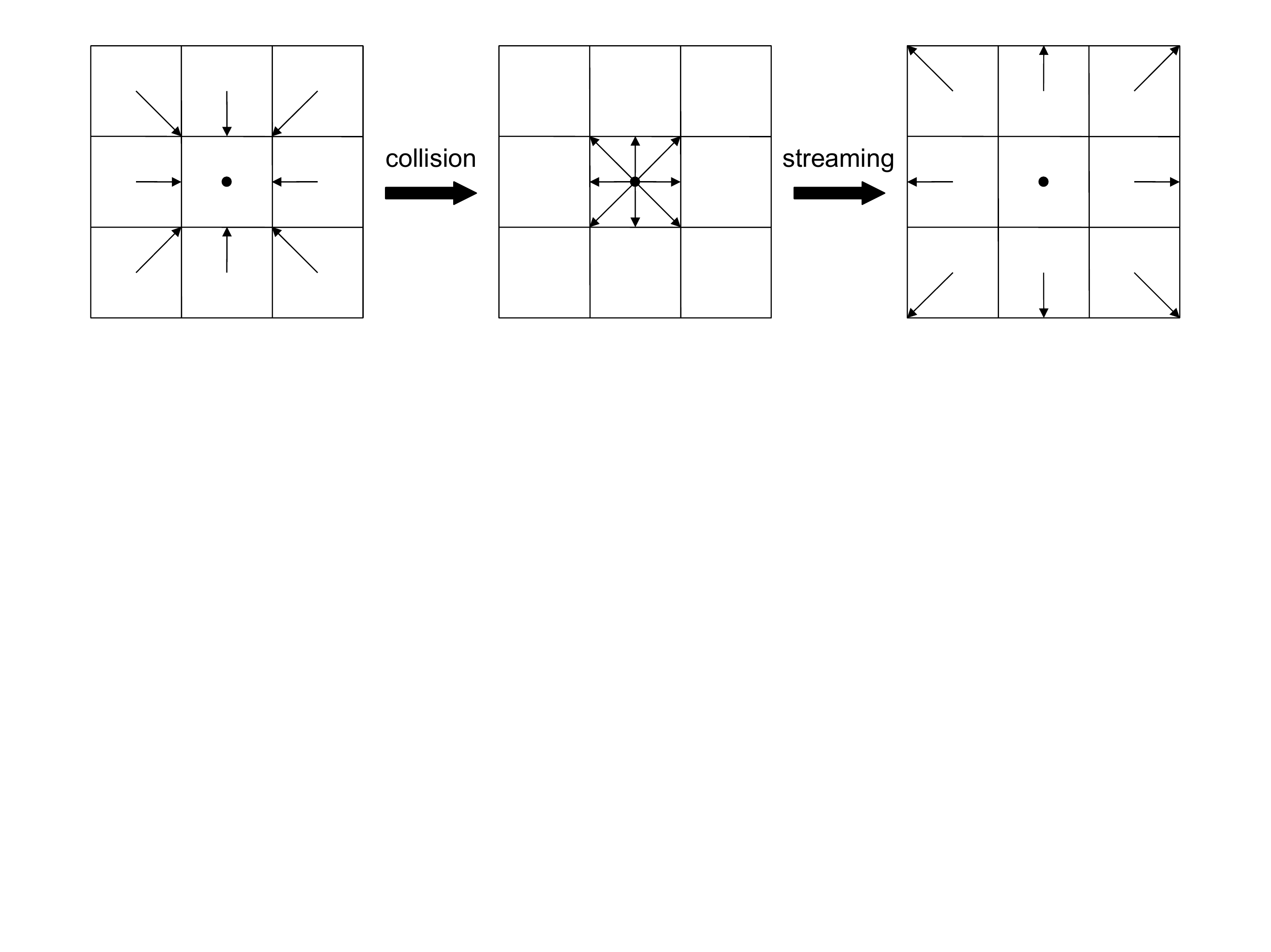}
\vspace{-4.5cm}
\caption{\small Collision-Streaming process}\label{fig:1}
\end{center}
\end{figure}

In order to facilitate the presentation we consider, in the rest of this section, the popular two-dimensional $D2Q9$ lattice model (see Fig.\ref{fig:1}), which involves 9 velocity vectors as an example (the arguments developed below can be extended to a 3-dimensional $D3Q19$ lattice without any substantial difficulties,  because this last has essentially  the same procedure as that for the $D2Q9$ one). 
According to the expression of microscopic velocities $\mathbf{e}_{i}$, we prove easily that 
\begin{equation*}
\sum\limits_{i=0}^8\mathbf{e}_i = \sum\limits_{i=0}^8\omega_i\mathbf{e}_i = \begin{pmatrix}0\\0\end{pmatrix},~
 \sum\limits_{i=0}^{8}\omega_i\mathbf{e}_i\otimes \mathbf{e}_i =C_s^2 \mathbf{I}_d ~\text{and}~\sum\limits_{i=0}^{8}\omega_i\mho_{i}=0.
\end{equation*}
The transportation matrix ${\cal M}$ can be written as ${\cal M}=E_{0}{\cal M}_{0}$ with ${\cal M}_{0}$ given by (see e.g. \cite{LL})
{\footnotesize{
\begin{eqnarray*}
{\cal M}_{0}=\begin{pmatrix}
1& 1& 1& 1& 1& 1& 1& 1& 1\\
-4& -1& -1& -1& -1& 2& 2& 2& 2\\
4& -2& -2& -2& -2& 1& 1& 1& 1\\
0& 1& 0& -1& 0& 1& -1& -1& 1\\
0& -2& 0& 2& 0& 1& -1& -1& 1\\
0& 0& 1& 0& -1& 1& 1& -1& -1\\
0& 0& -2& 0& 2& 1& 1& -1& -1\\
0& 1& -1& 1& -1& 0& 0& 0& 0\\
0& 0& 0& 0& 0& 1& -1& 1& -1
\end{pmatrix},
\end{eqnarray*}}}
and $E_{0}=\mathrm{diag}(1, c^{2}, c^{4}, c, c^{3}, c, c^{3}, c^{2}, c^{2})$ a diagonal matrix. We can deduce that ${\bf m}^{(eq)}_{\Theta}=(m^{(eq)}_{\Theta,i})_{i=0,8}$ is given by
\begin{equation*}
\begin{array}{lcr}
m^{(eq)}_{\Theta,0}=\Theta,~~ m^{(eq)}_{\Theta,1}=-4c^{2}\Theta+3(B^{(\Theta)}_{11}+B^{(\Theta)}_{22}), ~~ m^{(eq)}_{\Theta,2}=3c^{4}\Theta-3c^{2}(B^{(\Theta)}_{11}+B^{(\Theta)}_{22}),\\[0.2cm]
m^{(eq)}_{3}=\mathfrak{T}_{\Theta,1}(\Theta),~~ m^{(eq)}_{\Theta,4}=-c^{2}\mathfrak{T}_{\Theta,1}(\Theta),~~ m^{(eq)}_{\Theta,5}=\mathfrak{T}_{\Theta,2}(\Theta),\\[0.2cm]
m^{(eq)}_{6}=-c^{2}\mathfrak{T}_{\Theta,2}(\Theta), ~~ m^{(eq)}_{\Theta,7}=B^{(\Theta)}_{11}-B^{(\Theta)}_{22}, ~~ m^{(eq)}_{\Theta,8}=B^{(\Theta)}_{12},
\end{array}
\end{equation*}
where 
\begin{equation}\label{BC11}
{\bf B}_{\Theta}=(B^{(\Theta)}_{ij})_{i,j=1,2}=\sum_{k}(\mathbf{e}_k \otimes \mathbf{e}_k)\varphi_{\Theta,k}^{(eq)}=\left(\mathfrak{C}_{\Theta}+d_{\Theta}C_{s}^{2}\Theta{\bf I}_{d}\right). 
\end{equation}
In order to handle anisotropic diffusion problems, the relaxation matrix $\Pi_{\Theta}$ in the moment space is not a diagonal matrix, as in the isotropic case, but can be expressed as a block-diagonal invertible matrices :
$\Pi_{\Theta}=\mathrm{diag}(\Lambda_{\Theta,1}, \Lambda_{\Theta,2}, \Lambda_{\Theta,3})$, with
\begin{equation*}
\begin{array}{lcr}
\Lambda_{\Theta,1}=\mathrm{diag}(s_{\Theta,0},s_{\Theta,1},s_{\Theta,2}),~~\Lambda_{\Theta,2}=\begin{pmatrix}
k_{\Theta,11}& 0& k_{\Theta,12}\\
0&  s_{\Theta,3}& 0\\
k_{\Theta,21}&  0& k_{\Theta,22}
\end{pmatrix},
~~\Lambda_{\Theta,3}=\mathrm{diag}(s_{\Theta,4},s_{\Theta,5},s_{\Theta,6})
\end{array},
\end{equation*}
where the diagonal elements $s_{\Theta,j}$ ($j=0,2$), $k_{\Theta,11}$, $s_{\Theta,3}$   and $k_{\Theta,22}$, $s_{\Theta,j}$ ($j=4,6$) are the relaxation parameters corresponding to jth moment ${m}_{j}$, $j=0,8$, respectively, and the terms 
$k_{\Theta,12}$ and $k_{\Theta,21}$ correspond to the rotation of the principal axis of anisotropic diffusion (see e.g.\cite{YN}). 
\begin{req}
We prove easily that $E_{0}\Pi_{\Theta}=\Pi_{\Theta} E_{0}$ and then ${\cal M}^{-1}\Pi_{\Theta} {\cal M}={\cal M}_{0}^{-1}\Pi_{\Theta} {\cal M}_{0}$. Then from \eqref{eq:LBE} 
\begin{eqnarray}\label{eq:LBEM}
\begin{array}{lcr}
\displaystyle  \varphi_{\Theta, i}(\mathbf{x}+\mathbf{e}_i\tau,t+\tau)=&\displaystyle\varphi_{\Theta, i}(\mathbf{x},t)-\sum_{j}({\cal M}_{0}^{-1}\Pi_{\Theta} {\cal M}_{0})_{ij}\left( \varphi_{\Theta, j}(\mathbf{x},t)-\varphi_{\Theta, j}^ {eq}(\mathbf{x},t)\right)  \nonumber\\
\displaystyle &\displaystyle\displaystyle +\tau {\cal S}_{\Theta,i}(\mathbf{x},t)+  \tau\Phi_{\Theta, i}(\mathbf{x},t)+ \frac{\tau^2}{2}(\kappa_{t}\frac{\partial}{\partial t}+\kappa_{x}\mathbf{e}_{i}\cdot\nabla)\Phi_{\Theta, i}(\mathbf{x},t).
\end{array}
\end{eqnarray}
\end{req}
Based on the second relation of \eqref{MPV} we can deduce that\footnote{we consider only the first, fourth and sixth components which are useful in order to cover the state equation.} (according to \eqref{TF0})
\begin{equation}\label{C146}
\begin{array}{lcr}
\displaystyle \frac{\p\ \Theta}{\p\ t_{1}}+div_{1}(\mathfrak{T}_{\Theta})=\Phi^{(1)}_{\Theta},\\
\displaystyle \tau\left(\frac{\p\ \mathfrak{T}_{\Theta}}{\p\ t_{1}}+div_{1}({\bf B}_{\Theta})-\Phi_{\Theta}^{(1)}\vec{V}_{\Theta}-(A_{\Theta,1}\nabla \mathfrak{h}_{\Theta,1}+A_{\Theta,2}\nabla \mathfrak{h}_{\Theta,2})\right)=-{\bf K}_{\Theta}\tilde{\bf m}_{\Theta}^{(1)},
\end{array}%+{\cal S}^{(1)}_{z}
\end{equation}
where the invertible matrix ${\bf K}_{\Theta}$ is given by
\begin{eqnarray*}
{\bf K}_{\Theta}=\begin{pmatrix}
k_{\Theta,11}& k_{\Theta,12}\\
k_{\Theta,21}& k_{\Theta,22}
\end{pmatrix},~~\tilde{\bf m}^{(1)}_{\Theta}=(m_{\Theta,3}^{(1)},m_{\Theta,5}^{(1)})^{t}
\end{eqnarray*}
and on the third relation of \eqref{MPV} (i.e. the second order relation in $\varepsilon$) we can deduce (by considering only the first which is useful in order to cover the state equation)
\begin{equation}\label{PC1}
\begin{array}{lcr}
\displaystyle \frac{\p\ \Theta}{\p\ t_{2}}+div_{1}\left(({\bf I}_{d}-\frac{1}{2}{\bf K}_{\Theta})\tilde{\bf m}^{(1)}_{\Theta}\right)=\frac{\tau}{2}(\kappa_{t}-1)\frac{\p\ \Phi_{\Theta}^{(1)}}{\p\ t_{1}}+\Phi_{\Theta}^{(2)}\\
\hspace{2cm}\displaystyle+ \frac{\tau}{2}(\kappa_{x}-1)div_{1}(\vec{V}_{\Theta}\Phi_{\Theta}^{(1)})-\frac{\tau}{2}div_{1}(A_{\Theta,1}\nabla \mathfrak{h}_{\Theta,1}+A_{\Theta,2}\nabla \mathfrak{h}_{\Theta,2}).
\end{array}
\end{equation}
Then, from the second equation of \eqref{C146}, the equation \eqref{PC1} becomes
\begin{equation*}
\begin{array}{lcr}
\displaystyle \frac{\p\ \Theta}{\p\ t_{2}}-\tau div_{1}(({\bf K}_{\Theta}^{-1}-\frac{1}{2}{\bf I}_{d})\left(\frac{\p\ \mathfrak{T}_{\Theta}}{\p\ t_{1}}+div_{1}({\bf B}_{\Theta})-\Phi_{\Theta}^{(1)}\vec{V}_{\Theta}-(A_{\Theta,1}\nabla \mathfrak{h}_{\Theta,1}+A_{\Theta,2}\nabla \mathfrak{h}_{\Theta,2})\right )\\
%div_{1}\left(({\bf K}^{-1}-\frac{1}{2}{\bf I}_{d})(\frac{\p\ (\Theta \vec{\omega})}{\p\ t_{1}}+div_{1}({\bf B})-\vec{V}\Phi^{(1)})\right)\\
\hspace{2cm}\displaystyle=\Phi_{\Theta}^{(2)}+\frac{\tau}{2}(\kappa_{t}-1)\frac{\p\ \Phi_{\Theta}^{(1)} }{\p\ t_{1}}+ \frac{\tau}{2}(\kappa_{x}-1)div_{1}(\vec{V}_{\Theta}\Phi_{\Theta}^{(1)})
-\frac{\tau}{2}div_{1}(A_{\Theta,1}\nabla \mathfrak{h}_{\Theta,1}+A_{\Theta,2}\nabla \mathfrak{h}_{\Theta,2})
\end{array}
\end{equation*}
 and then 
\begin{equation}\label{PC12}
\begin{array}{lcr}
\displaystyle \frac{\p\ \Theta}{\p\ t_{2}}-\tau div_{1}\left(({\bf K}_{\Theta}^{-1}-\frac{1}{2}{\bf I}_{d})(\frac{\p\ \mathfrak{T}_{\Theta}}{\p\ t_{1}}+div_{1}({\bf B}_{\Theta}))\right)\\
\hspace{1cm}\displaystyle=\Phi_{\Theta}^{(2)}+\frac{\tau}{2}(\kappa_{t}-1)\frac{\p\ \Phi^{(1)}_{\Theta}}{\p\ t_{1}}- \tau div_{1}\left(({\bf K}_{\Theta}^{-1}-\frac{\kappa_{x}}{2}{\bf I}_{d})\vec{V}_{\Theta}\Phi^{(1)}_{\Theta}\right)\\
\hspace{1cm}\displaystyle-\tau div_{1}\left({\bf K}_{\Theta}^{-1}A_{\Theta,1}\nabla \mathfrak{h}_{\Theta,1}+{\bf K}_{\Theta}^{-1}A_{\Theta,2}\nabla \mathfrak{h}_{\Theta,2})\right).
\end{array}
\end{equation}
Calculate now the value $\frac{\p\ \mathfrak{T}_{\Theta}}{\p\ t_{1}}+div_{1}({\bf B}_{\Theta})$. From \eqref{BC11}, we can deduce
%Since (from \eqref{BC11}) 
\begin{equation*}
\begin{array}{lcr}
\displaystyle\frac{\p\ }{\p\ t_{1}}(\mathfrak{T}_{\Theta}(\mathbf{x},t;\Theta(\mathbf{x},t))+div_{1}({\bf B}_{\Theta}(\mathbf{x},t;\Theta(\mathbf{x},t)))=\frac{\p\ \mathfrak{T}_{\Theta}}{\p\ \Theta}(\mathbf{x},t;\Theta)\frac{\p\ \Theta}{\p\ t_{1}}+\frac{\p\ \mathfrak{T}_{\Theta}}{\p\ t_{1}}(\mathbf{x},t;\Theta)\\
\displaystyle \hspace{2cm}+\frac{\p\ \mathfrak{C}_{\Theta}}{\p\ \Theta}(\mathbf{x},t;\Theta)\nabla_{1}\Theta+div_{1}(\mathfrak{C}_{\Theta})(\mathbf{x},t;\Theta)+div_{1}(d_{\Theta}C_{s}^{2}\Theta {\bf I}_{d}).
\end{array}
\end{equation*}
According to  \eqref{tensorC}, we can deduce
\begin{equation*}
\begin{array}{lcr}
\displaystyle\frac{\p\ }{\p\ t_{1}}(\mathfrak{T}_{\Theta}(\mathbf{x},t;\Theta(\mathbf{x},t))+div_{1}({\bf B}_{\Theta}(\mathbf{x},t;\Theta(\mathbf{x},t)))=\frac{\p\ \mathfrak{T}_{\Theta}}{\p\ \Theta}(\mathbf{x},t;\Theta)\frac{\p\ \Theta}{\p\ t_{1}}+\frac{\p\ \mathfrak{T}_{\Theta}}{\p\ t_{1}}(\mathbf{x},t;\Theta)\\
\displaystyle \hspace{2cm}+((\frac{\p\ \mathfrak{T}_{\Theta}}{\p\ \Theta}(\mathbf{x},t;\Theta))\nabla_{1} \Theta)\frac{\p\ \mathfrak{T}_{\Theta}}{\p\ \Theta}(\mathbf{x},t;\Theta)
+div_{1}(\mathfrak{C}_{\Theta})(\mathbf{x},t;\Theta)+div_{1}(d_{\Theta}C_{s}^{2}\Theta {\bf I}_{d}).\\
\end{array}
\end{equation*}
Since 
$div_{1}(\mathfrak{T}_{\Theta}(\mathbf{x},t;\Theta(\mathbf{x},t))=div_{1}(\mathfrak{T}_{\Theta})(\mathbf{x},t;\Theta)+(\frac{\p\ \mathfrak{T}_{\Theta}}{\p\ \Theta}(\mathbf{x},t;\Theta))\nabla_{1} \Theta,$
we can deduce that 
\begin{equation*}
\begin{array}{lcr}
\displaystyle\frac{\p\ }{\p\ t_{1}}(\mathfrak{T}_{\Theta}(\mathbf{x},t;\Theta(\mathbf{x},t))+div_{1}({\bf B}_{\Theta}(\mathbf{x},t;\Theta(\mathbf{x},t)))\\
\displaystyle \hspace{2cm}=\frac{\p\ \mathfrak{T}_{\Theta}}{\p\ \Theta}(\mathbf{x},t;\Theta)(\frac{\p\ \Theta}{\p\ t_{1}}+div_{1}(\mathfrak{T}_{\Theta}(\mathbf{x},t;\Theta(\mathbf{x},t)))
+d_{\Theta}C_{s}^{2}\nabla_{1}\Theta \\
\displaystyle \hspace{2cm}+\frac{\p\ \mathfrak{T}_{\Theta}}{\p\ t_{1}}(\mathbf{x},t;\Theta)+(div_{1}(\mathfrak{C}_{\Theta})(\mathbf{x},t;\Theta)-div_{1}(\mathfrak{T}_{\Theta})(\mathbf{x},t;\Theta)\frac{\p\ \mathfrak{T}_{\Theta}}{\p\ \Theta}(\mathbf{x},t;\Theta)).
%+C_{s}^{2}\nabla_{1}\Theta+div_{1}(\vec{\omega}\Theta)\vec{\omega}+\frac{d\vec{\omega}}{dt_{1}},
\end{array}
\end{equation*}
Denotes now by $DV$ the following operator
\begin{equation}\label{deviation}
DV(\mathbf{x},t;\Theta)=\frac{\p\ \mathfrak{T}_{\Theta}}{\p\ t_{1}}(\mathbf{x},t;\Theta)+(div_{1}(\mathfrak{C}_{\Theta})(\mathbf{x},t;\Theta)-div_{1}(\mathfrak{T}_{\Theta})(\mathbf{x},t;\Theta)\frac{\p\ \mathfrak{T}_{\Theta}}{\p\ \Theta}(\mathbf{x},t;\Theta)).
\end{equation}
Then (from the first equation of \eqref{C146})
\begin{equation}\label{BC12}
\begin{array}{lcr}
\displaystyle\frac{\p\ }{\p\ t_{1}}(\mathfrak{T}_{\Theta}(\mathbf{x},t;\Theta(\mathbf{x},t))+div_{1}({\bf B}_{\Theta}(\mathbf{x},t;\Theta(\mathbf{x},t)))\\
\displaystyle \hspace{2cm}=\frac{\p\ \mathfrak{T}_{\Theta}}{\p\ \Theta}(\mathbf{x},t;\Theta(\mathbf{x},t))\Phi_{\Theta}^{(1)}+d_{\Theta}C_{s}^{2}\nabla_{1}\Theta+DV(\mathbf{x},t;\Theta).
\end{array}
\end{equation}
Consequently, equation \eqref{PC12} becomes
\begin{equation}\label{PC13}
\begin{array}{lcr}
\displaystyle \frac{\p\ \Theta}{\p\ t_{2}}-div_{1}\left(d_{\Theta}C_{s}^{2}\tau({\bf K}_{\Theta}^{-1}-\frac{1}{2}{\bf I}_{d})\nabla_{1}\Theta\right)
\displaystyle \hspace{0cm}=\Phi_{\Theta}^{(2)}+\frac{\tau}{2}(\kappa_{t}-1)\frac{\p\ \Phi_{\Theta}^{(1)}}{\p\ t_{1}}\\
\displaystyle\hspace{1.2cm}+ \tau div_{1}\left(\left(-({\bf K}_{\Theta}^{-1}-\frac{\kappa_{x}}{2}{\bf I}_{d})\vec{V}_{\Theta}+({\bf K}_{\Theta}^{-1}-\frac{1}{2}{\bf I}_{d})\frac{\p\ \mathfrak{T}_{\Theta}}{\p\ \Theta}(\mathbf{x},t;\Theta(\mathbf{x},t))\right)\Phi_{\Theta}^{(1)}\right)\\
\hspace{1.2cm}\displaystyle-\tau div_{1}\left({\bf K}_{\Theta}^{-1}A_{\Theta,1}\nabla \mathfrak{h}_{\Theta,1}+{\bf K}_{\Theta}^{-1}A_{\Theta,2}\nabla \mathfrak{h}_{\Theta,2}\right)
\displaystyle\hspace{0cm}-\tau div_{1}\left(({\bf K}_{\Theta}^{-1}-\frac{1}{2}{\bf I}_{d})DV(\mathbf{x},t;\Theta)\right).
\end{array}
\end{equation}
By taking $\kappa_{t}=1$, $\vec{V}_{\Theta}=({\bf K}_{\Theta}^{-1}-\frac{\kappa_{x}}{2}{\bf I}_{d})^{-1}({\bf K}_{\Theta}^{-1}-\frac{1}{2}{\bf I}_{d})\frac{\p\ \mathfrak{T}_{\Theta}}{\p\ \Theta}$, 
 $A_{\Theta,k}=\frac{1}{\tau}{\bf K}_{\Theta}D_{\Theta,k}$ (for $k=1,2$), and the matrix ${\bf K}_{\Theta}$ such the diffusion tensors $D_{\Theta}= d_{\Theta}C_{s}^{2}\tau({\bf K}_{\Theta}^{-1}-\frac{1}{2}{\bf I}_{d})$, we obtain
\begin{equation}\label{PC14}
\begin{array}{lcr}
\displaystyle \hspace{-0.4cm}\frac{\p\ \Theta}{\p\ t_{2}}-div_{1}\left(D_{\Theta}\nabla_{1}\Theta\right)\!=\!\Phi_{\Theta}^{(2)}- div_{1}\left(D_{\Theta,1}\nabla \mathfrak{h}_{\Theta,1}+D_{\Theta,2}\nabla \mathfrak{h}_{\Theta,2}\right)-\tau div_{1}\left(({\bf K}_{\Theta}^{-1}-\frac{1}{2}{\bf I}_{d})DV(\mathbf{x},t;\Theta)\right).\!\!
\end{array}
\end{equation}
In equation \eqref{PC14} there exists the deviation term  $\displaystyle  div_{1}\left(({\bf K}_{\Theta}^{-1}-\frac{1}{2}{\bf I}_{d})DV(\mathbf{x},t;\Theta)\right).$ 
If this term is null (e.g., $\mathfrak{T}_{\Theta}$ is independent on $(\mathbf{x},t)$ or $\mathfrak{T}_{\Theta}$ is independent on times with  $div(\mathfrak{C}_{\Theta})(\mathbf{x},t;\Theta)=0$ and $div(\mathfrak{T}_{\Theta})(\mathbf{x},t;\Theta)=0$) or is negligible compared to the other terms, then \eqref{PC14} becomes
\begin{equation}\label{PC14-B}
\begin{array}{lcr}
\displaystyle \frac{\p\ \Theta}{\p\ t_{2}}- div_{1}\left(D_{\Theta}\nabla_{1}\Theta\right)-(\Phi_{\Theta}^{(2)}- div_{1}\left(D_{\Theta,1}\nabla \mathfrak{h}_{\Theta,1}+D_{\Theta,2}\nabla \mathfrak{h}_{\Theta,2}\right))=(\text{~or~$\approx$}) 0.
\end{array}
\end{equation}
By combining the first equation of \eqref{C146} (multiplying by $\varepsilon$) and \eqref{PC14-B} (multiplying by $\varepsilon^{2}$), we can recover equation \eqref{eq:ReacDiff}. This completes the proof of the recovery of \eqref{eq:ReacDiff} from the developed LBM method. \hfill $\Box$
%
 %
%%%%%%%%%%%%%%%%%%%%%%%%%%%%%%%%%%%%  
%
%
%
%
\subsection{Discrete adjoint-based optimization}
Since the state ${\bf u}$ depends on the microscopic distribution function ${\cal G}=(G_{i})_{i=0,q-1}$ with $G_{i}=(\varphi_{u,i}, \varphi_{v,i},\varphi_{w,i})$, the continuous lattice Boltzmann type equations (LBE) with force (for $i=0,q-1$), which is associated to the three equations of system \eqref{eq:ReacDiff}, can be reformulated as follows (since the gradient of state functions can be computed locally as function of distributions functions) 
\begin{equation}\label{LIBEG}
\begin{array}{lcr}
\displaystyle \frac{\p\ \varphi_{u,i}}{\p\ t}+{\bf e}_{i}\nabla \varphi_{u,i}=-\sum_{j}(\hat{\bf \Lambda}_{u})_{ij}\left( \varphi_{u,j}-\varphi_{u,j}^ {eq}\right)+{\mathfrak{F}}_{u,i}(.;{\cal G},{\bf f}_{1})+\mathfrak{S}_{u,i}(.;{\cal G}),\\
\displaystyle\frac{\p\ \varphi_{v,i}}{\p\ t}+{\bf e}_{i}\nabla \varphi_{v,i}=-\sum_{j}(\hat{\bf \Lambda}_{v})_{ij}\left( \varphi_{v,j}-\varphi_{v,j}^ {eq}\right) +{\mathfrak{F}}_{v,i}(.;{\cal G}, {\bf f}_{2})+\mathfrak{S}_{v,i}(.;{\cal G}),\\
\displaystyle \frac{\p\ \varphi_{w,i}}{\p\ t}+{\bf e}_{i}\nabla \varphi_{w,i}=-\sum_{j}(\hat{\bf \Lambda}_{w})_{ij}\left( \varphi_{w,j}-\varphi_{w,j}^ {eq}\right) +{\mathfrak{F}}_{w,i}(.;{\cal G},{\bf f}_{3})+\mathfrak{S}_{w,i}(.;{\cal G}),
\end{array}
\end{equation}
where ${\mathfrak{F}}_{\Theta,i}$ (resp. $\mathfrak{S}_{\Theta,i}$), for $\Theta=u,v$ or $w$ are corresponding to $\Phi_{\Theta,i}$ (resp. to ${\cal S}_{\Theta,i}$).

The corresponding discrete form of \eqref{LIBEG} is (at $(\mathbf{x},t)\in {\cal L}_{h}\times \daleth_{\tau}$)
\begin{equation}\label{eq:LBEG}
\begin{array}{lcr}
\displaystyle \varphi_{u,i}(\mathbf{x}+\mathbf{e}_i\tau,t+\tau)=\varphi_{u,i}(\mathbf{x},t)%\\[0.1cm]
\hspace{0cm}\displaystyle-\sum_{j}({\cal M}^{-1}\Pi_{u} {\cal M})_{ij}\left(\varphi_{u,j}(\mathbf{x},t)-\varphi_{u,j}^ {eq}(\mathbf{x},t)\right)+\tau {\mathfrak{S}}_{u,i}(\mathbf{x},t;{\cal G})\\
\hspace{2.5cm}\displaystyle +  \tau {\mathfrak{F}}_{u,i}(\mathbf{x},t;{\cal G},{\bf f}_{1})+ \frac{\tau^2}{2}(\frac{\partial}{\partial t}+\kappa_{x}\mathbf{e}_{i}\cdot\nabla){\mathfrak{F}}_{u,i}(\mathbf{x},t;{\cal G},{\bf f}_{1}),\\[0.1cm]
\displaystyle \varphi_{v,i}(\mathbf{x}+\mathbf{e}_i\tau,t+\tau)=\varphi_{v,i}(\mathbf{x},t)%\\[0.1cm]
\hspace{0cm}\displaystyle-\sum_{j}({\cal M}^{-1}\Pi_{v} {\cal M})_{ij}\left(\varphi_{v,j}(\mathbf{x},t)-\varphi_{v,j}^ {eq}(\mathbf{x},t)\right)+\tau {\mathfrak{S}}_{v,i}(\mathbf{x},t;{\cal G})\\
\hspace{2.5cm}\displaystyle+  \tau {\mathfrak{F}}_{v,i}(\mathbf{x},t;{\cal G},{\bf f}_{2})+ \frac{\tau^2}{2}(\frac{\partial}{\partial t}+\kappa_{x}\mathbf{e}_{i}\cdot\nabla){\mathfrak{F}}_{v,i}(\mathbf{x},t;{\cal G},{\bf f}_{2}),\\[0.1cm]
\displaystyle \varphi_{w,i}(\mathbf{x}+\mathbf{e}_i\tau,t+\tau)=\varphi_{w,i}(\mathbf{x},t)%\\[0.1cm]
\hspace{0cm}\displaystyle-\sum_{j}({\cal M}^{-1}\Pi_{w} {\cal M})_{ij}\left(\varphi_{w,j}(\mathbf{x},t)-\varphi_{w,j}^ {eq}(\mathbf{x},t)\right)+\tau {\mathfrak{S}}_{w,i}(\mathbf{x},t;{\cal G})\\
\hspace{2.5cm}\displaystyle+ \tau {\mathfrak{F}}_{w,i}(\mathbf{x},t;{\cal G},{\bf f}_{3})+ \frac{\tau^2}{2}(\frac{\partial}{\partial t}+\kappa_{x}\mathbf{e}_{i}\cdot\nabla){\mathfrak{F}}_{w,i}(\mathbf{x},t;{\cal G},{\bf f}_{3}),
\end{array}
\end{equation}
where, for $\Theta=u,v$ or $w$ and for $i=0,q-1$, the equilibrium distribution function $\varphi^{eq}_{\Theta,i}(\mathbf{x},t)$, associated with the discrete velocity $\mathbf{e}_{i}$ 
at position $\mathbf{x}$ and time $t$, is defined by \eqref{phieq} and the discrete source terms ${\mathfrak{F}}_{\Theta,i}$ and ${\mathfrak{S}}_{\Theta,i}$ are defined, respectively, 
by \eqref{PhiS} and \eqref{S}.
After each time step, we obtain new particular repartitions which lead us to recover the macroscopic state ${\bf u}$ with the following summation of microscopic particle distribution functions $G_{i}$ 
\begin{equation}\label{UGS}
{\bf u}({\bf x},t)=\sum_{i}G_{i}({\bf x},t), \text{~for each node ${\bf x}$ of lattice ${\cal L}_{h}$ and each discrete values of time $t$~~}. 
\end{equation}
\begin{req}
In order to calculate $G_{i}$ by \eqref{LIBEG}, for a given ${\bf f}$, it is necessary to approach the third terms of the right hand side of system \eqref{LIBEG} which depend on the parameter $\kappa_{x}$. For simplicity, 
in the case when  $\kappa_{x}=0$ or 1, we can use an explicit characteristic method to approximate these terms  i.e., 
$$ (\frac{\partial}{\partial t}+\kappa_{x}\mathbf{e}_{i}\cdot\nabla){\mathfrak{F}}_{\Theta,i}(\mathbf{x},t)\approx \frac{{\mathfrak{F}}_{\Theta,i}(\mathbf{x},t)-{\mathfrak{F}}_{\Theta,i}(\mathbf{x}+\kappa_{x}\mathbf{e}_i\tau,t-\tau)}{\tau} \text{~(for $\Theta=u, v$ or $w$)}.$$
We can also use an implicit characteristic method (see e.g., \cite{BIJB,ZBZ}).\hfill $\Box$
\end{req}
\begin{req}
Since the developed method is explicit time-stepping scheme, then it is necessary to impose some conditions on discretization parameters, in the spirit of the CFL condition (as e.g. 
of type $\tau\leq C_{CFL}~ h^{2}$), which guarantee the convergence of the method.\hfill $\Box$
\end{req}
According to the previous development and using similar approach  as in \cite{BIJB},  we can now formulate and develop the adjoint-based multiple-relaxation-time lattice Boltzmann. The method consists in using the adjoint of the multiple-relaxation-time LBE \eqref{LIBEG} and its time discretization form \eqref{eq:LBEG}.
\subsubsection{Adjoint multiple-relaxation-time lattice Boltzmann method and algorithm}
Since the cost functional ${\cal J}$ depends on ${\bf f}$ and ${\bf u}=(u,v,w)$, and from LBM approximation of solution $\displaystyle {\bf u}=\sum_{i=0}^{q-1}G_{i}$, with $G_{i}=(\varphi_{u,i}, \varphi_{v,i},\varphi_{w,i})$, then ${\cal J}$ can be written as the form 
\begin{equation}\label{COSTF}
{\cal J}=\int_{\Omega}{\cal I}_{0}({\cal G}(T))dx+\int_{0}^{T}\int_{\Omega} {\cal I}({\cal G}, {\bf f})dxdt,
\end{equation}
where ${\cal G}=(G_{i})_{i=0,q-1}$ and $G_{i}$ is the solution of the continuous LBE  with force \eqref{LIBEG} (for $i=0,q-1$). 
Since the distribution function ${\cal G}$ depends on control variable  ${\bf f}$ implicitly, then a variation $\delta {\bf f}$ of ${\bf f}$ causes a variation $\delta {\cal G}$ of ${\cal G}$. Consequently a variation $\delta {\cal J}$ of cost functional ${\cal J}$ can be given as
\begin{equation}\label{DeltaJ} 
\begin{array}{lcr}
\displaystyle\delta {\cal J}=\int_{\Omega}\sum_{i}\frac{\p\ {\cal I}_{0}}{\p\ G_{i}}\delta G_{i}(T)dx+\int_{0}^{T}\int_{\Omega}(\sum_{i}\frac{\p\ {\cal I}}{\p\ G_{i}} \delta G_{i}+ \frac{\p\ {\cal I}}{\p\ {\bf f}} \delta {\bf f})dxdt\\
\displaystyle\hspace{0.5cm}=\int_{\Omega}\sum_{i}(\frac{\p\ {\cal I}_{0}}{\p\ \varphi_{u,i}} \delta \varphi_{u,i}(T)+\frac{\p\ {\cal I}_{0}}{\p\ \varphi_{v,i}} \delta \varphi_{v,i}(T)
+\frac{\p\ {\cal I}_{0}}{\p\ \varphi_{w,i}} \delta \varphi_{w,i}(T))dx\\
\displaystyle \hspace{1cm}+\int_{0}^{T}\int_{\Omega}\sum_{i}(\frac{\p\ {\cal I}}{\p\ \varphi_{u,i}} \delta \varphi_{u,i}+\frac{\p\ {\cal I}}{\p\ \varphi_{v,i}} \delta \varphi_{v,i}+\frac{\p\ {\cal I}}{\p\ \varphi_{w,i}} \delta \varphi_{w,i})dxdt \\
\displaystyle \hspace{1cm}+\int_{0}^{T}\int_{\Omega}\frac{\p\ {\cal I}}{\p\ {\bf f}_{1}} \delta {\bf f}_{1} dxdt+\int_{0}^{T}\int_{\Omega}\frac{\p\ {\cal I}}{\p\ {\bf f}_{2}} \delta {\bf f}_{2} dxdt+\int_{0}^{T}\int_{\Omega}\frac{\p\ {\cal I}}{\p\ {\bf f}_{3}} \delta {\bf f}_{3} dxdt.
\end{array}
\end{equation}
We can now show  the gradient of ${\cal J}$ by introducing intermediate sensitivity and adjoint LB equations.
\begin{theo}
The  components of the gradient of ${\cal J}$ with respect to control variable ${\bf f}$ is 
\begin{equation}\label{gradientLattice1}
\begin{array}{lcr}
\displaystyle\frac{\p\ {\cal J}}{\p\ {\bf f}_{1}}=\sum_{i}\psi_{u,i}\frac{\p\ {\mathfrak{F}}_{u,i}}{\p\ {\bf f}_{1}}+\frac{\p\ {\cal I}}{\p\ {\bf f}_{1}},~~
\displaystyle\frac{\p\ {\cal J}}{\p\ {\bf f}_{2}}=\sum_{i}\psi_{v,i}\frac{\p\ {\mathfrak{F}}_{v,i}}{\p\ {\bf f}_{2}}+\frac{\p\ {\cal I}}{\p\ {\bf f}_{2}},~~
\displaystyle \frac{\p\ {\cal J}}{\p\ {\bf f}_{1}}=\sum_{i}\psi_{w,i}\frac{\p\ {\mathfrak{F}}_{w,i}}{\p\ {\bf f}_{3}}+\frac{\p\ {\cal I}}{\p\ {\bf f}_{3}}.
\end{array}
\end{equation}
with $(\psi_{u,i},\psi_{v,i},\psi_{w,i})$ the solution of the following adjoint LB equations
\begin{equation}\label{CAE1}
\begin{array}{lcr}
\displaystyle-\frac{\p\ \psi_{u,i}}{\p\ t}-{\bf e}_{i}.\nabla \psi_{u,i}
=\sum_{k}(\psi_{u,k}\frac{\p\ \mathfrak{R}_{u,k}}{\p\ \varphi_{u,i}}+\psi_{v,k}\frac{\p\ {\mathfrak{R}}_{v,k}}{\p\ \varphi_{u,i}}+\psi_{w,k}\frac{\p\ {\mathfrak{R}}_{w,k}}{\p\ \varphi_{u,i}})\\
\displaystyle\hspace{3cm}-\sum_{j}(\hat{\bf \Lambda}^{*}_{u})_{ij}\left(\psi_{u,j}-\psi^{eq}_{u,j}\right)
\displaystyle +\frac{\p\ {\cal I}}{\p\ \varphi_{u,i}},\\ 
\displaystyle-\frac{\p\ \psi_{v,i}}{\p\ t}-{\bf e}_{i}.\nabla \psi_{v,i}
=\sum_{k}(\psi_{v,k}\frac{\p\ {\mathfrak{R}}_{v,k}}{\p\ \varphi_{v,i}}+\psi_{u,k}\frac{\p\ {\mathfrak{R}}_{u,k}}{\p\ \varphi_{v,i}}+\psi_{w,k}\frac{\p\ {\mathfrak{R}}_{w,k}}{\p\ \varphi_{v,i}})\\
\displaystyle\hspace{3cm}-\sum_{j}(\hat{\bf \Lambda}^{*}_{v})_{ij}\left(\psi_{v,j}-\psi^{eq}_{v,j}\right)
\displaystyle +\frac{\p\ {\cal I}}{\p\ \varphi_{v,i}},\\ 
%+
%
\displaystyle-\frac{\p\ \psi_{w,i}}{\p\ t}-{\bf e}_{i}.\nabla \psi_{w,i}
=\sum_{k}(\psi_{w,k}\frac{\p\ {\mathfrak{R}}_{w,k}}{\p\ \varphi_{w,i}}+\psi_{u,k}\frac{\p\ {\mathfrak{R}}_{u,k}}{\p\ \varphi_{w,i}}+\psi_{v,k}\frac{\p\ {\mathfrak{R}}_{v,k}}{\p\ \varphi_{w,i}})\\
\displaystyle\hspace{3cm}-\sum_{j}(\hat{\bf \Lambda}^{*}_{w})_{ij}\left(\psi_{w,j}-\psi^{eq}_{w,j}\right)
\displaystyle +\frac{\p\ {\cal I}}{\p\ \varphi_{w,i}},\\ 
\text{with final conditions}\\
\displaystyle(\psi_{u,i}, \psi_{v,i},\psi_{w,i})(T)=(\frac{\p\ {\cal I}_{0}}{\p\ \varphi_{u,i}}, \frac{\p\ {\cal I}_{0}}{\p\ \varphi_{v,i}}, \frac{\p\ {\cal I}_{0}}{\p\ \varphi_{w,i}}),\\
\text{and boundary conditions}\\
({\bf e}_{i}.{\bf n})\psi_{u,i}=0,~({\bf e}_{i}.{\bf n})\psi_{v,i}=0,~({\bf e}_{i}.{\bf n})\psi_{w,i}=0,~~~\text{on}~\Sigma
\end{array}
\end{equation}
where adjoint equilibrium distribution functions $\psi^{eq}_{\Theta,i}$ (for $\Theta=u$, $v$ or $w$ and  $i=0, q-1$) 
\begin{equation}\label{AEDs1}
\begin{array}{lcr}
\displaystyle\psi^{eq}_{\Theta,i}=\sum_{k}((\hat{\bf \Lambda}^{*}_{\Theta})^{-1}{\cal O}_{\Theta}\hat{\bf \Lambda}^{*}_{\Theta})_{ik}\psi_{\Theta,k}
\end{array}
\end{equation} 
and 
\begin{equation}\label{RFS1}
\begin{array}{lcr}
{\mathfrak{R}}_{\Theta,k}={\mathfrak{F}}_{\Theta,k}+{\mathfrak{S}}_{\Theta,k},~~\text{for $k=0, q-1$},\\
\hat{\bf \Lambda}^{*}_{\Theta}=(\hat{\bf \Lambda}_{\Theta})^{t},~~
{\cal O}_{\Theta}=(\frac{\p\ \varphi_{\Theta,j}^{eq}}{\p\ \varphi_{\Theta,i}})_{{ij}}.
\end{array}
\end{equation}
\end{theo}
{\bf Proof:}
Since $G_{i}=(\varphi_{u,i},\varphi_{v,i},\varphi_{w,i})$ is a solution of \eqref{LIBEG}, then the variation $\delta G_{i}=(\delta\varphi_{u,i}, \delta\varphi_{v,i},\delta\varphi_{w,i})$ is a solution of the following {\it sensitivity lattice Boltzmann equation}
\begin{equation}\label{LIBEG-V}
\begin{array}{lcr}
\displaystyle \frac{\p\ \delta\varphi_{u,i}}{\p\ t}+{\bf e}_{i}\nabla \delta\varphi_{u,i}=-\sum_{j}(\hat{\bf \Lambda}_{u})_{ij}\left(\delta \varphi_{u,j}-\sum_{k}\frac{\p\ \varphi_{u,j}^ {eq}}{\p\ \varphi_{u,k}} \delta\varphi_{u,k}\right)+\frac{\p\ {\mathfrak{F}}_{u,i}}{\p\ {\bf f}_{1}}\delta {\bf f}_{1}\\ 
\displaystyle \hspace{2cm}+\sum_{k}\frac{\p\ {\mathfrak{F}}_{u,i}}{\p\ G_{k}}\delta G_{k}+\sum_{k}\frac{\p\ {\mathfrak{S}}_{u,i}}{\p\ G_{k}}\delta G_{k},\\
\displaystyle\frac{\p\ \delta\varphi_{v,i}}{\p\ t}+{\bf e}_{i}\nabla \delta\varphi_{v,i}=-\sum_{j}(\hat{\bf \Lambda}_{v})_{ij}\left( \delta \varphi_{v,j}-\sum_{k}\frac{\p\ \varphi_{v,j}^ {eq}}{\p\ \varphi_{v,k}}\delta\varphi_{v,k}\right)+\frac{\p\ {\mathfrak{F}}_{v,i}}{\p\ {\bf f}_{2}}\delta {\bf f}_{2}\\
\displaystyle \hspace{2cm}+\sum_{k}\frac{\p\ {\mathfrak{F}}_{v,i}}{\p\ G_{k}}\delta G_{k}+\sum_{k}\frac{\p\ {\mathfrak{S}}_{v,i}}{\p\ G_{k}}\delta G_{k}
,\\
\displaystyle \frac{\p\ \delta\varphi_{w,i}}{\p\ t}+{\bf e}_{i}\nabla \delta\varphi_{w,i}=-\sum_{j}(\hat{\bf \Lambda}_{w})_{ij}\left( \delta \varphi_{w,j}-\sum_{k}\frac{\p\ \varphi_{w,j}^ {eq}}{\p\ \varphi_{w,k}}\delta\varphi_{w,k}\right)+\frac{\p\ {\mathfrak{F}}_{w,i}}{\p\ {\bf f}_{3}}\delta {\bf f}_{3}\\
\\ \displaystyle \hspace{2cm} +\sum_{k}\frac{\p\ {\mathfrak{F}}_{w,i}}{\p\ G_{k}}\delta G_{k}+\sum_{k}\frac{\p\ {\mathfrak{S}}_{w,i}}{\p\ G_{k}}\delta G_{k},\\
\text{with the initial conditions}\\
\delta G_{i}(t=0)=(\delta\varphi_{u,i}, \delta\varphi_{v,i},\delta\varphi_{w,i})(t=0)=(0,0,0).
\end{array}
\end{equation}
Let $H_{i}=(\psi_{u,i}, \psi_{v,i},\psi_{w,i})$ be sufficiently regular such that $\displaystyle H_{i}(T)=(\frac{\p\ {\cal I}_{0}}{\p\ \varphi_{u,i}}, \frac{\p\ {\cal I}_{0}}{\p\ \varphi_{v,i}}, \frac{\p\ {\cal I}_{0}}{\p\ \varphi_{w,i}})$. Multiplying \eqref{LIBEG-V} by $H_{i}$, summing on $i$ and integrating over $\Omega\times (0,T)$ this gives (since $\delta G_{i}(t=0)=(0,0,0)$)
\begin{equation*}%
\begin{array}{lcr}
\displaystyle -\int_{0}^{T}\!\!\!\int_{\Omega}\sum_{i}\frac{\p\ \psi_{u,i}}{\p\ t}\delta\varphi_{u,i}dxdt-\int_{0}^{T}\!\!\!\int_{\Omega}\sum_{i}({\bf e}_{i}\nabla \psi_{u,i})\delta\varphi_{u,i}dxdt\\[0.4cm]
\displaystyle \hspace{1cm}+\int_{0}^{T}\!\!\!\int_{\p\ \Omega}\sum_{i}({\bf e}_{i}.{\bf n})\psi_{u,i} \delta\varphi_{u,i}dxdt+\int_{\Omega}\sum_{i}\frac{\p\ {\cal I}_{0}}{\p\ \varphi_{u,i}}\delta\varphi_{u,i}(T)dx\\[0.3cm]
\displaystyle\hspace{1cm}+\int_{0}^{T}\!\!\!\int_{\Omega}\sum_{i}\left(\sum_{j}(\hat{\bf \Lambda}_{u})_{ji}\psi_{u,j}-\sum_{k,j}(\hat{\bf \Lambda}_{u})_{kj}\frac{\p\ \varphi_{u,j}^ {eq}}{\p\ \varphi_{u,i}}\psi_{u,k}\right) \delta\varphi_{u,i} dxdt\\[0.3cm] 
\displaystyle \hspace{1cm}-\int_{0}^{T}\!\!\!\int_{\Omega}\sum_{i}\left(\sum_{k}\psi_{u,k}(\frac{\p\ {\mathfrak{F}}_{u,k}}{\p\ G_{i}}+\frac{\p\ {\mathfrak{S}}_{u,k}}{\p\ G_{i}})\right).\delta G_{i} dxdt
\displaystyle \hspace{0cm}-\int_{0}^{T}\!\!\!\int_{\Omega}\sum_{i}\psi_{u,i}\frac{\p\ {\mathfrak{F}}_{u,i}}{\p\ {\bf f}_{1}}.\delta {\bf f}_{1} dxdt
=0,
\end{array}
\end{equation*}
\begin{equation*}
\begin{array}{lcr}
\displaystyle-\int_{0}^{T}\!\!\!\int_{\Omega}\sum_{i}\frac{\p\ \psi_{v,i}}{\p\ t}\delta\varphi_{v,i}dxdt-\int_{0}^{T}\!\!\!\int_{\Omega}\sum_{i}({\bf e}_{i}\nabla \psi_{v,i})\delta\varphi_{v,i}dxdt\\[0.4cm]
\displaystyle \hspace{1cm}+\int_{0}^{T}\!\!\!\int_{\p\ \Omega}\sum_{i}({\bf e}_{i}.{\bf n})\psi_{v,i} \delta\varphi_{v,i}dxdt+\int_{\Omega}\sum_{i}\frac{\p\ {\cal I}_{0}}{\p\ \varphi_{v,i}}\delta\varphi_{v,i}(T)dx\\[0.3cm]
\displaystyle\hspace{1cm}+\int_{0}^{T}\!\!\!\int_{\Omega}\sum_{i}\left(\sum_{j}(\hat{\bf \Lambda}_{v})_{ji}\psi_{v,j}-\sum_{k,j}(\hat{\bf \Lambda}_{v})_{kj}\frac{\p\ \varphi_{v,j}^ {eq}}{\p\ \varphi_{v,i}}\psi_{v,k}\right) \delta\varphi_{v,i} dxdt\\[0.3cm] 
\displaystyle \hspace{1cm}-\int_{0}^{T}\!\!\!\int_{\Omega}\sum_{i}\left(\sum_{k}\psi_{v,k}(\frac{\p\ {\mathfrak{F}}_{v,k}}{\p\ G_{i}}+\frac{\p\ {\mathfrak{S}}_{v,k}}{\p\ G_{i}})\right).\delta G_{i} dxdt
\displaystyle \hspace{0cm}-\int_{0}^{T}\!\!\!\int_{\Omega}\sum_{i}\psi_{v,i}\frac{\p\ {\mathfrak{F}}_{v,i}}{\p\ {\bf f}_{2}}.\delta {\bf f}_{2} dxdt=0,
\end{array}
\end{equation*}
and
\begin{equation*}
\begin{array}{lcr}
\displaystyle -\int_{0}^{T}\!\!\!\int_{\Omega}\sum_{i}\frac{\p\ \psi_{w,i}}{\p\ t}\delta\varphi_{w,i}dxdt-\int_{0}^{T}\!\!\!\int_{\Omega}\sum_{i}({\bf e}_{i}\nabla \psi_{w,i})\delta\varphi_{w,i}dxdt\\[0.4cm]
\displaystyle \hspace{1cm}+\int_{0}^{T}\!\!\!\int_{\p\ \Omega}\sum_{i}({\bf e}_{i}.{\bf n})\psi_{w,i} \delta\varphi_{w,i}dxdt+\int_{\Omega}\sum_{i}\frac{\p\ {\cal I}_{0}}{\p\ \varphi_{w,i}}\delta\varphi_{w,i}(T)dx\\[0.3cm]
\displaystyle\hspace{1cm}+\int_{0}^{T}\!\!\!\int_{\Omega}\sum_{i}\left(\sum_{j}(\hat{\bf \Lambda}_{w})_{ji}\psi_{w,j}-\sum_{k,j}(\hat{\bf \Lambda}_{w})_{kj}\frac{\p\ \varphi_{w,j}^ {eq}}{\p\ \varphi_{w,i}}\psi_{w,k}\right) \delta\varphi_{w,i} dxdt\\[0.3cm] 
\displaystyle \hspace{1cm}-\int_{0}^{T}\!\!\!\int_{\Omega}\sum_{i}\left(\sum_{k}\psi_{w,k}(\frac{\p\ {\mathfrak{F}}_{w,k}}{\p\ G_{i}}+\frac{\p\ {\mathfrak{S}}_{w,k}}{\p\ G_{i}})\right).\delta G_{i} dxdt
\displaystyle \hspace{0cm}-\int_{0}^{T}\!\!\!\int_{\Omega}\sum_{i}\psi_{w,i}\frac{\p\ {\mathfrak{F}}_{w,i}}{\p\ {\bf f}_{3}}.\delta {\bf f}_{3} dxdt=0.
\end{array}
\end{equation*}
By setting for $\Theta=u,v$ or $w$, $({\bf e}_{i}.{\bf n})\psi_{\Theta,i}=0$ on $\p\ \Omega$, and noting by %$F^{(\rho)}=(\varphi_{\Theta,i})_{i}$, 
\begin{equation}\label{EQAdj}
\text{$\Psi^{(\Theta)}=(\psi_{\Theta,i})_{i}$, $\hat{\bf \Lambda}^{*}_{\Theta}=(\hat{\bf \Lambda}_{\Theta})^{t}$ and ${\cal O}_{\Theta}=(\frac{\p\ \varphi_{\Theta,j}^{eq}}{\p\ \varphi_{\Theta,i}})_{{ij}}$}, 
\end{equation}
we can deduce  
\begin{equation}\label{LIBEG-VII}
\hspace{-0.2cm}
\begin{array}{lcr}
\displaystyle -\int_{0}^{T}\!\!\!\int_{\Omega}\sum_{i}\frac{\p\ \psi_{u,i}}{\p\ t}\delta\varphi_{u,i}dxdt-\int_{0}^{T}\!\!\!\int_{\Omega}\sum_{i}({\bf e}_{i}\nabla \psi_{u,i})\delta\varphi_{u,i}dxdt\\
\displaystyle+\int_{0}^{T}\!\!\!\int_{\Omega}\sum_{i}\left(\sum_{j}(\hat{\bf \Lambda}^{*}_{u})_{ij}\left(\psi_{u,j}-\sum_{k}((\hat{\bf \Lambda}^{*}_{u})^{-1}{\cal O}_{u}\hat{\bf \Lambda}^{*}_{u})_{jk}\psi_{u,k})\right)\right) \delta\varphi_{u,i} dxdt\\ 
\displaystyle \hspace{1cm}+\int_{\Omega}\sum_{i}\frac{\p\ {\cal I}_{0}}{\p\ \varphi_{u,i}}\delta\varphi_{u,i}(T)dx-\int_{0}^{T}\!\!\!\int_{\Omega}\sum_{i}\psi_{u,i}\frac{\p\ {\mathfrak{F}}_{u,i}}{\p\ {\bf f}_{1}}.\delta {\bf f}_{1} dxdt\\
\displaystyle \hspace{1cm}-\int_{0}^{T}\!\!\!\int_{\Omega}\sum_{i}\left(\sum_{k}\psi_{u,k}(\frac{\p\ {\mathfrak{F}}_{u,k}}{\p\ G_{i}}+\frac{\p\ {\mathfrak{S}}_{u,k}}{\p\ G_{i}})\right).\delta G_{i} dxdt=0,\\ 
\displaystyle-\int_{0}^{T}\!\!\!\int_{\Omega}\sum_{i}\frac{\p\ \psi_{v,i}}{\p\ t}\delta\varphi_{v,i}dxdt-\int_{0}^{T}\!\!\!\int_{\Omega}\sum_{i}({\bf e}_{i}\nabla \psi_{v,i})\delta\varphi_{v,i}dxdt\\
\displaystyle+\int_{0}^{T}\!\!\!\int_{\Omega}\sum_{i}\left(\sum_{j}(\hat{\bf \Lambda}^{*}_{v})_{ij}\left(\psi_{v,j}-\sum_{k}((\hat{\bf \Lambda}^{*}_{v})^{-1}{\cal O}_{v}\hat{\bf \Lambda}^{*}_{v})_{jk}\psi_{v,k})\right)\right) \delta\varphi_{v,i} dxdt\\ 
\displaystyle \hspace{1cm}+\int_{\Omega}\sum_{i}\frac{\p\ {\cal I}_{0}}{\p\ \varphi_{v,i}}\delta\varphi_{v,i}(T)dx-\int_{0}^{T}\!\!\!\int_{\Omega}\sum_{i}\psi_{v,i}\frac{\p\ {\mathfrak{F}}_{v,i}}{\p\ {\bf f}_{2}}.\delta {\bf f}_{2} dxdt\\ 
\displaystyle \hspace{1cm}-\int_{0}^{T}\!\!\!\int_{\Omega}\sum_{i}\left(\sum_{k}\psi_{v,k}(\frac{\p\ {\mathfrak{F}}_{v,k}}{\p\ G_{i}}+\frac{\p\ {\mathfrak{S}}_{v,k}}{\p\ G_{i}})\right).\delta G_{i} dxdt=0,\\ 
\displaystyle -\int_{0}^{T}\!\!\!\int_{\Omega}\sum_{i}\frac{\p\ \psi_{w,i}}{\p\ t}\delta\varphi_{w,i}dxdt-\int_{0}^{T}\!\!\!\int_{\Omega}\sum_{i}({\bf e}_{i}\nabla \psi_{w,i})\delta\varphi_{w,i}dxdt\\
\displaystyle+\int_{0}^{T}\!\!\!\int_{\Omega}\sum_{i}\left(\sum_{j}(\hat{\bf \Lambda}^{*}_{w})_{ij}\left(\psi_{w,j}-\sum_{k}((\hat{\bf \Lambda}^{*}_{w})^{-1}{\cal O}_{w}\hat{\bf \Lambda}^{*}_{w})_{jk}\psi_{w,k})\right)\right) \delta\varphi_{w,i} dxdt\\ 
\displaystyle \hspace{1cm}+\int_{\Omega}\sum_{i}\frac{\p\ {\cal I}_{0}}{\p\ \varphi_{w,i}}\delta\varphi_{w,i}(T)dx-\int_{0}^{T}\!\!\!\int_{\Omega}\sum_{i}\psi_{w,i}\frac{\p\ {\mathfrak{F}}_{w,i}}{\p\ {\bf f}_{3}}.\delta {\bf f}_{3} dxdt\\ 
\displaystyle \hspace{1cm}-\int_{0}^{T}\!\!\!\int_{\Omega}\sum_{i}\left(\sum_{k}\psi_{w,k}(\frac{\p\ {\mathfrak{F}}_{w,k}}{\p\ G_{i}}+\frac{\p\ {\mathfrak{S}}_{w,k}}{\p\ G_{i}})\right).\delta G_{i} dxdt
=0,
\end{array}
\end{equation}
Introduce the following adjoint equilibrium distribution functions ($\text{for~} \Theta=u,~v~\text{or~} w$ and $i=0, q-1$)
\begin{equation}\label{AEDs}
\begin{array}{lcr}
\displaystyle\psi^{eq}_{\Theta,i}=\sum_{k}((\hat{\bf \Lambda}^{*}_{\Theta})^{-1}{\cal O}_{\Theta}\hat{\bf \Lambda}^{*}_{\Theta})_{ik}\psi_{\Theta,k}.
\end{array}
\end{equation}
Otherwise
\begin{equation}\label{AEDsM}
\begin{array}{lcr}
\displaystyle \tilde{m}^{eq}_{\Theta,i}=\sum_{j}({\cal O}_{\Theta})_{ij}\tilde{m}_{\Theta,j} ~(i.e.~ \tilde{\bf m}^{eq}_{\Theta}={\cal O}_{\Theta}\tilde{\bf m}_{\Theta}),
\end{array}
\end{equation}
where $\Psi_{\Theta}=(\psi_{\Theta,i})_{i}$, $\Psi^{eq}_{\Theta}=(\psi^{eq}_{\Theta,i})_{i}$, $\tilde{\bf m}_{\Theta}=\hat{\bf \Lambda}^{*}_{\Theta}\Psi_{\Theta}$, 
$\tilde{\bf m}^{eq}_{\Theta}=\hat{\bf \Lambda}^{*}_{\Theta}\Psi^{eq}_{\Theta}$. Then \eqref{LIBEG-VII} becomes
\begin{equation}\label{LIBEG-VIII}
\hspace{-0.1cm}
\begin{array}{lcr}
\displaystyle \!-\!\!\int_{0}^{T}\!\!\!\int_{\Omega}\sum_{i}\left[\frac{\p\ \psi_{u,i}}{\p\ t}-{\bf e}_{i}.\nabla \psi_{u,i}
\displaystyle+\sum_{j}(\hat{\bf \Lambda}^{*}_{u})_{ij}\left(\psi_{u,j}-\psi^{eq}_{u,j}\right) \right.
\displaystyle \hspace{0cm}\left.\!-\!\sum_{k}\psi_{u,k}(\frac{\p\ {\mathfrak{F}}_{u,k}}{\p\ \varphi_{u,i}}+\frac{\p\ {\mathfrak{S}}_{u,k}}{\p\ \varphi_{u,i}})\right] \delta\varphi_{u,i} dxdt\!\\
\displaystyle \hspace{0.5cm}-\int_{0}^{T}\!\!\!\int_{\Omega}\sum_{i}\left[\sum_{k}\psi_{u,k}\left(\frac{\p\ {\mathfrak{F}}_{u,k}}{\p\ \varphi_{v,i}}+\frac{\p\ {\mathfrak{S}}_{u,k}}{\p\ \varphi_{v,i}}\right)\right]\delta\varphi_{v,i} dxdt\\
\displaystyle \hspace{0.5cm}-\int_{0}^{T}\!\!\!\int_{\Omega}\sum_{i}\left[\sum_{k}\psi_{u,k}\left(\frac{\p\ {\mathfrak{F}}_{u,k}}{\p\ \varphi_{w,i}}+\frac{\p\ {\mathfrak{S}}_{u,k}}{\p\ \varphi_{w,i}}\right)\right]\delta\varphi_{w,i} dxdt\\
\displaystyle \hspace{0.5cm}+\int_{\Omega}\sum_{i}\frac{\p\ {\cal I}_{0}}{\p\ \varphi_{u,i}}\delta\varphi_{u,i}(T)dx-\int_{0}^{T}\int_{\Omega}\sum_{i}\left(\psi_{u,i}\frac{\p\ {\mathfrak{F}}_{u,i}}{\p\ {\bf f}_{1}}.\delta {\bf f}_{1}\right) dxdt=0,\\ 
\displaystyle \!-\!\!\int_{0}^{T}\!\!\!\int_{\Omega}\sum_{i}\left[\frac{\p\ \psi_{v,i}}{\p\ t}-{\bf e}_{i}.\nabla \psi_{v,i}
\displaystyle+\sum_{j}(\hat{\bf \Lambda}^{*}_{v})_{ij}\left(\psi_{v,j}-\psi^{eq}_{v,j}\right)\right.
\displaystyle \hspace{0cm}\left.\!-\!\sum_{k}\psi_{v,k}(\frac{\p\ {\mathfrak{F}}_{v,k}}{\p\ \varphi_{v,i}}+\frac{\p\ {\mathfrak{S}}_{v,k}}{\p\ \varphi_{v,i}})\right] \delta\varphi_{v,i} dxdt\!\\ 
\displaystyle \hspace{0.5cm}-\int_{0}^{T}\!\!\!\int_{\Omega}\sum_{i}\left[\sum_{k}\psi_{v,k}\left(\frac{\p\ {\mathfrak{F}}_{v,k}}{\p\ \varphi_{u,i}}+\frac{\p\ {\mathfrak{S}}_{v,k}}{\p\ \varphi_{u,i}}\right)\right]\delta\varphi_{u,i} dxdt\\
\displaystyle \hspace{0.5cm}-\int_{0}^{T}\!\!\!\int_{\Omega}\sum_{i}\left[\sum_{k}\psi_{v,k}\left(\frac{\p\ {\mathfrak{F}}_{v,k}}{\p\ \varphi_{w,i}}+\frac{\p\ {\mathfrak{S}}_{v,k}}{\p\ \varphi_{w,i}}\right)\right]\delta\varphi_{w,i} dxdt\\
\displaystyle \hspace{0.5cm}+\int_{\Omega}\sum_{i}\frac{\p\ {\cal I}_{0}}{\p\ \varphi_{v,i}}\delta\varphi_{v,i}(T)dx-\int_{0}^{T}\int_{\Omega}\sum_{i}\left(\psi_{v,i}\frac{\p\ {\mathfrak{F}}_{v,i}}{\p\ {\bf f}_{2}}.\delta {\bf f}_{2}\right) dxdt=0,\\ 
\displaystyle \!-\!\!\int_{0}^{T}\!\!\!\int_{\Omega}\sum_{i}\left[\frac{\p\ \psi_{w,i}}{\p\ t}-{\bf e}_{i}.\nabla \psi_{w,i}
\displaystyle+\sum_{j}(\hat{\bf \Lambda}^{*}_{w})_{ij}\left(\psi_{w,j}-\psi^{eq}_{w,j}\right)\right.
\displaystyle \hspace{0cm}\left.\!-\!\sum_{k}\psi_{w,k}(\frac{\p\ {\mathfrak{F}}_{w,k}}{\p\ \varphi_{w,i}}+\frac{\p\ {\mathfrak{S}}_{w,k}}{\p\ \varphi_{w,i}})\right] \delta\varphi_{w,i} dxdt\!\\ 
\displaystyle \hspace{0.5cm}-\int_{0}^{T}\!\!\!\int_{\Omega}\sum_{i}\left[\sum_{k}\psi_{w,k}\left(\frac{\p\ {\mathfrak{F}}_{w,k}}{\p\ \varphi_{u,i}}+\frac{\p\ {\mathfrak{S}}_{w,k}}{\p\ \varphi_{u,i}}\right)\right]\delta\varphi_{u,i} dxdt\\
\displaystyle \hspace{0.5cm}-\int_{0}^{T}\!\!\!\int_{\Omega}\sum_{i}\left[\sum_{k}\psi_{w,k}\left(\frac{\p\ {\mathfrak{F}}_{w,k}}{\p\ \varphi_{v,i}}+\frac{\p\ {\mathfrak{S}}_{w,k}}{\p\ \varphi_{v,i}}\right)\right]\delta\varphi_{v,i} dxdt\\
\displaystyle \hspace{0.5cm}+\int_{\Omega}\sum_{i}\frac{\p\ {\cal I}_{0}}{\p\ \varphi_{w,i}}\delta\varphi_{w,i}(T)dx-\int_{0}^{T}\int_{\Omega}\sum_{i}\left(\psi_{w,i}\frac{\p\ {\mathfrak{F}}_{w,i}}{\p\ {\bf f}_{3}}.\delta {\bf f}_{3}\right) dxdt=0.
\end{array}
\end{equation}
Assuming that $(\psi_{u,i},\psi_{v,i},\psi_{w,i})$ is a solution of the following {\it adjoint lattice Boltzmann equation}
\begin{equation}\label{CAE}
\begin{array}{lcr}
\displaystyle-\frac{\p\ \psi_{u,i}}{\p\ t}-{\bf e}_{i}.\nabla \psi_{u,i}
=\sum_{k}(\psi_{u,k}\frac{\p\ \mathfrak{R}_{u,k}}{\p\ \varphi_{u,i}}+\psi_{v,k}\frac{\p\ {\mathfrak{R}}_{v,k}}{\p\ \varphi_{u,i}}+\psi_{w,k}\frac{\p\ {\mathfrak{R}}_{w,k}}{\p\ \varphi_{u,i}})\\
\displaystyle\hspace{3cm}-\sum_{j}(\hat{\bf \Lambda}^{*}_{u})_{ij}\left(\psi_{u,j}-\psi^{eq}_{u,j}\right)
\displaystyle +\frac{\p\ {\cal I}}{\p\ \varphi_{u,i}},\\ 
\displaystyle-\frac{\p\ \psi_{v,i}}{\p\ t}-{\bf e}_{i}.\nabla \psi_{v,i}
=\sum_{k}(\psi_{v,k}\frac{\p\ {\mathfrak{R}}_{v,k}}{\p\ \varphi_{v,i}}+\psi_{u,k}\frac{\p\ {\mathfrak{R}}_{u,k}}{\p\ \varphi_{v,i}}+\psi_{w,k}\frac{\p\ {\mathfrak{R}}_{w,k}}{\p\ \varphi_{v,i}})\\
\displaystyle\hspace{3cm}-\sum_{j}(\hat{\bf \Lambda}^{*}_{v})_{ij}\left(\psi_{v,j}-\psi^{eq}_{v,j}\right)
\displaystyle +\frac{\p\ {\cal I}}{\p\ \varphi_{v,i}},\\ 
\displaystyle-\frac{\p\ \psi_{w,i}}{\p\ t}-{\bf e}_{i}.\nabla \psi_{w,i}
=\sum_{k}(\psi_{w,k}\frac{\p\ {\mathfrak{R}}_{w,k}}{\p\ \varphi_{w,i}}+\psi_{u,k}\frac{\p\ {\mathfrak{R}}_{u,k}}{\p\ \varphi_{w,i}}+\psi_{v,k}\frac{\p\ {\mathfrak{R}}_{v,k}}{\p\ \varphi_{w,i}})\\
\displaystyle\hspace{3cm}-\sum_{j}(\hat{\bf \Lambda}^{*}_{w})_{ij}\left(\psi_{w,j}-\psi^{eq}_{w,j}\right)
\displaystyle +\frac{\p\ {\cal I}}{\p\ \varphi_{w,i}},\\ 
\text{with the final conditions}: 
\displaystyle(\psi_{u,i}, \psi_{v,i},\psi_{w,i})(T)=(\frac{\p\ {\cal I}_{0}}{\p\ \varphi_{u,i}}, \frac{\p\ {\cal I}_{0}}{\p\ \varphi_{v,i}}, \frac{\p\ {\cal I}_{0}}{\p\ \varphi_{w,i}}),
\end{array}
\end{equation}
where (for $\Theta=u$, $v$ or $w$) : 
${\mathfrak{R}}_{\Theta,k}={\mathfrak{F}}_{\Theta,k}+{\mathfrak{S}}_{\Theta,k}.$\\
By adding the first part, second part and third part of \eqref{LIBEG-VIII} we have (according to \eqref{CAE})
\begin{equation}\label{LIBEG-VIV}
\begin{array}{lcr}
\displaystyle \int_{0}^{T}\int_{\Omega}\sum_{i}\left[\left(\psi_{u,i}\frac{\p\ {\mathfrak{F}}_{u,i}}{\p\ {\bf f}_{1}}\right).\delta {\bf f}_{1} +\left(\psi_{v,i}\frac{\p\ {\mathfrak{F}}_{v,i}}{\p\ {\bf f}_{2}}\right).\delta {\bf f}_{2}+\left(\psi_{w,i}\frac{\p\ {\mathfrak{F}}_{w,i}}{\p\ {\bf f}_{3}}\right).\delta {\bf f}_{3}\right]dxdt \\
\displaystyle\hspace{1cm}=\int_{0}^{T}\int_{\Omega}\sum_{i}\left(\frac{\p\ {\cal I}}{\p\ \varphi_{u,i}} \delta \varphi_{u,i}+\frac{\p\ {\cal I}}{\p\ \varphi_{v,i}} \delta \varphi_{v,i}+\frac{\p\ {\cal I}}{\p\ \varphi_{w,i}} \delta \varphi_{w,i}\right)dxdt\\
\displaystyle\hspace{1.5cm}+\int_{\Omega}\sum_{i}\left(\frac{\p\ {\cal I}_{0}}{\p\ \varphi_{u,i}}\delta\varphi_{u,i}(T)+\frac{\p\ {\cal I}_{0}}{\p\ \varphi_{v,i}}\delta\varphi_{v,i}(T)+\frac{\p\ {\cal I}_{0}}{\p\ \varphi_{w,i}}\delta\varphi_{w,i}(T)\right)dx.
\end{array}
\end{equation}
Then \eqref{DeltaJ} becomes
\begin{equation}\label{DeltaJS} 
\begin{array}{lcr}
\displaystyle\delta {\cal J}&=\displaystyle \int_{0}^{T}\!\!\!\int_{\Omega}\left(\sum_{i}\psi_{u,i}\frac{\p\ {\mathfrak{F}}_{u,i}}{\p\ {\bf f}_{1}}+\frac{\p\ {\cal I}}{\p\ {\bf f}_{1}}\right) \delta {\bf f}_{1}dxdt\\
  &\displaystyle+\int_{0}^{T}\!\!\!\int_{\Omega}\left(\sum_{i}\psi_{v,i}\frac{\p\ {\mathfrak{F}}_{v,i}}{\p\ {\bf f}_{2}}+\frac{\p\ {\cal I}}{\p\ {\bf f}_{2}}\right) \delta {\bf f}_{2}dxdt \\
\displaystyle &\displaystyle +\int_{0}^{T}\!\!\!\int_{\Omega}\left(\sum_{i}\psi_{w,i}\frac{\p\ {\mathfrak{F}}_{w,i}}{\p\ {\bf f}_{3}}+\frac{\p\ {\cal I}}{\p\ {\bf f}_{3}}\right) \delta {\bf f}_{3}dxdt.
\end{array}
\end{equation}
Consequently the components of the gradient of ${\cal J}$ with respect to variable ${\bf f}$ is then 
\begin{equation}\label{gradientLattice}
\begin{array}{lcr}
\displaystyle\frac{\p\ {\cal J}}{\p\ {\bf f}_{1}}=\sum_{i}\psi_{u,i}\frac{\p\ {\mathfrak{F}}_{u,i}}{\p\ {\bf f}_{1}}+\frac{\p\ {\cal I}}{\p\ {\bf f}_{1}},~~
\displaystyle\frac{\p\ {\cal J}}{\p\ {\bf f}_{2}}=\sum_{i}\psi_{v,i}\frac{\p\ {\mathfrak{F}}_{v,i}}{\p\ {\bf f}_{2}}+\frac{\p\ {\cal I}}{\p\ {\bf f}_{2}},~~
\displaystyle \frac{\p\ {\cal J}}{\p\ {\bf f}_{1}}=\sum_{i}\psi_{w,i}\frac{\p\ {\mathfrak{F}}_{w,i}}{\p\ {\bf f}_{3}}+\frac{\p\ {\cal I}}{\p\ {\bf f}_{3}}.
\end{array}
\end{equation}
This complete the proof. \hfill $\Box$

To close the adjoint lattice Boltzmann equation \eqref{CAE} and sensitivity lattice Boltzmann equation \eqref{LIBEG-V}, it is necessary to evaluate the adjoint distribution function $\psi^{eq}_{\Theta,i}$ (fo $\Theta=u,v~\text{or~} w$). For this, it is sufficient (according to \eqref{AEDsM} and \eqref{EQAdj}) to derive $\displaystyle\frac{\p\ \varphi^{eq}_{\Theta,j}}{\p\ \varphi_{\Theta,i}}$.  
According to \eqref{TF0} and \eqref{phieq}, we can compute each $\displaystyle\frac{\p\ \varphi^{eq}_{\Theta,j}}{\p\ \varphi_{\Theta,i}}$ as follows:
\begin{equation}
\begin{array}{lcr}
\displaystyle\frac{\p\ \Theta}{\p\ \varphi_{\Theta,i}}=1,\\
\displaystyle\frac{\p\ \varphi^{eq}_{\Theta,j}}{\p\ \varphi_{\Theta,i}}=\frac{\p\ \varphi^{eq}_{\Theta,j}}{\p\ \Theta}\frac{\p\ \Theta}{\p\ \varphi_{\Theta,i}}
\\
\displaystyle \hspace{1cm}=\omega_{j}\left[1+\frac{1}{C_{s}^{2}}\frac{\p\ \mathfrak{T}_{\Theta}}{\p\ \Theta}.{\bf e}_{j}+\frac{1}{2C_{s}^{4}}(\frac{\p\ \mathfrak{C}_{\Theta}}{\p\ \Theta}+(d_{\Theta}-1)C_{s}^{2}{\bf I}_{d}):\mho_{j}\right]=\omega_{j}\aleph_{j,\Theta}.
\end{array}
\end{equation}
The expression of $\aleph_{j,\Theta}$ can be derived easily from the expressions of $\frac{\p\ \mathfrak{C}_{\Theta}}{\p\ \Theta}$ and $\mho_{j}$. 
Consequently, from \eqref{AEDsM}, the adjoint equilibrium function $\tilde{\bf m}^{eq}_{\Theta}$  satisfies
\begin{equation}%\label{AdEF}
\displaystyle\tilde{m}^{eq}_{\Theta,i}={\cal T}_{\Theta},%+{\cal Y}_{\rho}.\frac{({\bf e}_{i}-\vec{\omega}^{(\rho)})}{C_{s}^{2}},
\end{equation}
where 
$\displaystyle{\cal T}_{\Theta}=\sum_{j}\omega_{j}\tilde{m}_{j}\aleph_{j,\Theta}$.

We can now present numerical algorithm to solve the  optimal control problem. 
The continuous adjoint LBE equation \eqref{CAE1} in time and space with force can be written as follows
\begin{equation}\label{CAEL}
\begin{array}{lcr}
\displaystyle-\frac{\p\ \psi_{u,i}}{\p\ t}-{\bf e}_{i}.\nabla \psi_{u,i}
=-\sum_{j}(\hat{\bf \Lambda}^{*}_{u})_{ij}\left(\psi_{u,j}-\psi^{eq}_{u,j}\right)+\tilde{\mathfrak{R}}_{u,i},\\
\displaystyle-\frac{\p\ \psi_{v,i}}{\p\ t}-{\bf e}_{i}.\nabla \psi_{v,i}
=-\sum_{j}(\hat{\bf \Lambda}_{v}^{*})_{ij}\left(\psi_{v,j}-\psi^{eq}_{v,j}\right)+ \tilde{\mathfrak{R}}_{v,i},\\
%
%(\hat{\bf \Lambda}_{u})
\displaystyle-\frac{\p\ \psi_{w,i}}{\p\ t}-{\bf e}_{i}.\nabla \psi_{w,i}=
-\sum_{j}(\hat{\bf \Lambda}_{w}^{*})_{ij}\left(\psi_{w,j}-\psi^{eq}_{w,j}\right)+\tilde{\mathfrak{R}}_{w,i},\\
\text{with the final conditions}\\
\displaystyle(\psi_{u,i}, \psi_{v,i},\psi_{w,i})(T)=(\frac{\p\ {\cal I}_{0}}{\p\ \varphi_{u,i}}, \frac{\p\ {\cal I}_{0}}{\p\ \varphi_{v,i}}, \frac{\p\ {\cal I}_{0}}{\p\ \varphi_{w,i}}),
\end{array}
\end{equation}
where 
\begin{equation}\label{CAE-L}
\begin{array}{lcr}
\displaystyle\tilde{\mathfrak{R}}_{u,i}=\sum_{k}(\psi_{u,k}\frac{\p\ \mathfrak{R}_{u,k}}{\p\ \varphi_{u,i}}+\psi_{v,k}\frac{\p\ \mathfrak{R}_{v,k}}{\p\ \varphi_{u,i}}+\psi_{w,k}\frac{\p\ \mathfrak{R}_{w,k}}{\p\ \varphi_{u,i}})+\frac{\p\ {\cal I}}{\p\ \varphi_{u,i}}\\
\displaystyle\tilde{\mathfrak{R}}_{v,i}=\sum_{k}(\psi_{v,k}\frac{\p\ \mathfrak{R}_{v,k}}{\p\ \varphi_{v,i}}+\psi_{u,k}\frac{\p\ \mathfrak{R}_{u,k}}{\p\ \varphi_{v,i}}+\psi_{w,k}\frac{\p\ \mathfrak{R}_{w,k}}{\p\ \varphi_{v,i}})+\frac{\p\ {\cal I}}{\p\ \varphi_{v,i}},\\
\displaystyle\tilde{\mathfrak{R}}_{w,i}=\sum_{k}(\psi_{w,k}\frac{\p\ \mathfrak{R}_{w,k}}{\p\ \varphi_{w,i}}+\psi_{u,k}\frac{\p\ \mathfrak{R}_{u,k}}{\p\ \varphi_{w,i}}+\psi_{v,k}\frac{\p\ \mathfrak{R}_{v,k}}{\p\ \varphi_{w,i}})+\frac{\p\ {\cal I}}{\p\ \varphi_{w,i}}.
%\hat{\bf \Lambda}^{*}_{\rho}=\hat{\bf \Lambda}_{\rho}^{t} \text{~~~(for $\Theta=u,v$ or w)}.
%
\end{array}
\end{equation}
This problem is similar the continuous LBE in time and space with force \eqref{LIBEG}. The main difference of these systems is their streaming steps. In LBE \eqref{LIBEG}, the distributions $\varphi_{\Theta,i}$ for $\Theta=u,v$ and $w$ in a lattice node streams to the adjacent node along $i$ lattice line with ${\bf e}_{i}$ as forward in time.  Conversely, in adjoint LBE \eqref{CAEL}, the adjoint distributions $\psi_{\Theta,i}$ for $\Theta=u,v$ and $w$ in a lattice node streams to the adjacent node along $i$ lattice line with $-{\bf e}_{i}$ as backward in time. Consequently we can derive, in same way as to obtain  \eqref{eq:LBEG}, the discrete adjoint LBE corresponding to \eqref{CAEL} (when ${\bf x}\in {\cal L}_{h}$ and $t=t_{n}$, n=N,N-1,...,1 (backward in time)). 

In the sequel we denote by $G_{i}^{(h,n)}=(\varphi^{h,n}_{u,i},\varphi^{h,n}_{v,i},\varphi^{h,n}_{w,i})$ (respectively $H_{i}^{(h,n)}=(\psi^{h,n}_{u,i},\psi^{h,n}_{v,i},\psi^{h,n}_{w,i})$) the solution of the discrete LBE \eqref{eq:LBEG} (respectively of the discrete adjoint LBE  corresponding to \eqref{CAEL}), ${\cal G}^{h,n}=(G_{i}^{(h,n)})_{i=0,q-1}$, ${\cal H}^{h,n}=(H_{i}^{(h,n)})_{i=0,q-1}$ and ${\bf u}^{h,n}=(u^{h,n},v^{h,n},w^{h,n})=\sum_{i}G_{i}^{(h,n)}$. 

In the objective functional ${\cal J}$, the integrals with respect to time can be approximated by the composite trapezoidal rule, where we adopt the point $t_{n}$ from the time discretization of the state and adjoint systems, as follows (for ${\bf f}^{h}\in {\cal K}^{h}_{ad}$)
\begin{equation}\label{EJD-H}
\begin{array}{lcr}
\displaystyle {\cal J}_{h}({\bf f}^{h})=\int_{\Omega}{\cal I}_{0}({\cal G}^{h,N})dx+
\tau\sum_{n=0}^{N}\!\!\theta_{n}\int_{\Omega}{\cal I}({\cal G}^{h,n},{\bf f}^{h,n})dx
\end{array}
\end{equation} 
where finite dimensional subspace ${\cal K}^{h}_{ad}$ is the internal approximation of ${\cal K}_{ad}$. Here $\theta_{0}=\theta_{N}={1}/{2}$ and $\theta_{n}=1$ for $n\neq 0$ and $n\neq N$.
The gradient of ${\cal J}_{h}$ at point ${\bf f}^{h}=({\bf f}_{1}^{h},{\bf f}_{2}^{h}, {\bf f}_{3}^{h})$ in direction of ${\bf g}^{h}=({\bf g}_{1}^{h},{\bf g}_{2}^{h}, {\bf g}_{3}^{h})$ is given by (analogous to the one of ${\cal J}$)
\begin{equation}\label{JD-H-B}
\begin{array}{lcr}
\displaystyle {\cal J}'_{h}({\bf f}^{h}).{\bf g}^{h}=\tau\sum_{n=0}^{N}\!\!\theta_{n}\int_{\Omega}\left[
\left(\sum_{i}\psi_{u,i}\frac{\p\ {\mathfrak{F}}_{u,i}}{\p\ {\bf f}_{1}}({\cal G}^{h,n},{\bf f}^{h,n})+\frac{\p\ {\cal I}}{\p\ {\bf f}_{1}}({\cal G}^{h,n},{\bf f}^{h,n})\right){\bf g}_{1}^{h,n}\right. \\
\displaystyle \hspace{3cm}+\left(\sum_{i}\psi_{v,i}\frac{\p\ {\mathfrak{F}}_{v,i}}{\p\ {\bf f}_{2}}({\cal G}^{h,n},{\bf f}^{h,n})+\frac{\p\ {\cal I}}{\p\ {\bf f}_{2}}({\cal G}^{h,n},{\bf f}^{h,n})\right){\bf g}_{2}^{h,n}\\
\displaystyle \left.\hspace{3cm}+\left(\sum_{i}\psi_{w,i}\frac{\p\ {\mathfrak{F}}_{w,i}}{\p\ {\bf f}_{3}}({\cal G}^{h,n},{\bf f}^{h,n})+\frac{\p\ {\cal I}}{\p\ {\bf f}_{3}}({\cal G}^{h,n},{\bf f}^{h,n})\right){\bf g}_{3}^{h,n}\right]dx.
\end{array}
\end{equation}
According to the previous discrete formula \eqref{EJD-H} and \eqref{JD-H-B} we can now present the optimization process using the adjoint LBE  in order to solve the lattice Boltzmann discrete optimal control problem 
(in ${\cal K}^{h}_{ad}$). 

For $k=0,...$  (iteration index) we denote by $X_{k}$ the numerical approximation of all functions $X$ (involved in algorithm) at $k$th iteration of algorithm.
\begin{enumerate}
\item Initialization : ${\bf f}^{h}_{k}$ (given)
\item Compute primal problem \eqref{eq:ReacDiff} with source term ${\bf f}^{h}_{k}$ by using LBE \eqref{eq:LBEG} (with initial conditions), gives
the distribution function ${\cal G}^{h,n}_{k}$ and primal solution ${\bf u}_{k}^{h,n}$, for $n=0,..., N$ (forward in time).\label{OPSTEP2}
\item Compute discrete adjoint LBE corresponding to \eqref{CAEL} (with final conditions), based on ${\cal G}^{h}_{k}$, ${\bf u}_{k}^{h}$ and  ${\bf f}^{h}_{k}$, gives 
the adjoint distribution function  ${\cal H}_{k}^{h,n}$, for $n=N,..., 0$ (backward in time). 
\item Compute ${\cal J}'_{h}({\bf f}_{k}^{h})$ from \eqref{JD-H-B} (based on ${\cal H}_{k}^{h}$, ${\cal G}^{h}_{k}$, ${\bf u}_{k}^{h}$ and  ${\bf f}^{h}_{k}$).
\item Perform an appropriate optimization algorithm (Gradient, Conjugate Gradient,...), gives the new control variable  ${\bf f}^{h}_{k+1}$.
\item If gradient of ${\cal J}_{h}$ is sufficiently small (convergence) then end; else set k:=k+1 and goto step \ref{OPSTEP2}. Discrete approximation of optimal solution $({\bf u}^{*}, {\bf f}^{*})$ is  
$({\bf u}_{k}^{h},{\bf f}^{h}_{k})$.
\end{enumerate}
\subsubsection{Continuous adjoint problem and multiple-relaxation-time LBM}
The treated control problem, can be solved also numerically from the discretization of the continuous adjoint problem \eqref{FVM-PF-Adj-K} by 
using the proposed multiple-relaxation-time lattice Boltzmann method. Indeed, the adjoint problem \eqref{FVM-PF-Adj-K} (with final conditions), which is backward in time, can be transformed  into an initial-boundary value problem by reversing the sense of time i.e., $t :=T-t$ and the derived system is similar to \eqref{eq:ReacDiff} (with initial conditions). The algorithm can be formulated as follows.

For $k=0,...$  (the iteration index) we denote by $X_{k}$ the numerical approximation of all functions $X$ (involved in the algorithm) at $k$th iteration of the algorithm.
\begin{enumerate}
\item Initialization : ${\bf f}^{h}_{k}$ (given)
\item Compute primal problem \eqref{eq:ReacDiff} with the source term ${\bf f}^{h}_{k}$ by using LBE \eqref{eq:LBEG} (with initial conditions), gives
the primal solution ${\bf u}_{k}^{h,n}$, for $n=0,..., N$.\label{OLPSTEP2}
\item Compute adjoint \eqref{FVM-PF-Adj}, based on ${\bf u}_{k}^{h}$ and  ${\bf f}^{h}_{k}$, gives adjoint solution $\tilde{\bf u}_{k}^{h,n}$, for $n=N,..., 0$ (by reversing the sense of time and using the LBE \eqref{eq:LBEG}). 
\item Compute ${\cal J}'_{h}({\bf f}_{k}^{h})$ from \eqref{GJA} (based on $\tilde{\bf u}_{k}^{h}$, ${\bf u}_{k}^{h}$ and  ${\bf f}^{h}_{k}$).
\item Perform an appropriate optimization algorithm (Gradient, Conjugate Gradient,...), gives the new control variable  ${\bf f}^{h}_{k+1}$.
\item If gradient of ${\cal J}_{h}$ is sufficiently small (convergence) then end; else set k:=k+1 and goto step \ref{OLPSTEP2}. Discrete approximation of optimal solution $({\bf u}^{*}, {\bf f}^{*})$ is  
$({\bf u}_{k}^{h},{\bf f}^{h}_{k})$.
\end{enumerate}
%%%%%%%%%%%%%%%%%%%%%%%%%%%%%%%%%%%%%%%%%%%%%%%%%%%%%%%%%%%%%%%%%%%%%%%%%%%%%%%%%%%%%%%%%%%%%%%%%%%%%%%%%%%%%%%%%%%%%%%%%%%%%%%%%%%%%%%%%%%%%%%%%
\section{Conclusions}
This paper presents a complete formulation for the continuous and discrete adjoint approach for automatic control of some complex real-life systems, which appear in diverse biochemical, biological and biosocial criminology problems, governed by a class of nonlinear coupled anisotropic convection-diffusion chemotaxis-type system (NCACDCS). The optimal control problem subject to primal system NCACDCS has been formulated within a general framework. First-order necessary optimality conditions are established by using sensitivity and adjoint calculus. These conditions allows easy access to gradients of objective functions (even with a large number of estimate parameters and/or constraints), that is very useful in numerical optimization algorithms. A general adjoint-based sensitivity methodology for these inverse optimization problems using a general multiple-relaxation-time Lattice Boltzmann (MRT) is developed and an adjoint multiple-relaxation-time lattice Boltzmann model (AMRT), which is found to be as simple as MRT model with also highly-efficient parallel nature, is derived. 

Despite its apparent complexity, the developed discrete algorithm is quite easy to implement with low implementation cost (even in the case of large time horizons $T$ or fine space-time meshes).
Due to the local nature of the adjoint and primal collision operators and the fact that nearest neighbor nodes only interact during the streaming stages (allowing for straightforward and efficient parallelization and use of GPU computing), this proposed method is a reliable and efficient approach to solve this type of practical real-life problems (without forgetting that the ability of LBM to approximate complicated geometries with simple algorithmic modifications). Moreover, to help reduce further the computational complexity, CPU and memory costs, we can coupled the developed method and Proper Orthogonal Decomposition Reduced Order Model (POD).

\section{Conflict of interest}
Author declared that there are no conflicts of interests


\addcontentsline{toc}{section}{Bibliographie}
\begin{thebibliography}{9}
\renewcommand{\baselinestretch}{1.53}
%%
%%
\bibitem{A1} \sc {R.A. Adams}, {\it Sobolev spaces,} Academic Press, New-York, 1975.
%
\bibitem{ABP} \sc {D. Ambrosi, F. Bussolino, L. Preziosi}, {\it A review of vasculogenesis models,} Computational and Mathematical Methods in Medicine, 2005;6:1-19.
%
\bibitem{AC} \sc{A.R.A. Anderson, M.A.J. Chaplain}, {\it Continuous and discret mathematical models of tumor-induced angiogenesis,} Bull. Math. Biol., 1998;60:857-899.
%
\bibitem{ABS} \sc{ B. Andreianov, M. Bendahmane, M. Saad,} {\it Finite volume methods for degenerate chemotaxis model,} J. Comput. Appl. Math., 2011;235:4015-4031.
%
\bibitem{AM} \sc {R.P. Araujo, D.L.S. McElwain}, {\it A history of the study of solid tumour growth: the contribution of mathematical modelling,} Bull. Math. Biol., 2004;66:1039-1091.
%
%
\bibitem{ASV} \sc{E.E.E. Arenas, A. Stevens, J.J.L. Vel'azquez,} {\it Simultaneous finite time blow-up in a two-species model for chemotaxis,} Analysis, 2009;29:317-338.
%
\bibitem{BIJB} \sc{A. Belmiloudi}, {\it Mathematical modeling and optimal control problems in brain tumor targeted drug delivery strategies,} International Journal of Biomathematics, 
2017;10:1750056[62 pages].
%
\bibitem{BDI} \sc{A. Belmiloudi}, {\it Dynamical behavior of nonlinear impulsive abstract partial differential equations on networks with multiple time-varying delays and mixed boundary conditions involving time-varying delays,} J. Dynam. Control Systems, 2015;21:95-146.
%
\bibitem{B1} \sc{A. Belmiloudi}, {\it Robust control problem of uncertain bidomain models in cardiac electrophisiology,} Journal of Coupled Systems and Multiscale Dynamics, 2013;19:332-350.
%
\bibitem{B2} \sc{A. Belmiloudi}, {\it Stabilization, Regulation and Robust Control of Uncertain processes and parameters in Porous Medium systems with Applications,} In Edited Book: Focus on Porous Media Research, C. Zhao ed., Nova Science Publishers, New York, 2013, pp. 165-228. 
%
\bibitem{BTT} \sc{A. Belmiloudi,} {\it  Thermal Therapy: Stabilization and Identification,} In Heat Transfer Mathematical Modelling, Numerical Methods and Information Technology, A. Belmiloudi ed., INTECH: Vienna, 2011, pp. 33-76.
%
\bibitem{B8} \sc{A. Belmiloudi},  {\it Parameter identification problems and analysis of the impact of porous media in biofluid heat transfer in biological tissues during thermal therapy.} Nonlinear Analysis: Real World Applications, 2010;11:1345-1363.
%
\bibitem{B0} \sc{A. Belmiloudi}, {\it Stabilization, optimal and robust control. Theory and applications in biological and physical sciences,} Communications and Control Engineering. Springer, London, 2008.  
%
\bibitem{Ba} \sc{M. Ben Amar,} {\it Chemotaxis migration and morphogenesis of living colonies,} Eur. Phys. J. E. Soft. Matt., 2013;36(68):1-13.
%
\bibitem{BL1} \sc{E. Ben-Jacob, I. Cohen, H. Levine,} {\it Cooperative self-organization of microorganisms,} Adv. Phys., 2000;49:395-554.
%
\bibitem{BL2} \sc{E. Ben-Jacob, H. Levine,} {\it Self-engineering capabilities of bacteria,}  J. R. Soc. Interface, 2006;3:197-214.
%
\bibitem{BBI} \sc{L. Berezansky, E. Braverman,  L. Idels,} {\it Effect of treatment on the global dynamics of delayed pathological angiogenesis models,} J. Theor. Biol., 2014;363:13-21. 
%
\bibitem{BGK} \sc{P.L. Bhatnagar, E.P. Gross, M. Krook,} {\it A model for Collision Processes in Gases Small Amplitude Processes in Charged and Neutral One-Component Systems}, Phys. Rev., 1954;94:511-525.
 %
\bibitem{BKZ}  \sc{P. Biler, G. Karch, J. Zienkiewicz,} {\it  Optimal criteria for blowup of radial and N-symmetric solutions of chemotaxis systems,} Nonlinearity,2015;28:43-69. 
%
 \bibitem{Boltz1872} \sc{L. Boltzmann}, {\it Weitere Studien uber das Warme gleichgenicht unfer Gasmolakuler}, Sitzungsberichte Akad. Wiss., 1872;66:275.
%
\bibitem{BFD}  \sc{M. Burger, M. Di Francesco, Y. Dolak-Struss,} {\it The Keller–Segel model for chemotaxis with prevention of overcrowding: linear vs. nonlinear diffusion,} SIAM J. Math. Anal., 
2008;38:1288-1315.
%
\bibitem{CL}  \sc{M.A.J. Chaplain, G. Lolas,} {\it  Mathematical Modeling of Cancer of Tissue: Dynamic Heterogeneity,} Networks and Heterogeneous Media, 2006;1:399-439.
%
\bibitem{CO} \sc{J.O. Campos, R.S. Oliveira, R.W. dos Santos, M. Rocha,} {\it Lattice Boltzmann method for parallel simulations of cardiac electrophysiology using GPUs,} J. Comput. Appl. Math., 2016;295:70-82.
%
\bibitem{ChapEnsk} \sc{S. Chapman S., T.G. Cowling,} {\it The mathematical theory of non-uniform gases}, Cambridge University Press, 1990.
%
\bibitem{CB} \sc{ S. Corre, A. Belmiloudi,} {\it Coupled Lattice Boltzmann Modeling of Bidomain type models in Cardiac Electrophysiology,} In Mathematical and Computational Approaches 
in Advancing Modern Science and Engineering, eds. J. B\'elair, I. Frigaard, H. Kunze, R., Melnik, J. Spiteri,  Springer-Verlag, Switzerland, 2016, pp. 209-221.
%
\bibitem{CBA} \sc{S. Corre, A. Belmiloudi,} {\it 	Coupled lattice Boltzmann method for numerical simulations of fully coupled heart and torso bidomain system in electrocardiology,} Journal of Coupled Systems and Multiscale Dynamics, 2016;3:207-229.
%
%
\bibitem{DHG} \sc{D. d'Humi\`eres, I. Ginzburg, M. Krafczyk, P. Lallemand, L.S. Luo,} {\it Multiple relaxation-time lattice Boltzmann models in three-dimensions,} Proc. Roy. Soc.
London A, 2002;360:437-451.
 %
\bibitem{EJDW} \sc{H.J. Eberl, E.M. Jalbert, A. Dumitrache, G.M. Wolfaardt,} {\it A spatially explicit model of inverse colony formation of cellulolytic biofilms,} Biochemical Engineering Journal, 
2017;122:141-151. 
%
\bibitem{Ef}  \sc{M. Efendiev,} {\it Attractors for Degenerate Parabolic Type Equations,} American Mathematical Society, Providence, RI, 2013. 
%
\bibitem{EZ}  \sc{M. Efendiev, A. Zhigun,} {\it  On an exponential attractor for a class of PDEs with degenerate diffusion and chemotaxis,} Preprint.
%
\bibitem{Eal} \sc{M. Eisenbach (Ed.),} {\it Chemotaxis,} Imperial College Press, London, 2004.
%
\bibitem{Ep} \sc{Y. Epshteyn,} {\it Upwind-difference potentials method for Patlak–Keller–Segel chemotaxis model,} J. Sci. Comput., 2012;53,689-713.
%
\bibitem{EK} \sc{Y. Epshteyn, A. Kurganov,} {\it New interior penalty discontinuous Galerkin methods for the Keller–Segel chemotaxis model,} SIAM J. Numer. Anal., 2008;47:386-408.
%
\bibitem{Fi} \sc{F. Filbet,} {\it A finite volume scheme for the Patlak–Keller–Segel chemotaxis model,} Numer. Math., 2006;104:457-488.
%
\bibitem{FS}  \sc{K. Fujie, T. Senba,} {\it Global existence and boundedness of radial solutions to a two dimensional fully parabolic chemotaxis system with general sensitivity,}  Nonlinearity, 
2016;29:2417. 
%
\bibitem{GC} \sc{A. Gerisch, M.A. Chaplain,} {\it  Mathematical modelling of cancer cell invasion of tissue: local and non-local models and the effect of adhesion,} J. Theor.Biol., 2008;250:684-704.
%
%
\bibitem{HM}  \sc{M. Hekmat, M. Mirzaei}, {\it Extraction of macroscopic and microscopic adjoint concepts using a lattice Boltzmann method and discrete adjoint approach,} Phys. Rev. E., 2015;91:23-42.
%
\bibitem{HP1} \sc{T. Hillen, K. J. Painter,}  {\it  A user's guidence to PDE models for chemotaxis,} J. Math. Biol., 2009;58:183-217.
%
\bibitem{Ho2} \sc{D. Horstmann,} {\it Generalizing the Keller-Segel Model: Lyapunov Functionals, Steady State Analysis, and Blow-Up Results for Multi-species Chemotaxis Models in the Presence of Attraction and Repulsion Between Competitive Interacting Species,} J. Nonlinear Sci., 2011;21:231-270 .
%
\bibitem{Ho3}  \sc{D. Horstmann,}  {\it From 1970 until present: the Keller-Segel model in chemotaxis and its consequences: I,} Jahresber. Deutsch. Math.-Verein., 2003;105:103-165.
\bibitem{Ho4}  \sc{D. Horstmann,} {\it From 1970 until present: the Keller-Segel model in chemotaxis and its consequences: II,} Jahresber. Deutsch. Math.-Verein., 2004;106:51-69.
%
\bibitem{HW}  \sc{R. Huang, H. Wu,} {\it A modified multiple-relaxation-time lattice Boltzmann model for convection-diffusion equation,} J. Comput. Physics, 2014;274:50-63. 
%
\bibitem{I1} \sc{R. Insall,} {\it The interaction between pseudopods and extracellular signalling during chemotaxis and directed migration,} Curr. Opin. Cell Biol., 2013;25:526-31.
%
\bibitem{KK} \sc{T.V. Kasyap, D.L. Koch,} {\it Chemotaxis driven Instability of a confined bacterial Suspension,} Phys. Rev. Lett., 2012;108:03101.
%
\bibitem{KS1} \sc{E.F. Keller, L.A. Segel,}  {\it  Initiation of slime mold aggregation viewed as an instability,} J. Theor. Biol., 1970;26:399-415
\bibitem{KS2} \sc{E.F. Keller, L.A. Segel,} {\it  Travelling bands of chemotactic bacteria: a theoretical analysis,} J. Theor. Biol., 1971;30:235-248
%
\bibitem{KTH}  \sc{M. J. Krause,	G. ThaTer , V. Heuveline,} {\it Adjoint-based fluid flow control and optimisation with lattice Boltzmann methods,} Comput. Math. Appl., 2013;65:945-960.
%
\bibitem{LL} \sc{P. Lallemand, L.S. Luo,} {\it Theory of the lattice Boltzmann method: dispersion, dissipation, isotropy, Galilean invariance and stability,} Phys. Rev. E., 2000;61:6546-6562.
%
\bibitem{LR} \sc{L. Laniewski-Wollk, J. Rokicki,} {\it Adjoint Lattice Boltzmann for topology optimization on multi-GPU architecture,} Comput. Math. Appl., 2016;71:833-848.
%
\bibitem{LMa} \sc{D. Lebiedz D, H. Maurer,} {\it External optimal control of self-organisation dynamics in a chemotaxis reaction diffusion system,} Syst. Biol. (Stevenage), 2004;1:222-229.
%
\bibitem{LLFS} \sc{C. Lin, K.H. Luo, L. Fei, S. Succi,}  {\it A multi-component discrete Boltzmann model for nonequilibrium reactive flows,} Scientific reports, 2017;7:1-12.
%
\bibitem{LCEM} \sc{M. Luca, A. Chavez-Ross, L. Edelstein-Keshet, A. Mogilner,} {\it Chemotactic signalling, microglia and Alzheimer's disease senile plaques: Is there a connection?} Bull. Math. Biol., 
2003;65:693-730.
%
\bibitem{LCA} \sc{P.M. Lushnikov, N. Chen, M. Alber,} {\it Macroscopic dynamics of biological cells interacting via chemotaxis and direct contact,} Phys. Rev. E., 2008;78:061904.
%
\bibitem{MM} \sc{P.K. Maini and J. D. Murray,} {\it A nonlinear analysis of a mechanical model for biological pattern formation,} SIAM J. Appl. Math., 1988;48:1064-1072.
%
\bibitem{Mg} \sc{G. Marinoschi,} {\it  A control problem for a cross-diffusion system in a nonhomogeneous medium,} Journal of Biological Dynamics, 2013;7:88-107.
%
\bibitem{MWZZ}  \sc{C. Mu, L. Wang, P. Zheng, Q. Zhang,} {\it  Global existence and boundedness of classical solutions to a parabolic-parabolic
chemotaxis system,} Nonlinear Analysis: Real World Applications, 2013;14:1634-1642.
%
\bibitem{Mu} \sc{J.D. Murray,} {\it Mathematical Biology,} Springer-Verlag, New York, 1993.
%
\bibitem{Mu1} \sc{J.D. Murray,} {\it  Mathematical Biology I: An Introduction,} Springer, Berlin, 3rd edition, 2002.
%
\bibitem{NO}  \sc{E. Nakaguchi, K. Osaki.} {\it  Global solutions and exponential attractors of a parabolic-parabolic system for chemotaxis with subquadratic degradation,} Discrete Contin. Dyn. Syst.,
Ser. B., 2013;18:2627-2646.
%
\bibitem{OS} \sc{H. Othmer, A. Stevens,} {\it Aggregation, blowup and collapse: The abc’s of generalized taxis,} SIAM J. Appl. Math., 1997;57:1044-1081.
%
\bibitem{PS1} \sc{K. Painter, J. A. Sherratt,}  {\it Modelling the movement of interacting cell populations,} J. Theor. Biol., 2003;225:327-339.
%
\bibitem{PH2} \sc{K. Painter, T. Hillen,} {\it Spatio-temporal chaos in a chemotaxis model,} Phys. D, 2011;240:363-375.
%
\bibitem{PX}  \sc{Y. Peng, Z. Xiang,} {\it Global existence and boundedness in a 3D Keller-Segel–Stokes system with nonlinear diffusion and rotational flux,}	
Zeitschrift für angewandte Mathematik und Physik, 2017;68:68.
%
\bibitem{Pe} \sc{ B. Perthame,} {\it Transport equations in biology,} Birkh\"auser, Basel, 2007.
%
%
\bibitem{RST} \sc{P. Roca-Cusachs, R. Sunyer, X. Trepat,}  {\it Mechanical guidance of cell migration: lessons from chemotaxis,} Curr. Opin. Cell Biol., 2013;25:543-549.
%
\bibitem{RB}  \sc{N. Rodriguez, A. Bertozzi,} {\it  Local existence and uniqueness of solutions to a PDE model for criminal behavior,} Math. Mod. Meth. Appl. Sci., 2010;20:1425-1457.
%
\bibitem{Sa}  \sc{N. Saito,} {\it Error analysis of a conservative finite-element approximation for the Keller-Segel system of chemotaxis,} Commun. Pure Appl. Anal., 2012;11:339-364.
%
\bibitem{Sal} \sc{G. Serini, D. Ambrosi, E. Giraudo, A. Gamba,L. Preziosi,F. Bussolino,}  {\it Modeling the early stages of vascular network assembly,} EMBO J., 2003;22:1771-1779.
%
\bibitem{SOPT} \sc{M.B. Short, M.R. D'Orsogna, V.B. Pasour, G.E. Tita, P.J. Brantingham, A.L. Bertozzi, L. Chayes,} {\it A statistical model of criminal behavior,} Math. Models Methods Appl. Sci., 2008;18:1249-1267.
%
\bibitem{SSK} \sc{R. Strehl, A. Sokolov, D. Kuzmin, D. Horstmann, S. Turek,} {\it A positivity-preserving finite element method for chemotaxis problems in 3D,} J. Comput. Appl. Math., 2013;239:290-303.
%
\bibitem{SST} \sc{R. Strehl, A. Sokolov, S. Turek,} {\it Efficient, accurate and flexible finite Element solvers for chemotaxis problems,} Comput. Math. Appl., 2012;64:175-189.
%
\bibitem{SY}  \sc{Y. Sugiyama, Y. Yahagi,} {\it Extinction, decay and blowup for Keller-Segel systems of fast diffusion type,} Journal of Differential Equations, 2011;250:3047-3087.
%
\bibitem{TW}  \sc{Y. Tao, M. Winkler,} {\it Eventual smoothness and stabilization of large-data solutions in a three-dimensional chemotaxis system with consumption of chemoattractant,} Journal of Differential Equations, 2012;252:2520-2543.
%
\bibitem{UYa} \sc{S. Uk Ryu, A. Yagi,} {\it  Optimal control of Keller-Segel equations,} J. Math. Anal. Appl., 2001;256:45-66.
%
\bibitem{VM} \sc{R. Van Houdt C.W. Michiels,} {\it  Role of bacterial cell surface structures in Escherichia coli biofilm formation, Research in Microbiology,} 2005;156:626-633.
%
\bibitem{Wi}  \sc{M. Winkler,} {\it  Finite-time blow-up in the higher-dimensional parabolic-parabolic Keller-Segel system,} Journal de Math\'ematiques Pures et Appliqu\'ees, 2013;100:748-767.
%
\bibitem{Wr1} \sc{D. Wrzosek,} {\it Long time behavior of solutions to a chemotaxis model with volume filling effects,} Proc. R. Soc. Edinburgh A: Math., 2006;136:431-444.

\bibitem{Wr2} \sc{D. Wrzosek,} {\it Model of chemotaxis with threshold density and singular diffusion,} Nonlinear Anal. TMA, 2010;73:338-349
%
\bibitem{XLZ} \sc{ A. Xu, C. Lin, G. Zhang, Y. Li,} {\it Multiple-relaxation-time lattice Boltzmann kinetic model for combustion,} Phys. Rev. E, 2015;91:043306.
%
\bibitem{XO} \sc{C. Xue, H. G. Othmer,} {\it Multiscale models of taxis-driven patterning in bacterial populations,} SIAM J. Appl. Math., 2009;70:133-167.
%
\bibitem{YSC} \sc{X. Yang, B. Shi, Z. Chai,} {\it Coupled lattice Boltzmann method for generalized Keller-Segel chemotaxis model}, Comput.  Math. Appl., 2014;12:1653-1670.
%
\bibitem{ImprovedBounceBack} \sc{X. Yin, J. Zhang,} {\it An improved bounce-back scheme for complex boundary conditions in lattice Boltzmann method}, Journal of Computational Physics, 
2012;231:4295-4303.
%
\bibitem{YY}  \sc{K. Yajia, T. Yamadaa, M. Yoshino, T. Matsumoto, K. Izui, S. Nishiwaki,} {\it Topology optimization in thermal-fluid flow using the lattice Boltzmann method,}  J. Comput. Phys., 2016;307:355-377.
%
\bibitem{YN} \sc{H. Yoshida, M. Nagaoka,} {\it  Multiple-relaxation-time lattice Boltzmann model for the convection and anisotropic diffusion equation,} J. Comput. Phys., 2010;229:7774-7795.
%
\bibitem{ZhangBounceBack12} \sc{Zhang, B. Shi, Z. Guo, Z. Chai, J. Lu,} {\it General bounce-back scheme for concentration boundary condition in the lattice-Boltzmann method}, American Physical Society, 2012;85:016701.
%
\bibitem{ZZL} \sc{R. Zhang, J. Zhu, A. Loula, X. Yu,} {\it Operator splitting combined with positivity-preserving discontinuous Galerkin method for the chemotaxis model,} J. Comput. Appl. Math., 2016;302:312-326.
%
\bibitem{GuoBoundary2002} \sc{G. Zhao-Li, Z. Chu-Guang, S. Bao-Chang,} {\it Non-equilibrium extrapolation method for velocity and pressure boundary conditions in the lattice Boltzmann method}, Chinese Physics, 2002;4:366-374.
%
\bibitem{ZBZ}  \sc{C. Zhenhua, S. Baochang, G. Zhaoli,} {\it A Multiple-Relaxation-Time Lattice Boltzmann Model for General Nonlinear Anisotropic Convection-Diffusion Equations,}, J. Sci. Comput., 2016;69:355-390.
%
\bibitem{ZC}  \sc{Y. Zhu, F. Cong} {\it Global Existence to an Attraction-Repulsion Chemotaxis Model with Fast Diffusion and Nonlinear Source,} Discrete Dynamics in Nature and Society, 2015:143718 [8 pages].


\end{thebibliography}
\end{document}